\documentclass[12pt]{article}
\usepackage{amsfonts}
\usepackage{enumerate}
\usepackage{CJK}
\usepackage{graphicx,psfrag,color}
\usepackage{fancyhdr}
\usepackage{amsmath,amsthm,amssymb}
\usepackage{mathrsfs}
\usepackage{indentfirst}
\usepackage{epsfig, epsfig, graphics, lscape, amsbsy, amstext, natbib}
\newtheorem{theorem}{THEOREM}
\newtheorem{corollary}{COROLLARY}

\newtheorem{lemma}{LEMMA}
 \newtheorem{assumption}{ASSUMPTION}
\newcommand{\T}{\!\top\!}

\begin{document}
\begin{titlepage}

\title{Exponentially tilted likelihood inference on growing dimensional unconditional moment models 
{\small \author{Nian-Sheng Tang\footnote{Correspondence to:
Dr. Nian-Sheng Tang, Key Lab of Statistical Modeling and Data Analysis of Yunnan Province, Yunnan University,
Kunming 650091, P. R. of China. Tel: 86-871-5032416~ Fax: 86-871-5033700~ E-mail: nstang@ynu.edu.cn}  ~~Xiao-Dong Yan~~and~~ Pu-Ying Zhao\\
{\small {\small {\it Key Lab of Statistical Modeling and Data Analysis of Yunnan Province,}}}\\[-3mm]
{\small {\small {\it Yunnan University, Kunming 650091, China} }} }}}
\date{}
\maketitle

\thispagestyle{empty}

\noindent {\bf Abstract}:
Growing-dimensional data with likelihood unavailable are often encountered in various fields. This paper presents a penalized exponentially tilted likelihood (PETL) for variable selection and parameter estimation for growing dimensional unconditional moment models
in the presence of correlation among variables and model misspecification. Under some regularity conditions, we investigate the consistent and oracle properties of the PETL estimators of parameters,
and show that the constrainedly PETL ratio statistic for testing contrast hypothesis asymptotically follows the central chi-squared distribution.
Theoretical results reveal that the PETL approach is robust to model misspecification. We also study high-order asymptotic properties of the proposed PETL estimators.
Simulation studies are conducted to investigate the finite performance of the proposed methodologies. An example from the Boston Housing Study is illustrated.\\

\begin{CJK*}{GBK}{song}

\end{CJK*}
\end{titlepage}

\newpage
\section{Introduction}

\noindent Exponentially tilted (ET) likelihood (Imbens, Spady and Johnson, 1998) is a useful nonparametric approach to evaluate estimates and confidence regions of unknown parameters in
unconditional moment models  of the form $E\{g(x;\theta)\}=0$, which provides a unified approach for parameter estimation in a class of statistical models with likelihood function unavailable, where $g(x;\theta)$ is a vector-valued nonlinear function of a random vector $x$ and a parameter vector $\theta$.
The merits of the ET likelihood include (i) it behaves better than empirical likelihood under model misspecification (Schennach, 2007), that is, the ET likelihood is robust to model misspecification, (ii) it allows a computationally convenient treatment of misspecified models (Kitamura, 2000), and (iii) it is flexible in incorporating auxiliary information. Hence, several authors, for example, Schennach (2005, 2007), Zhu et al. (2009) and Caner (2010), discussed its properties and applications when the number of parameters is fixed and less than or equal to sample size.

Growing-dimensional parametric or semiparametric models are widely used to make statistical inference on complicated data sets such as longitudinal and panel data in econometrics (Fan and Peng, 2004). It is commonly assumed that only a small number of covariates actually contribute to the considered models, which leads to the well-known sparse models for helping interpretation and improving prediction accuracy (Bradic, Fan and Wang, 2011). To this end, many penalized methods have been developed to simultaneously select the important covariates and estimate parameters in various statistical models when the number of parameters diverges. For example, Fan and Peng (2004) investigated the nonconcave penalized likelihood with a growing number of nuisance parameters in a linear regression model; Lam and Fan (2008) presented a profile-kernel likelihood inference in a linear regression model; Wang, Li and Leng (2009) studied shrinkage tuning parameter selection; Zou and Zhang (2009) proposed an adaptive elastic-net procedure for a linear regression model; Li, Peng and Zhu (2011) investigated asymptotic properties of a nonconcave penalized M-estimator in a sparse, diverging-dimensional, linear regression model; Caner and Zhang (2014) extended the least squares based adaptive elastic net estimator of Zou and Zhang (2009) to generalized method of moments (GMMs); Caner, Han and Lee (2016) presented an adaptive elastic net GMM estimation for many invalid moment conditions. Recently, Leng and Tang (2012) presented a penalized empirical likelihood method in estimating equations, which can be used to improve the efficiency of parameter estimation by incorporating some auxiliary information when likelihood function is unavailable, with a diverging number of parameters, but their empirical likelihood method is sensitive to model misspecification. Also, to the best of our knowledge, there is little work done on extending the above mentioned approaches to unconditional moment models with a diverging number of parameters in the presence of model misspecification. More importantly, this extension is challenging in the presence of model misspecification and high correlation among variables because (i) the number of Lagrange multipliers used to obtain the solution to minimizing the ET likelihood function increases with sample size, (ii) the nonconvex optimization is involved (Leng and Tang, 2012), and (iii) there is a well-known ill-posed problem, i.e., the resulting estimator has very slow rate of convergence (see, e.g., Ai and Chen, 2003; Hall and Horowitz, 2005; Darolles, Fan, Florens and Renault, 2011; Chen and Pouzo, 2012).

In this paper, we develop a penalized ET (PET) likelihood procedure for parameter estimation, variable selection and statistical inference for unconditional moment models
with a diverging number of parameters in the presence of model misspecification and high correlation among variables via the sieve method (Ai and Chen, 2003).
With a proper penalty function and diverging rate of dimensionality,
we demonstrate that (i) the resulting estimator possesses the advantages of the penalized likelihood approach, i.e.,
the PET method has the oracle properties (Fan and Li, 2001) that it identifies the true sparse structure of the considered model with probability tending to one and
with the optimal efficiency; (ii) the resulting estimator has the advantages of the ET likelihood method, i.e., the PET method behaves better than the penalized empirical likelihood approach
in the presence of model misspecification;
(iii) the constrainedly profiled PET likelihood ratio statistic is asymptotically distributed as the chi-squared distribution indicating that the Wilks' theorem holds,
which can be used to test hypotheses and construct confidence regions of parameters of interest. In addition, we extend the high-order asymptotic properties of the ET estimator given in Schennach (2007) for a fixed number of parameters to the case that the number of parameters diverges; and we also establish selection consistency for NP dimension case.

The rest of this paper is organized as follows. Section 2 first introduces the PET likelihood, and then investigates the oracle properties of the proposed PET estimators,
asymptotic chi-squared distribution, high-order asymptotic properties and selection consistency.
Simulation studies are given in Section 3. An example from the Boston Housing Study is analyzed in Section 4. Some concluding remarks are given in Section 5.
Proofs of Theorems are presented in Appendix.

\section{Methods}
\renewcommand{\theequation}{2.\arabic{equation}}
\setcounter{equation}{0}
\subsection{Exponentially tilted likelihood}

\noindent Suppose that $X_1,\ldots,X_n$ are independent and identically distributed (i.i.d.) random vectors from an unknown distribution $F(x)$ with $x\in \mathcal{X}\subset \mathcal{R}^{\iota}$. Without assuming a specific form of $F(x)$,
we are interested in making inference on a $p\times 1$ vector of unknown parameters of interest, denoted by $\theta$, based on $r$ ($r\geq p$) functionally independent estimating functions   $g(X_i;\theta)=(g_1(X_i;\theta),\ldots,g_r(X_i;\theta))^{\T}$ that satisfy the unconditional moment condition of the form: $E_{F_x}\{g(X_i;\theta_0)\}=0$ for $\theta_0\in\Theta\subset\mathcal{R}^p$ and $i=1,\ldots,n$, which is usually referred to as general estimating equations or unconditional moment models (Owen, 2001), where $\theta_0$ is the unique true value of $\theta$ and $E_{F_x}$ denotes the expectation taken with respect to $F(x)$. The selection of $g(X;\theta)$ is flexible and the details can refer to Leng and Tang (2012).

When $r=p$, one can obtain estimation of $\theta$ by solving the following unconditional moment conditions: $n^{-1}\sum_{i=1}^ng(X_i;\theta)=0$ (Leng and Tang, 2012).
When $r>p$ and $p$ is fixed,
one can employ empirical likelihood approach to obtain more efficient estimation of $\theta$ by combining available information (Qin and Lawless, 1994).
However, when $r>p$ and $p$ is large, it is commonly assumed that only a small number of variables actually contribute to unconditional moment conditions,
which leads to the sparsity pattern in unknown parameter vector $\theta$ and thus makes variable selection crucial (Bradic, Fan and Wang, 2011).
To this end, Leng and Tang (2012) studied growing dimensional unconditional moment models via empirical likelihood approach, and presented a
penalized empirical likelihood procedure for parameter estimation and variable selection. In what follows, we present an ET approach to
investigate parameter estimation and variable selection for unconditional moment models with a growing number of parameters because the
ET likelihood is a robust nonparametric tool to make statistical inference on unconditional moment
models (Imbens et al., 1998; Owen, 2001) when unconditional moment models are misspecified.

For $i=1,\ldots,n$, let $w_i=dF(X_i)={\rm Pr}(\mathbb{X}_i=X_i)$, where $X_i$ is the observation of random vector $\mathbb{X}_i$. The ET likelihood can be defined as the Kullback-Leibler divergence between the empirical frequencies $1/n$ and $w_i$ subject to some restrictions. Following Imbens et al. (1998), the ET estimator $\hat\theta_{ET}$ of $\theta$ is the solution to the following Kullback-Leibler information criterion: ${\rm inf}_{w_1,\ldots,w_n,\theta}\sum_{i=1}^nw_i\log w_i$ subject to $\sum_{i=1}^nw_i=1$, $w_i\geq 0$ and $\sum_{i=1}^nw_ig(X_i;\theta)=0$. To wit, the ET likelihood for $\theta$ can be defined as
\begin{equation}\label{PET1}
L(\theta)={\rm inf}\left\{\prod_{i=1}^nw_i\log(w_i): w_i\geq 0, \sum_{i=1}^nw_i=1, \sum_{i=1}^nw_ig(X_i;\theta)=0\right\},
\end{equation}
which is minimized at $w_i=\exp\{\nu^{\T}g(X_i;\theta)\}/\sum_{j=1}^n\exp\{\nu^{\T}g(X_j;\theta)\}$, where $\nu$
is the lagrange multiplier. Under some regularity conditions, it is easily shown that the profiled log-ET likelihood ratio function can be expressed as
\begin{equation}\label{PETL22}
\ell(\nu,\theta)=-\{\log L(\theta)+\log(n)\}=\log\left[\frac{1}{n}\sum\limits_{i=1}^{n}\exp\{\nu^{\T}g(X_i;\theta)\}\right],
\end{equation}
where $\nu$ satisfies $Q_{n1}(\nu,\theta)=n^{-1}\sum_{i=1}^n\exp\{\nu^{\T}g(X_i;\theta)\}g(X_i;\theta)=0$. Thus, $\hat\theta_{ET}$
is the solution to the following nonlinear optimization problem: $\hat\theta_{ET}=\arg\max_{\theta\in\Theta}\inf_{\nu\in\widehat{\mathcal{V}}_n(\theta)} \ell(\nu,\theta)$, where $\widehat{\mathcal{V}}_n(\theta)=\{\nu: \nu^{\T}g(X_i;\theta)\in\mathcal{E},i=1,\cdots,n\}$ in which $\mathcal{E}$ is an open interval containing zero. Generally, under some regularity conditions, $\hat\theta_{ET}$ can be obtained by simultaneously solving the following equations:
\begin{equation} \label{ITEQ}
Q_{n1}(\nu,\theta)=0~~{\rm and}~~Q_{n2}(\nu,\theta)=n^{-1}\sum\limits_{i=1}^n\exp\{\nu^{\T}g(X_i;\theta)\}\nu^{\T}\partial_\theta g(X_i;\theta)=0,
\end{equation}
where $\partial_\theta$ represents the partial derivative with respect to $\theta$.

\subsection{Penalized exponentially tilted likelihood}

\noindent To identify the important covariates in growing-dimensional data analysis,
following Fan and Li (2001), we consider the following profiled PET likelihood ratio function by combining Equation (\ref{PETL22}) and some proper penalty function:
\begin{equation}\label{PETL13}
\ell_p(\theta)=\log\left[\frac{1}{n}\sum_{i=1}^{n}\exp\{\nu^{\T}(\theta)g(X_i;\theta)\}\right]-\sum_{j=1}^{p}p_{\gamma}(|\theta_j|),
\end{equation}
where $\nu(\theta)=\inf_{\nu\in\widehat{\mathcal{V}}_n(\theta)} \ell(\nu,\theta)$, $p_{\gamma}(t)$ is some proper penalty function with a tuning parameter $\gamma$ controlling the trade-off between bias and model complexity.

When the number of parameters diverges, there is a well-known ill-posed
problem (Chen and Pouzo, 2012) for our considered unconditional moment models, i.e., for any $\mathbb{C}> 0$,
there are sequences $\{\theta^{(k)}\}_{k=1}^{\infty}$ of $\theta$
in $\Theta$ such that $\liminf_{k\rightarrow\infty}||\theta^{(k)}-\theta_0||\ge \mathbb{C}$ but
$\liminf_{k\rightarrow\infty}E\{||g(X;\theta^{(k)})||^2\}=0$. To solve the ill-posed problem,
we incorporate two types of regularization methods (e.g., the regularization by sieves and the regularization by penalization) into
the PET procedure (\ref{PETL13}). The commonly used sieves with sparsity constraints can be expressed as
\begin{equation}\label{sieve}
\Theta_{s(n)}=\Big\{\theta\in\Theta: |\mathbb{J}_n|\leq s(n)\Big\},
\end{equation}
where $\mathbb{J}_n=\{j:\theta_j\neq 0\}$, $|\mathbb{J}_n|$ denotes the cardinality of $\mathbb{J}_n$, and $s(n)$ is some given positive integer
associated with sample size $n$.
The constraint  $|\mathbb{J}_n|\leq s(n)$ reflects the prior sparsity information on $\theta_0\in\Theta$. Following Chen and Pouzo (2012),
the penalty function $p_\gamma(t)$ in (\ref{PETL13}) is typically some convex and/or lower semicompact.
Here $p_\gamma(t)$ is taken to be the smoothly clipped absolute deviation (SCAD) penalty because the SCAD penalty function satisfies three desirable conditions for variable selection, i.e. asymptotic unbiasedness, sparsity and continuity of the estimated parameters (Fan and Li, 2001). The SCAD penalty is a function whose first derivative has the following form
$$p_{\gamma}'(t)=\gamma\left\{I(t\leq \gamma)+\frac{(a\gamma-t)_+}{(a-1)\gamma}I(t>\gamma)\right\} $$
for some $a>2$, where $(s)_+=s$ for $s>0$ and $0$ otherwise. The corresponding PET likelihood ratio function is referred to as the SCAD-PET likelihood ratio function.
Similar to Fan and Li (2001), we take $a=3.7$ in our simulation studies.  Thus, under Equations (\ref{PETL13}) and (\ref{sieve}), the SCAD-PET estimator (denoted by $\hat\theta$) of $\theta$ can be obtained by minimizing $\ell_p(\theta)$, i.e.,
\begin{equation}\label{sievepet}
\hat{\theta}=\min_{\theta\in\Theta_{s(n)}}\tilde{\ell}_p(\theta_{\mathbb{J}_n}),
\end{equation}
where $\tilde{\ell}_p(\theta_{\mathbb{J}_n})=\log\left[\frac{1}{n}\sum_{i=1}^{n}\exp\{\tilde{\nu}^{\T}(\theta_{\mathbb{J}_n})
\tilde{g}(\tilde{X}_i;\theta_{\mathbb{J}_n})\}\right]-\sum_{j\in\mathbb{J}_n}p_{\gamma}(|\theta_j|)$,
and $\tilde{\nu}$, $\tilde{g}(\cdot)$, $\tilde{X}_i$ are their corresponding reduced forms of Lagrange, estimating functions and covariates
under the sieve space $\Theta_{s(n)}$ after ignoring the auxiliary variables excluding $\theta_1$, respectively.
Thus, under the sparsity assumption, we can write $\theta=(\theta_1^{\T},\theta_2^{\T})^{\T}$, where $\theta_1\in \mathcal {R}^q$ and $\theta_2\in \mathcal {R}^{p-q}$
correspond to the nonzero and zero components of $\theta$,
respectively.  Once we have the prior sparsity structure of $\theta$,
we can make statistical inference on $\theta_1$ based on the reduced estimating equations $\psi(Z_i;\theta_1)=(\psi_1(Z_i;\theta_1)$
$,\ldots,\psi_k(Z_i;\theta_1))^{\T}$ satisfying the following unconditional moment restrictions: $E_{F_z}\{\psi(Z_i;\theta_{10})\}=0$
for $\theta_{10}\in \mathcal{R}^q$, where $\psi(Z_i;\theta_{1})$ is some reduced version of $g(X_i;\theta)$ under $\theta_2=0$, $\theta_{10}$ is the true value of $\theta_1$,
and $Z_i$ is the selected important covariates (i.e., $Z_i\subset X_i$) for $i=1,\ldots,n$.
It is assumed that $Z_1,\ldots,Z_n$ are independent and identically distributed as an unknown
distribution $F(z)$ with $z\in \mathcal{Z}\subset \mathcal{R}^q$. In this case,
the corresponding constrainedly profiled PET likelihood ratio function can be defined as
\begin{equation}\label{cPET}
\bar{\ell}_p(\theta_1)=\log\left[\frac{1}{n}\sum_{i=1}^{n}\exp\{\lambda^{\T}(\theta_1)\psi(Z_i;\theta_1)\}\right]-\sum_{j=1}^{q}p_{\gamma}(|\theta_{1j}|),
\end{equation}
where $\lambda(\theta_1)\in\widehat{\Lambda}_n(\theta_1)=\{\lambda:\lambda^{\T}\psi(Z_i;\theta_1)\in\mathcal{E},i=1,\ldots,n\}$ in which
$\mathcal{E}$ is an open interval containing zero.

{\bf Example 1}. As an illustration, we consider a mean regression model: $E(X_i)=\theta$ for $i=1,\ldots,n$.
In this case, estimating equations: $g(X_i,\theta)=X_i-\theta$ can be used to make statistical inference on $\theta$.
Clearly, estimating equations $g(\cdot)$ satisfy unconditional moment model: $E\{g(X_i,\theta)\}=0$.
Under the assumption: $\theta_2^T=(\theta_{21},\ldots,\theta_{2,p-q})=0$, we can obtain the following reduced estimating equations:
$\psi(Z_i,\theta_1)=(X_{i1}-\theta_{11},\ldots, X_{iq}-\theta_{1q})^{\T}$, where $Z_i=(X_{i1},\ldots,X_{iq})^T$ and $\theta_1=(\theta_{11},\ldots,\theta_{1q})^T$.

{\bf Example 2}. We consider a linear regression model: $Y_i=\tilde{X}_i^T\theta+\varepsilon_i$ with $E(\varepsilon_i)=0$, where $\tilde{X}_i=(X_{i1},\ldots,X_{ip})^T$ and $\theta=(\theta_{11},\ldots,\theta_{1q},\theta_{21},\ldots,\theta_{2,p-q})^T$. Let $X_i=\{Y_i,\tilde{X}_i\}$. Thus, estimating equations:
$g(X_i,\theta)=(X_{i1}(Y_i-\tilde{X}_i^T\theta),\ldots,X_{ip}(Y_i-\tilde{X}_i^T\theta))^T$ can be employed
to make statistical inference on $\theta$.
When $\theta_{21}=\cdots=\theta_{2,p-q}=0$, the reduced estimating equations are given by
$\psi(Z_i,\theta_1)=(X_{i1}(Y_i-\tilde{Z}_i^T\theta_1),\ldots, X_{iq}(Y_i-\tilde{Z}_i^T\theta_1))^{\T}$, which satisfies $E\{\psi(Z_i,\theta_1)\}=0$,
where $\tilde{Z}_i=(X_{i1},\ldots,X_{iq})^T$, $Z_i=\{Y_i,\tilde{Z}_i\}$ and $\theta_1=(\theta_{11},\ldots,\theta_{1q})^T$.

\subsection{Selection consistency}

\noindent In this subsection, we investigate the consistency of the above presented model selection.
Let $\mathbb{J}=\{j:\theta_{0j}\neq 0\} $ be the index set of nonzero components of the true parameter
vector $\theta_0$, where $\theta_{0j}$ is the $j$th component of $\theta_0$ for $j=1,\ldots,p$. Denote the cardinality of
$\mathbb{J}$ as $q=|\mathbb{J}|$, which is unknown.
The true parameter vector $\theta_0$ has the following form $\theta_0=(\theta_{10}^{\T},0^{\T})^{\T}$, where $\theta_{10}$ is the true value of $\theta_1$.
The corresponding decomposition of $\hat\theta$ can be written as $\hat\theta=(\hat\theta_1^{\T},\hat\theta_2^{\T})^{\T}$.
Let $\mathcal{D}_n=\{\theta:||\theta-\theta_{0}||\leq \mathbb{C}\sqrt{r/n}\}$ be neighborhoods of $\theta_{0}$ for some constant $\mathbb{C}>0$.
$\Sigma(\theta)=E\{(g(X_i;\theta)-Eg(X_i;\theta))(g(X_i;\theta)-Eg(X_i;\theta))^{\T}\}$,
$\Gamma(\theta)=E\{\partial_{\theta}g(X_{i};\theta)\}$, and $\Sigma=
\Sigma(\theta_0)$, $\Gamma=\Gamma(\theta_0)$.

Lv and Fan (2009) pointed out that the increase of  $p^{'}_{\gamma}(s)$ with respect to $\gamma$ allows for $\gamma$ effectively controlling the overall strength of penalty. Therefore, we can take the tuning parameter $\gamma$ to be our required threshold. Specifically, the selection criterion is $\hat{\mathbb{J}}=\{j: |\hat{\theta}_{j}|>\gamma\}$, where $\hat{\theta}_{j}$ is the $j$th component of the PET estimator $\hat\theta$. Here, our main purpose is to show the selection consistency, i.e., ${\rm Pr}(\hat{\mathbb{J}}= \mathbb{J})\rightarrow 1$ as $n\rightarrow \infty$ even if $p$ diverges with $n$. Note that
the event $\{\hat{\mathbb{J}}\neq \mathbb{J}\}$ is equivalent to the event $\{|\hat{\theta}_{j}|\leq \gamma$ for some $j\in \mathbb{J}\}\cup\{|\hat{\theta}_{j}|> \gamma$ for some $j\in \mathbb{J}^c\}$. Also, we have
\[
\begin{array}{llll}
{\rm Pr}(\{|\hat{\theta}_{j}|\leq \gamma\; {\rm for\; some\; j} \in \mathbb{J}\})&\leq& \sum_{j\in \mathbb{J}}{\rm Pr}(|\hat{\theta}_{j}|\leq \gamma)\\
&=&\sum_{j\in \mathbb{J}}{\rm Pr}(|\theta_{0j}|-|\hat{\theta}_{j}|\ge |\theta_{0j}|-\gamma)\\
&\leq&\sum_{j\in \mathbb{J}}{\rm Pr}(|\theta_{0j}-\hat{\theta}_{j}|\ge |\theta_{0j}|-\gamma)\\
&\leq&\sum_{j\in \mathbb{J}}{\rm Pr}(|\theta_{0j}-\hat{\theta}_{j}|\ge \min_{j\in \mathbb{J}}|\theta_{0j}|-\gamma)\\
&\leq&\sum_{j\in\mathbb{J}} {\rm Pr}(|\hat{\theta}_{j}-\theta_{0j}|\ge \gamma).
\end{array}
\]
The last inequality holds because Assumption \ref{ass4}(i) leads to $\min_{j\in \mathbb{J}}|\theta_{0j}|>2\gamma$. Similarly, we can obtain
${\rm Pr}(\{|\hat{\theta}_{j}|> \gamma$ for some
$j\in \mathbb{J}^c\})\leq\sum_{j\in \mathbb{J}^c}{\rm Pr}(|\theta_{0j}-\hat{\theta}_{j}|\ge \gamma)
\leq\sum_{j\in\mathbb{J}^c}{\rm Pr}(|\hat{\theta}_{j}-\theta_{0j}|\ge \gamma)$.
Combining the above equations leads to ${\rm Pr}(\hat{\mathbb{J}}\neq \mathbb{J})\leq \sum_{j=1}^pP(|\hat{\theta}_{j}-\theta_{0j}|\ge \gamma)$.
Thus, the error probability of selection consistency is
affected by the inconsistency of parameter estimation. Let $\mathbb{E}(A)$ be the eigenvalue of a positive definite matrix $A$.
The following assumptions are required to make statistical inference on $\hat{\theta}$.
\begin{assumption}\label{ass1}
 (Identification, Sieves) (i) The support $\Theta$ of $\theta$ is a compact set in $\mathcal{R}^p$, and $\theta_{0}=(\theta_{10}^{\T},0^{\T})^{\T}\in \Theta$ is the unique solution to $E\{g(X_i;\theta)\}=0$ for $i=1,\ldots,n$; (ii)
$\{\Theta_{s}: s\ge 1\}$  is a sequence of nonempty closed subsets satisfying $\Theta_{s}\subseteq\Theta_{s+1}\subseteq\Theta$,
and there is $\Pi_n\theta_0\in\Theta_{s(n)}$
such that $||\Pi_n\theta_0-\theta_0||=O(\sqrt{r/n})$; (iii) $E\{(||g(X_i;\Pi_n\theta_0)||D(g)^{-1/2})^{\delta}\}<\infty$
for some $\delta>2$ and $D(g)^2n^{2/\delta-1}=o(1)$, where $D(g)$ is the number of moment conditions $g(\cdot)$.
\end{assumption}

\begin{assumption}\label{ass3}
(Sample Moment Criterion) Let $a_0$ and $b_0$ be constants. (i) $a_0\leq \mathbb{E}\{\frac{1}{n}\sum_{i=1}^ng(X_i;\Pi_n\theta_0)g^{\T}(X_i;\Pi_n\theta_0)\}\leq b_0<\infty$ w.p.a.1; (ii) $a_0\leq \sup_{j}E\{g_j(X_i;\Pi_n\theta_0)\}^2\leq b_0<\infty$, $a_0\leq \sup_{j,l,}E\{g_j(X_i;\Pi_n\theta_0)
g_l(X_i;\Pi_n\theta_0)\}^2\leq b_0<\infty$ for $j,l=1,\ldots,r$; (iii) there are $\mathcal{K}_1<\infty$ and $\mathbb{K}_{1}(X_i)$
such that $\sup_{j,l,\theta\in\Theta_{s(n)}}|\partial g_j(X_i;\theta)/\partial\theta_{l}|\leq \mathbb{K}_{1}(X_i)$ and
$E\{\mathbb{K}^2_{1}(X_i)\}\leq \mathcal{K}_1$ for $j=1,\ldots,r$ and $l=1,\ldots,p$;
(iv) there are $\mathcal{K}_2<\infty$ and $\mathbb{K}_{2}(X_i)$ such that
$\sup_{j,l_1,l_2,\theta\in\Theta_{s(n)}}|\partial^2g_j(X_i;\theta)/\partial\theta_{l_1}\partial\theta_{l_2}|\leq \mathbb{K}_{2}(X_i)$ and
$E\{\mathbb{K}^2_{2}(X_i)\}\leq \mathcal{K}_2$ for $j=1,\ldots,r$ and $l_1=q+1,\ldots,p,l_2=1,\ldots,q$.
\end{assumption}

Assumption \ref{ass1}(i) ensures the existence and consistency of the maximizer of Equation (\ref{PETL13}).
It also implies that $\theta_{10}\in \Theta_1$ is the unique solution to $E\{\psi(Z_i;\theta_1)\}=0$ with the sparsity assumption $\theta_{20}=0$ for $i=1,\ldots,n$,
where $\theta_{20}$ is the true value of $\theta_2$.
Although, following Chen and Pouzo (2012), we give the similar definition for the sieve method in Assumption \ref{ass1}(ii),
the PET procedure is quite a good fit for
semiparametric estimation under  slowly growing dimension.
Assumption \ref{ass1}(iii) is proposed to control the tail probability behavior of unconditional moment models by considering the diverging rate of data.
Similar to Chang, Chen and Chen (2015), we can use some function $h(r)>0$ to replace the factor $r^{1/2}$.
Assumption \ref{ass1}(iii) also implies that $E\{\sup_{\theta_1\in\Theta_1}(||\psi(Z_i;\theta_1)||k^{-1/2})^{\delta}\}<\infty$ for some $\delta>2$.
Assumption \ref{ass3}(i) allows for bounding eigenvalues of the corresponding sample matrices in probability.
Assumption \ref{ass3}(ii) can be applied to unconditional  moment constraints $\psi(Z_i;\theta_1)$.
Assumption \ref{ass3}(iii) implies that for any $\theta_1\in\Theta_1$, there exist $\tilde{\mathcal{K}}_1<\infty$ and $\tilde{K}_{1}(Z_i)$ such that $|\partial \psi_j(Z_i;\theta_1)/\partial\theta_{1l}|\leq \tilde{K}_{1}(Z_i)$ and
$E\{\tilde{K}^2_{1}(Z_i)\}\leq \tilde{\mathcal{K}}_1$ for $j=1,\ldots,k$ and $l=1,\ldots,q$. Assumption \ref{ass3}(iv) indicates that only $rq(p-q)$ twice derivatives of moments are bounded, which is less than the number assumed in Leng and Tang (2012). Under Assumption \ref{ass3}, we can estimate the true model well due to its well-posedness.

{\bf Example 3} (Linear instrumental variable regression model). Following Staiger and Stock (1997), a linear instrumental variable (IV) regression model has the following structural equation: $y_i=Y_i^{\T}\theta+\epsilon_i$ for $i=1,\ldots,n$,
where  $Y_i$ is a $p\times1$ vector of endogenous regressors and $\theta$ is a $p\times 1$ vector of parameters of interest,
together with the following reduced equation for $Y_i$:
$Y_i=\mathbb{B} D_i+w_i$,  where $D_i$ is a $K\times 1$ ($K\ge p$) vector of instrument variables, and $\mathbb{B}$ is a $p\times K$ matrix of nuisance parameters.
It is assumed that $v_i =(\epsilon_i, w_i^{\T})^{\T}$ satisfies moment conditions $E(v_i|D_i)=0$ for $i=1,\ldots,n$,
and $v_1|D_1,\ldots,v_n|D_n$ are i.i.d.. Here we consider the following estimating equations: $g(X_i;\theta)=\tilde{\mathcal{X}}_i^{\T}(y_i-\tilde{\mathcal{X}}_i\theta)$, where $\tilde{\mathcal{X}}_i=D_i^{\T}(D^{\T}D)^{-1}D^{\T}Y$,
$D=(D_1,\ldots,D_n)^{\T}$, $Y=(Y_1,\ldots,Y_n)^{\T}$, and $X_i=\{y_i,\tilde{\mathcal{X}}_i\}$. Under the above assumption together with $E(D_i)=0$, we obtain $E\{g(X_i;\theta)\}=0$.
To provide more primitive and transparent regular conditions, we assume ${\rm var}(D_iD_i^{\T})=I_K$. Thus, the identification condition and moment criterion
corresponding to Assumptions $\ref{ass1}$ and $\ref{ass3}$ reduce to
$E||K^{-1/2}D_i||^{4\kappa}<\infty$, $E|w_i^{\T}\theta_0+\epsilon_i|^{4\kappa}<\infty$, and $\mathbb{E}_{\max}(\mathbb{B}\mathbb{B}^{\T})<\infty$,
where  $\kappa$ is some positive integer.

\begin{theorem}\label{th6}
Under Assumptions \ref{ass1}-\ref{ass3}, $r^2=o(n)$ and $rp=o(n)$, for $\theta\in \mathcal{D}_n$, we have
\begin{equation}\label{selec2}
\ell_p(\theta)=-\frac{1}{2}(\theta-\hat{\theta})^{\T}\mathfrak{J}(\theta-\hat{\theta})+R_n,
\end{equation}
where $\mathfrak{J}=J_{0}+J$, $\hat{\theta}=\mathfrak{J}^{-1}(J_{0}\theta_{0}+J\hat{\theta}_{ET})$, $R_n=o_p(r/n)$, $J=\Gamma\Sigma^{-1}\Gamma^{\T}$, and $\hat{\theta}_{ET}$ is the maximum ET likelihood estimator of parameter vector $\theta$, and $J_0$ is a diagonal matrix with the $j$th diagonal element being $J_{0}^{jj}$ for $j=1,\ldots,p$.
\end{theorem}

The $j$th nonzero diagonal element $J_0^{jj}$ can be obtained from the following quadratic approximation of $p_\gamma(|\theta_j|)$ at $\theta_{0j}$:
$p_{\gamma}(|\theta_j|)=p_{\gamma}(|\theta_{0j}|)+\frac{p_{\gamma} '(|(\Pi_n\theta_{0})_j|)(\Pi_n\theta_{0})_j}{|(\Pi_n\theta_{0})_j|}(\theta_j-\theta_{0j})
\approx p_{\gamma}(|\theta_{0j}|)+\frac{p_{\gamma} '(|(\Pi_n\theta_{0})_j|)}{2|(\Pi_n\theta_{0})_j|}(\theta_j^2-\theta_{0j}^2)$
for $\theta_j\in\{\theta_j: |\theta_j-\theta_{0j}|\leq \mathbb{C}\sqrt{r/n}\}$, which leads to $J_0^{jj}=\frac{p_{\gamma} '(|(\Pi_n\theta_{0})_j|)}{|(\Pi_n\theta_{0})_j|}$. 
By Theorem \ref{th6} and the selection criterion $\hat{\mathbb{J}}=\{j: |\hat{\theta}_{j}|>\gamma\}$, if the eigenvalues of matrix $J$ are limited to some finite interval,
the condition $\frac{p_{\gamma} '(|(\Pi_n\theta_{0})_j|)}{|(\Pi_n\theta_{0})_j|}<\frac{p_{\gamma} '(|(\Pi_n\theta_{0})_j|)}{\min_{j\in\mathbb{J}}|(\Pi_n\theta_{0})_j|}\ll\frac{p_{\gamma} '(|(\Pi_n\theta_{0})_j|)}{\gamma}=o(\mathbb{E}_{\min}(J))=o(1)$ implies that $\hat{\theta}_{j}$ can be derived from $\hat{\theta}_{ET}^j$ for $j\in\mathbb{J}$ under Assumption \ref{ass4}(i), where $\hat\theta_{ET}^j$ is the $j$th element of $\hat\theta_{ET}$. This fact shows that we can obtain selection consistency of $\hat{\theta}_{\mathbb{J}}$ from Theorem \ref{th6} and the selection criterion when the ET likelihood dominates the penalty function. On the other hand, when the penalty function dominates the ET likelihood and $\frac{p_{\gamma} '(|(\Pi_n\theta_{0})_j|)}{|(\Pi_n\theta_{0})_j|}>\frac{p_{\gamma} '(|(\Pi_n\theta_{0})_j|)}{\gamma}\gg\mathbb{E}_{\max}(J)$ for $j\in\mathbb{J}^c$, which leads to $\hat{\theta}_{j}=0$ because of the sparsity assumption of $\theta_{\mathbb{J}^c}=0$, we can obtain selection consistency of $\hat{\theta}_{\mathbb{J}^c}$. In what follows, we will discuss the error probability of the event $\hat{\mathbb{J}}\neq\mathbb{J}$, and bound it by using some special moment conditions.
To this end, similar to Bondell and Reich (2012), we denote $Q^*$ as a $p\times p$ matrix whose column vectors are eigenvectors of matrix $\Gamma\Sigma^{-1}\Gamma^{\T}$ corresponding to its $p$ eigenvalues $d_1^*\ge\ldots\ge d_t^*>0=d_{t+1}^*=\ldots=d_p^*$. Denote $Q^*=(Q_1^*,Q_2^*)$, where $Q_1^*$ denotes the first $t$ columns of $Q^*$, those corresponding
to the nonzero eigenvalues and $Q_2^*$ are the remaining $p-t$ columns of $Q^*$. Since $Q^*$ is an orthonormal basis, we have $\theta_0=Q^*\eta=Q_1^*\eta_1+Q_2^*\eta_2$, where $\eta_1$ and $\eta_2$ are the corresponding partition of $\eta$, i.e., $\eta=(\eta_1^{\T},\eta_2^{\T})^{\T}$. To obtain selection consistency, we need the following assumption.

\begin{assumption}\label{ass4}
 (Penalty Criterion)  The penalty function $p_\gamma(t)$ is lower semicompact, and
 (i) $\max_j|\theta_{0j}|<\infty$, $\gamma$ satisfies $\gamma/\min_{j\in\mathbb{J}}|(\Pi_n\theta_{0})_j|\rightarrow 0$;
 (ii) $\gamma\rightarrow 0$ together with $\frac{\log p}{\gamma\min_{j\in\mathbb{J}^c}p_{\gamma}^{'}(|(\Pi_n\theta_{0})_j|)}\rightarrow 0$ as $n\rightarrow \infty$;
(iii) $\gamma\rightarrow 0$ together with $\max_{j\in\mathbb{J}}p_{\gamma}^{'}(|(\Pi_n\theta_{0})_j|)\frac{\sqrt{q}}{\gamma^2}\rightarrow 0$  as $n\rightarrow \infty$; (iv) $||Q_2^*\eta_2||_{\infty}=$
$O(\frac{\max_{j\in\mathbb{J}}p_{\gamma}^{'}(|(\Pi_n\theta_{0})_j|)\sqrt{q}}{\gamma})$  as $n\rightarrow \infty$.
\end{assumption}

Assumption \ref{ass4}(i) assumes that the true values of parameters are finite,
and the rate of the threshold $\gamma$ decreasing to zero is faster than the rate on the magnitude of the true nonzero coefficients, which is guaranteed to remain identifiable from zero. Assumption \ref{ass4}(ii) gives the rate at which
the threshold may decrease to zero, while still allowing for the
exclusion of the unimportant predictors with enough large $\min_{j\in\mathbb{J}^c}p_{\gamma}^{'}(|(\Pi_n\theta_{0})_j|)\gg \log(p)/\min_{j\in\mathbb{J}}|(\Pi_n\theta_{0})_j|$
by combining this Assumption and Assumption \ref{ass4}(i). If the threshold diverges
too quickly, the bias of estimator will not enough quickly vanish. Because Assumption \ref{ass4}(iii) implies $\max_{j\in\mathbb{J}}p_{\gamma}^{'}(|(\Pi_n\theta_{0})_j|)\frac{\sqrt{q}}{\gamma}\rightarrow 0$, we have $\min_{j\in\mathbb{J}}|(\Pi_n\theta_{0})_j|\gg(\max_{j\in\mathbb{J}}p_{\gamma}^{'}(|(\Pi_n\theta_{0})_j|)\sqrt{q})^{1/2}$
under Assumptions \ref{ass4}(i) and \ref{ass4}(iii).
Assumptions \ref{ass4}(iii) and \ref{ass4}(iv) show that the true parameter lies in a linear space
spanned by $\Gamma\Sigma\Gamma^{\T}$, which is a basic identifiability condition for identifying the true nonzero parameters.
Note that true parameters actually lie in a $q$-dimensional subspace. Thus, only if one function of true parameters is
estimable within the linear space, we allow that $q/n$ may even diverge. We can allow for the case that
$\theta$ may actually be sparse in some linear transformed space, but it is
not sparse in the original space. Note that Assumption \ref{ass4}(iv) is only a possible formulation
for which the assumption holds, but it is not the only way to
satisfy the assumption.

\begin{theorem}\label{th7}
Under Theorem \ref{th6} and Assumption \ref{ass4}, as $n\rightarrow \infty$, we have
$$P(\hat{\mathbb{J}}\neq \mathbb{J})\leq2\frac{p}{\sqrt{\gamma\min_{j\in\mathbb{J}^c}p_{\gamma}^{'}(|(\Pi_n\theta_{0})_j|)}}\exp\left\{-\frac{\gamma\min_{j\in\mathbb{J}^c}p_{\gamma}^{'}(|(\Pi_n\theta_{0})_j|)}{8}\right\}\rightarrow 0.$$
\end{theorem}

Theorem \ref{th7} shows that the above proposed parameter selection procedure is consistent, i.e., $P(\hat{\mathbb{J}}= \mathbb{J})\rightarrow 1$. It also implies
that $\hat{\theta}_2$=0 with probability tending to 1.
In the following corollary, we obtain the rate of $\max_{j\in\mathbb{J}^c}p_{\gamma}^{'}(|(\Pi_n\theta_{0})_j|)$ (hence the rate of $\gamma$) when Assumptions \ref{ass4}(i)-(iv) hold and $p$ diverges at its fastest possible rate.

\begin{corollary}\label{cor2}
Suppose that $\frac{\max_{j\in\mathbb{J}}p_{\gamma}^{'}(|(\Pi_n\theta_{0})_j|)}{\gamma}=O(\frac{\log(p)}{n})$,
$\max_{j\in\mathbb{J}^c}p_{\gamma}^{'}(|(\Pi_n\theta_{0})_j|)=O(\frac{n}{\sqrt{q}})$ and Assumptions \ref{ass4}(i)-(iv) hold.
The PET thresholding parameter selection procedure is also consistent when $p$ satisfies $\log(p)=O((n/\sqrt{q})^c)$ for some $0<c<1$.
\end{corollary}

When $q$ does not grow with $n$, i.e., the true number of nonzero parameters is fixed, we can allow for an exponential growing case, for example,
$\log(p)=O(n^c)$ for $0<c<1$. If $\max_{j\in\mathbb{J}}p_{\gamma}^{'}(|(\Pi_n\theta_{0})_j|)/\gamma=O(\log(p)/n)$ and
$\max_{j\in\mathbb{J}^c}p_{\gamma}^{'}(|(\Pi_n\theta_{0})_j|)=O(n/\sqrt{q})$, Assumptions \ref{ass4}(ii) and \ref{ass4}(iii)
are equivalent to $\log(p)\sqrt{q}/(n\gamma)$ $\rightarrow 0$,
which corresponds to the rate of threshold $\gamma$.

\subsection{Asymptotic properties of the SCAD-PET estimator}

\noindent In this subsection, we discuss consistent and oracle properties of the SCAD-PET estimator $\hat{\theta}_1$ , and show that the constrainedly profiled PET likelihood ratio function like ET likelihood ratio function is
asymptotically distributed as the chi-squared distribution.

Denote $\bar\psi(\theta_1)=n^{-1}\sum_{i=1}^n\psi(Z_i;\theta_1)$, $\Gamma_1(\theta_1)=E\{\partial_{\theta_1}\psi(Z_{i};\theta_1)\}$,
$\Sigma_1(\theta_1)=E\{(\psi(Z_i;\theta_1)$ $-E\psi(Z_i;\theta_1))(\psi(Z_i;\theta_1)-E\psi(Z_i;\theta_1))^{\T}\}$, $\Sigma(\theta)=E\{(g(X_i;\theta)-Eg(X_i;\theta))(g(X_i;\theta)-Eg(X_i;\theta))^{\T}\}$,
$K_1(\theta_1)=\{\Gamma_1^{\T}(\theta_1)\Sigma^{-1}_1(\theta_1)\Gamma_1(\theta_1)\}^{-1}$.
Let $\mathcal{D}_{1n}=\{\theta_1:||\theta_1-\theta_{10}||\leq \mathbb{C}\sqrt{k/n}\}$ be neighborhoods of $\theta_{10}$ for some constant $\mathbb{C}>0$.

\begin{assumption}\label{ass5}
(i) There are two positive constants $a_0$ and $b_0$ such that (i) the eigenvalue of $\Gamma_1^{\T}(\theta_1)\Gamma_1(\theta_1)$ satisfies $a_0\leq \mathbb{E}(\Gamma_1^{\T}(\theta_1)\Gamma_1(\theta_1))\leq b_0<\infty$ for all $\theta_1\in \mathcal{D}_{1n}$;
(ii) for any $\theta_1\in \mathcal{D}_{1n}$, there are $\mathcal{K}_3<\infty$ and $\mathbb{K}_{3}(Z_i)$ such that $|\partial^3\psi_j(Z_i;\theta_1)/\partial\theta_{1l_1}\partial\theta_{1l_2}\partial\theta_{1l_3}|\leq \mathbb{K}_{3}(Z_i)$ and
$E\{\mathbb{K}^2_{3}(Z_i)\}\leq \mathcal{K}_3$ for $j=1,\ldots,k$ and $l_1,l_2,l_3=1,\ldots,q$.
\end{assumption}

\begin{assumption}\label{ass6}
(i) There are two positive functions $\zeta_1(k,q)$ and $\zeta_2(\varepsilon)$ such that for any $\varepsilon$,
$\inf_{\{\theta_1\in\Theta_1:||\theta_1-\theta_{10}||\ge\varepsilon\}}||E\psi(Z_i;\theta_1)||\ge\zeta_1(k,q)
\zeta_2(\varepsilon)>0$,
where $\lim\inf_{k,q\rightarrow\infty}\zeta_1(k,q)>0$;
(ii) $\sup_{\theta_1\in\Theta_1}||\bar{\psi}(\theta_1)-E\psi(Z_i;\theta_1)||=o_p\{\zeta_1(k,q)\}$.
\end{assumption}

\begin{assumption}\label{ass7}
 Suppose the penalty function $p_\gamma(t)$ is lower semicompact and satisfies $\max_{j\in\mathbb{J}}p_{\gamma}(|\theta_{0j}|)\leq \mathbb{C}\gamma$, where $\mathbb{C}$ is some constant and $\gamma=O(k/\{nq\})$.

 \end{assumption}
\vspace{0.1cm}

Assumption \ref{ass5}(i) shows that the constaints on eigenvalues of matrix $\Gamma_1(\theta_1)\Gamma_1^{\T}(\theta_1)$ is a relaxed condition of Chang et al. (2015).
Assumption \ref{ass5}(ii) is used to control the order of the remainder term when taking the third-order expansion of the objective function.
Assumption \ref{ass6}(i) is the population identification condition for the diverging parameter space.
Assumption \ref{ass6}(ii) is an extension of the uniform convergence, whose detailed interpretation can refer to Chang et al. (2015).
The lower semicompact penalty in Assumption \ref{ass7} has been used by Chen and Pouzo (2012),
and implies that the effective parameter space converts an ill-posed problem into a well-posed one.
Assumptions \ref{ass4} and \ref{ass7} hold for many penalty functions such as the SCAD penalty function and the hard-threshold penalty.
However, for $L_1$ penalty, $\gamma=p'_\gamma(|\theta_{0j}|)=O(\frac{k}{nq})$ in Assumption \ref{ass7} is in conflict
with $\gamma$ supposed in Assumption \ref{ass4}, which implies that the
PET likelihood estimator with the $L_1$ penalty generally cannot achieve the consistency rate of $O_p(\sqrt{k/n})$ established in Theorem \ref{th1},
and has not the oracle property established in Theorem \ref{th2} when the number of parameters diverges with sample size $n$.
In fact, the above mentioned issue has been pointed out by Fan and Li (2001) and Zou (2006) even for the finite $p$.

\begin{theorem}\label{th1}{\rm (}Consistency of PET Estimator of $\hat{\theta}_1${\rm )}~
Under Assumptions \ref{ass1}-\ref{ass7},
there is a strict local maximizer $\hat{\theta}=(\hat{\theta}_1^{\T},\hat{\theta}_2^{\T})$ of the PET likelihood $\ell_p(\theta)$ such that $\hat{\theta}_2=0$ with probability tending to 1  as $n\rightarrow \infty$ and
$||\hat{\theta}_1-\theta_{10}||=O_p(\sqrt{k/n})$. 
\end{theorem}

Theorem \ref{th1} establishes the consistency of the proposed PET estimator $\hat{\theta}_1$, that is, there is a root-($n/k$)-consistent PET estimator of $\theta_1$.
It also shows that the sparsity property of the proposed PET estimator $\hat\theta$ is still valid, that is,
zero components in $\theta_0$ are estimated as zero with probability tending to one under Theorem \ref{th1}. Generally,
it is right to assume $q/k<1$ in Theorem \ref{th1} because our main goal is to select nonzero parameters.

\begin{assumption}\label{ass8}
 The tuning parameter $\gamma$ and second derivatives of the penalty function $p_\gamma(t)$ satisfy $\gamma=o(n^{1/4})$ and
 $\max_{j\in\mathbb{J}}p_{\gamma}^{''}(|\theta_{10j}|)=o(1/\sqrt{kq})$, respectively.
 \end{assumption}

Under Assumption \ref{ass4}(iii) and $\gamma=o(n^{1/4})$, we have
$\max_{j\in\mathbb{J}}p_{\gamma}^{'}(|\theta_{10j}|)=o(1/\sqrt{nq})$, which is a useful conclusion in the proof of Theorem \ref{th2}.
Clearly, Assumption \ref{ass8} holds for the SCAD penalty function.

\begin{theorem}\label{th2}{\rm (}Oracle Property{\rm )}~ Under Assumptions \ref{ass1}-\ref{ass8}, we have

{\rm (i) (}Sparsity{\rm )} $\hat\theta_2=0$ with probability tending to one.

{\rm (ii) (}Asymptotic normality{\rm )} $\sqrt{n}G_n\mathcal{K}^{-1/2}(\hat{\theta}_1-\theta_{10})\stackrel{{\cal L}}{\rightarrow}\mathcal {N}(0,V)$ when $k^2(k+q)^3=o(n)$, where
$G_n$ is a $d\times q$ matrix such that $G_nG_n^{\T}\rightarrow V$, $V$ is a $d\times d$ nonnegative symmetric matrix with the fixed $d$,
$\mathcal{K}=K_1(\theta_{10})$, and $\stackrel{\cal L}{\rightarrow}$ denotes convergence in distribution.

\end{theorem}

Theorem \ref{th2} shows that the sparsity and asymptotic normality of the proposed PET estimator still hold when the number of parameters diverges.
Under different assumptions, we can obtain different diverging rates in Theorems \ref{th2}(i) and \ref{th2}(ii)
by controlling the remainder term of the Taylor's expansion.
Note that when $k$ grows with $q$ and $q/k\rightarrow \kappa\in(0,1]$, Theorem \ref{th2}(ii) is similar
to that given in Leng and Tang (2012) for the penalized empirical likelihood estimator with $q=o(n^{1/5})$.

Let $\mathcal{P}=\Theta\times\widehat{\Lambda}_n$, $S(\Delta)=\bar{\ell}_p(\theta_1)$ in which $\Delta=\{\lambda,\theta_1\}$, $S_1(\Delta)=\partial S(\Delta)/\partial\lambda$, $S_2(\Delta)=\partial S(\Delta)/\partial\theta_1$, $S_{l,j}(\Delta)$ be the $j$th component of $S_l(\Delta)$ for $l=1,2$ and $j=1,\ldots,k+q$, and $\mathbb{E}_{\max}(A)$ be the maximum eigenvalue of matrix $A$. If Assumptions \ref{ass1}-\ref{ass4}, \ref{ass6}-\ref{ass8},
$\sup_{\Delta\in\mathcal{P}}\mathbb{E}_{\max}$ $\{\partial_\Delta^{2}S_{l,j}(\Delta)\}=O_p(1)$ in probability for $l=1,2$ and $j=1,\ldots,k+q$ hold
and $q/k\rightarrow \kappa\in(0,1]$, the asymptotic normality of Theorem \ref{th2}(ii) is still valid when $k^2(k+q)=o(n)$ or $q=o(n^{1/3})$, which is assumed in Fan and Lv (2011).

The above results are established on the basis of the corrected specification of unconditional moment models. However, in some applications, unconditional moment models may be misspecified.
Hence, it is necessary to study the asymptotic properties of the PET estimators of parameters in the presence of misspecified unconditional moment models.
To this end, we need the following regularity conditions.

\begin{assumption}\label{ass9}
The function $Q(\theta)=\log E\exp\{\nu^{*\T}(\theta)g(X_i;\theta)\}-\sum_{j=1}^{p}p_{\gamma}(|\theta_j|)
$ is maximized at a unique ``pseudo-true" value $\theta^{*}=(\theta_1^{*\T},\theta_2^{*\T})^{\T}\in {\rm int}(\Theta)$ of $\theta$, where $\theta_2^*=0$
and ${\rm int}(\Theta)$ is the inner set of the compact set $\Theta$.
\end{assumption}

\begin{assumption}\label{ass10} Functions $g(X_i;\theta)$ and $p_\gamma(|\theta|)$ are  continuous with respect to parameter vector $\theta\in\Theta$ (or components of $\theta$).
\end{assumption}
\begin{assumption}\label{ass11}  There are a function $H(X_i)$ and a finite constant $H^*<\infty$ such that $\sup_{\theta\in\Theta}\sup_{\nu\in\widehat{\mathcal{V}}_n(\theta)}\exp\{\nu^{\T}g(X_i;$ $\theta)\}<H(X_i)$ and $E(H(X_i))^2\leq H^*$, where $\widehat{\mathcal{V}}_n(\theta)$ is a compact set such that $\nu^{*}(\theta)\in {\rm int}(\widehat{\mathcal{V}}_n(\theta))$.
\end{assumption}
\begin{assumption}\label{ass12} Functions $\Omega_{jl}(X_i;\theta)=\partial^2g(X_i;\theta)/\partial\theta_j\partial\theta_l$ are continuous with respect to $\theta$ in the neighborhood
$\mathcal{Q}^*$ of $\theta^*$ for $j,l=1,\ldots,p$.
\end{assumption}

\begin{assumption}\label{ass13}  As $n\rightarrow \infty$, $\lim \inf_{\gamma\rightarrow 0}\lim \inf_{\theta\rightarrow 0^+}p_{\gamma}^{'}(\theta)/\gamma>0$ and $\gamma$
satisfies $||\nu(\theta)||=o(\gamma)$.
\end{assumption}

\begin{assumption}\label{ass14}  There is a function $f(X_i)$ satisfying
$$E\left\{\sup_{\theta\in\cal Q^*}\sup_{\nu\in\widehat{\mathcal{V}}_n(\theta)}\exp\{k_1\nu^{\T}g(X_i;\theta)\}(f(X_i))^{k_2}\right\}<\infty~~{\rm for}~k_1, k_2=0,1,2$$
such that $||g(X_i;\theta)||\leq f(X_i)$, $||\Gamma(X_i;\theta)||\leq f(X_i)$,
$||\Omega_{jl}(X_i;\theta)||\leq f(X_i)$ for $\forall X_i\in\mathcal{R}^{\iota}$, $\forall \theta\in\cal Q^*$ and $j,l=1,\ldots,p$, where $\Gamma(X_i;\theta)=\partial_{\theta}g(X_i;\theta)$. The first and second derivatives of the penalty function $p_\gamma(t)$ satisfy $\max_{j\in\mathbb{J}}p_{\gamma}^{'}(|\theta_{0j}|)\leq \mathbb{C}\gamma$ and $\max_{j\in\mathbb{J}}p_{\gamma}^{''}(|\theta_{0j}|)\leq \mathbb{C}\gamma$.
\end{assumption}

Assumption \ref{ass9} is employed to ensure the existence of the PET estimator maximizing the objective function $\ell_p(\theta)$ with the sparsity condition.
The similar sparsity condition has been adopted in Lu, Goldberg and Fine (2012),
which showed that there are zero components in the ``pseudo-true" value, and the notion of an oracle
estimator is defined in terms of such sparseness.
Assumptions \ref{ass10}-\ref{ass12} for the continuity and boundness satisfying the conditions of Slutsky Theorem are adopted
to ensure the consistency of the proposed PET estimator under misspecification.
Assumption \ref{ass13} shows that the order of the lagrange multiplier of the PET likelihood is controlled by that of the tuning parameter in penalty function,
which plays an important role on the local concave optimization problem of the PET likelihood.
This ensures the sparsity property of the PET estimator in the presence of misspecified unconditional moment restrictions.
The boundness in Assumption \ref{ass14} is designed to satisfy the conditions given in Theorem \ref{th3}.4 of Newey and McFadden (1994) for asymptotic normality under the just-identified case.

\begin{theorem}\label{th3}
(large sample properties under misspecification).  Under Assumptions \ref{ass9}-\ref{ass14}, as $n\rightarrow \infty$, we have

{\rm (i) (}Consistency{\rm )} $\hat{\theta}\stackrel{{\cal P}}{\rightarrow}\theta^{*}$, where $\stackrel{\cal P}{\rightarrow}$ denotes convergence in probability;

{\rm (ii) (}Sparsity{\rm )} $\hat\theta_2=0$ with probability tending to one;

{\rm (iii) (}Asymptotic normality{\rm )} Let $\mathcal{G}=E\{\partial_\phi\Psi(X_i;\phi)\}_{\phi=\phi^{*}}$,
$\Xi=E\{\Psi(X_i;\phi^{*})\Psi^{\T}(X_i;\phi^{*})\}$. If $\mathcal{G}$ is a nonsingular matrix, we have $n^{1/2}(\hat{\phi}-\phi^{*})\stackrel{{\cal L}}{\rightarrow}\mathcal {N}(0,\mathcal{G}^{-1}\Xi\mathcal{G}^{-\T})$, where $\Psi$ and $\phi$ are defined in Lemma \ref{lem5} of the Appendix.
\end{theorem}

Theorem \ref{th3}(i) indicates that the nonexistence of convergence rate is due to misspecified unconditional moment restrictions.
The consistency holds only when the objective PET likelihood function converges to its population form satisfying some continuity and boundness assumptions.
The sparsity of Theorem \ref{th3}(ii) is derived from the existence of zero components in ``pseudo-true" value of $\theta$.
Combining Theorems \ref{th3}(i) and \ref{th3}(ii) yields $\hat{\theta}_1\stackrel{{\cal P}}{\rightarrow} \theta_1^*$.
Theorem \ref{th3}(iii) establishes the asymptotic normality of nonzero components of $\theta$ or $\phi$ based on the objective function defined
in Lemma \ref{lem5}.

\subsection{Constrainedly PET likelihood ratio test}

\noindent In this subsection, similar to Fan and Peng (2004) and Leng and Tang (2012), we consider testing linear hypotheses of the form
$$H_0: B_n\theta_{10}=0 ~~{\rm versus}~~ H_1: B_n\theta_{10}\neq 0,$$
where $B_n$ is a user-specified $d\times q$ matrix such that $B_nB_n^{\T}=I_d$ for a fixed $d$, and $I_d$ is a $d\times d$ identity matrix.
The above hypotheses include testing individual and multiple components of $\theta_{10}$ as special cases (Fan and Peng, 2004; Leng and Tang, 2012). For example, the null hypothesis $H_{0j}: \theta_{10j}=0$ for some $j\in\{1,\ldots,q\}$ can be written as the null hypothesis $H_{0j}: B_n\theta_{10}=0$ in which $B_n$ is a $1\times q$ vector with the $j$th element being 1 and 0 elsewhere, where $\theta_{10j}$ is the $j$th component of $\theta_{10}$.
Inspired by Fan and Peng (2004) and Leng and Tang (2012), we consider the following constrainedly PET likelihood ratio statistic
\begin{equation}\label{PETL23}
\hat{\ell}(B_n)=2n\{\bar{\ell}_p(\hat{\theta}_1)-\max_{\theta_{1},B_n\theta_{1}=0}\bar{\ell}_p(\theta_{1})\}
\end{equation}
for testing $H_0$: $B_n\theta_{10}=0$.

\begin{theorem}\label{th4}
Under the null hypothesis and Assumptions \ref{ass1}-\ref{ass8}, we have $\hat{\ell}(B_n)\stackrel{{\cal L}}{\rightarrow} \chi_d^2$ as $n\rightarrow \infty$, where $\chi_d^2$ denotes the chi-squared distribution with $d$ degrees of freedom.
\end{theorem}

Theorem \ref{th4} establishes the asymptotic distribution of the above presented test statistic $\hat\ell(B_n)$ under the null hypothesis $H_0$: $B_n\theta_{10}=0$, which indicates that the well-known Wilk's phenomenon for the likelihood, empirical likelihood (Owen, 2001) and adjusted ET likelihood (Zhu et al., 2009) functions is still valid for the PET likelihood with a growing number of parameters. Thus, we extend the result given in Zhu et al. (2009) to unconditional moment models with a growing number of parameters.
The above asymptotic result can be used to simultaneously test statistically significance of several covariates by taking some specific matrix $B_n$ (Fan and Peng, 2004).
To wit, we can use the above asymptotic result to identify zero and nonzero components of $\theta_1$.
Also, it can be adopted to construct asymptotic confidence region of $B_n\theta_1$. The PET-likelihood-ratio-test-based $100(1-\alpha)\%$ approximate confidence region for $B_n\theta_{1}$ is given by
\begin{equation}\label{PETCI}
R_{\alpha}=\left\{\xi:2n\{\bar{\ell}_p(\hat{\theta}_1)-\max_{\theta_1,B_n\theta_1=\xi}\bar{\ell}_p(\theta_1)\}\leq\mathcal {\chi}_d^2(1-\alpha)\right\},
\end{equation}
where $\mathcal {\chi}_d^2(1-\alpha)$ is the $(1-\alpha)$-quantile of the chi-squared distribution with $d$ degrees of freedom.

\subsection{High-order asymptotic properties}

\noindent In this subsection, we present the high-order asymptotic property of the proposed PET estimator under some proper regularity conditions,
which is an extension of the high-order results for the ET estimator in Schennach (2007).

Let $\eta=(\theta_1^{\T},\lambda^{\T})^{\T}$, $\eta_0=(\theta_{10}^{\T},0^{\T})^{\T}$, $\hat{\eta}$ be the constrained PET estimator of $\eta$, $\Gamma_{1i}(\theta_1)=\partial \psi(Z_i;\theta_1)/\partial\theta_1^{\T}$ and $\psi_i(\theta_1)=\psi(Z_i;\theta_1)$.  Theorem \ref{th2} has established the following first-order conditions:
 \begin{equation}\label{high0}
0=\frac{1}{n}\sum\limits_{i=1}^nm(Z_i;\hat{\eta})=\frac{1}{n}\sum\limits_{i=1}^n\rho_i(\hat{\eta})
\begin{pmatrix}
\Gamma_{1i}(\hat{\theta}_1)\hat{\lambda} \\
\psi_i(\hat{\theta}_1)
\end{pmatrix}
-\begin{pmatrix}
 w(\hat{\theta}_1)\\
0
\end{pmatrix},
\end{equation}
where $\rho_i(\eta)=n\pi_i=n\exp\{\lambda^{\T}\psi(Z_{i};\theta_1)\}/\sum_{j=1}^{n}\exp\{\lambda^{\T}\psi(Z_j;\theta_1)\}$, and the $j$th component of vector $w(\theta_1)$ has the form of  $p_{\gamma}^{'}(|\theta_{1j}|){\rm sign}(\theta_{1j})$ for $j=1,\ldots,q$. By (\ref{high0}), we derive the stochastic expansion of the constrained PET estimator based on the following additional smoothness and moment conditions.

\begin{assumption}\label{ass15}
(i) Let $\mathcal{S}=k+q$, and $m(Z_i;\eta)$ be four-order continuously differentiable with respect to $\eta\in\{\eta:||\eta-\eta_0||\leq \sqrt{k/n}\}$; (ii) there are
$\mathcal{M}_t<\infty$ and $M_{t}(Z_i)$ such that $\partial^tm_j(Z_i;\eta)/\partial\eta_{l_1}\cdots\partial\eta_{l_t}\leq M_{t}(Z_i)$ and
$E\{M^2_{t}(Z_i)\}\leq \mathcal{M}_t$ for $j=1,\ldots,\mathcal{S}$, $l_1,\ldots,l_t=1,\ldots,\mathcal{S}$, $t=1,\ldots,4$;
(iii) $\mathbb{M}=E\partial m(Z_i;\eta_0)/\partial\eta^{\T}$ exists and has finite eigenvalue.
\end{assumption}
Although Assumption \ref{ass15} is a general condition in restricting the boundness of derivatives and eigenvalues, it has an overlap with Assumptions \ref{ass3} and \ref{ass5}.  Then, we obtain

\begin{theorem}\label{th5}
Under Theorems \ref{th1} and \ref{th2}, Assumption \ref{ass15} and $q/k\rightarrow \kappa$, where $\kappa$ is some constant, as $n\rightarrow \infty$, we have
\begin{equation}\label{high01}
\hat{\eta}-\eta_0=\frac{\tilde{\upsilon}}{\sqrt{n}}+\frac{Q_1(\tilde{\upsilon})+Q_2(\tilde{\upsilon},\tilde{A})}{n}
+\frac{Q_3(\tilde{\upsilon},\tilde{A})+
Q_4(\tilde{\upsilon},\tilde{A})+Q_5(\tilde{\upsilon})+Q_6(\tilde{\upsilon})}{n\sqrt{n}}+R_n,
\end{equation}
where $\tilde{\upsilon}=-\frac{1}{\sqrt{n}}\sum_{i=1}^n\mathbb{M}^{-1}m(Z_i;\eta_0)$, $\tilde{A}=\frac{1}{\sqrt{n}}\sum_{i=1}^n\partial m(Z_i;\eta_0)/\partial\eta^{\T}-\sqrt{n}\mathbb{M}$, $Q_1(\tilde{\upsilon})=-(2\mathbb{M})^{-1}\sum_{j=1}^{\mathcal{S}}\tilde{\upsilon}_jM_j^*\tilde{\upsilon}$, $\tilde{\upsilon}_j$ is the $j$th element of $\tilde{\upsilon}$ and $M_j^*=E\{\partial^2m(Z_i;\eta_0)/\partial\eta_j\partial\eta^{\T}\}$,
$Q_2(\tilde{\upsilon},\tilde{A})=-\mathbb{M}^{-1}\tilde{A}\tilde{\upsilon}$, $Q_3(\tilde{\upsilon},\tilde{A})=-(4\mathbb{M})^{-1}$
$\{-4\tilde{A} \mathbb{M}^{-1}\tilde{A}\tilde{\upsilon}-\sum_{j=1}^{\mathcal{S}}(\tilde{a}_j+2\tilde{b}_j)M_j^*\tilde{\upsilon}+2\sum_{j=1}^{\mathcal{S}}\tilde{\upsilon}_j\tilde{B}_j\tilde{\upsilon}\}$,
$\tilde{B}_j=\frac{1}{\sqrt{n}}\sum_{i=1}^n\partial^2m(Z_i;\eta_0)/\partial\eta_j\partial\eta^{\T}-\sqrt{n}M_j^*$, $\tilde{a}_j$ and $\tilde{b}_j$ are the $j$th elements of $\tilde{a}=\mathbb{M}^{-1}\sum_{j=1}^{\mathcal{S}}\tilde{\upsilon}_jM_j^*\tilde{\upsilon}$ and $\tilde{b}=\mathbb{M}^{-1}\tilde{A}\tilde{\upsilon}$, respectively,
$Q_4(\tilde{\upsilon},\tilde{A})=(2\mathbb{M})^{-1}\{\tilde{A}\mathbb{M}^{-1}$ $\sum_{j=1}^{\mathcal{S}} \tilde{\upsilon}_jM_j^*\tilde{\upsilon}+\sum_{j=1}^{\mathcal{S}}\tilde{\upsilon}_jM_j^*\tilde{b}\}$,
$Q_5(\tilde{\upsilon})=-(6\mathbb{M})^{-1}\sum_{j,l=1}^{\mathcal{S}}\tilde{\upsilon}_j\tilde{\upsilon}_lM_{jt}^*\tilde{\upsilon}$ with
$M_{jt}^*=E\{\partial^3m(Z_i;$ $\eta_0)/\partial\eta_j\partial\eta_t\partial\eta^{\T}\}$,
$Q_6(\tilde{\upsilon})=-(2\mathbb{M})^{-1}\sum_{j=1}^{\mathcal{S}}\tilde{\upsilon}_jM_j^*$ $Q_1(\tilde{\upsilon})$ and
$R_n=O_p(q^{9/2}/n^2)$.
\end{theorem}

By Theorem \ref{th5}, we have $||Q_1(\tilde{\upsilon})||=||Q_3(\tilde{\upsilon},\tilde{A})||=O_p(q^{5/2})$, $||Q_2(\tilde{\upsilon},\tilde{A})||=O_p(q^{3/2})$, $||Q_4(\tilde{\upsilon},\tilde{A})||=O_p(q^{7/2})$, $||Q_5(\tilde{\upsilon},\tilde{A})||=O_p(q^{4})$, and $||Q_6(\tilde{\upsilon},\tilde{A})||=O_p(q^{9/2})$.
Because $||Q_1(\tilde{\upsilon})/n||$ has the largest order $O_p(q^{5/2}/n)$ among all terms except for $\tilde{\upsilon}/\sqrt{n}$, the conditions $q=o(n^{1/5})$ and $q/k\rightarrow \kappa<1$ given in remark of Theorem \ref{th2}(ii) become the sufficient condition of $O_p(q^{5/2}/n)=o_p(n^{-1/2})$.
Furthermore, it follows from Theorem \ref{th5} that asymptotic (higher-order) bias of the proposed PET estimator for nonzero parameter vector is given by
$${\rm Bias}(\hat{\theta}_1)=E\{Q_1(\tilde{\upsilon})+Q_2(\tilde{\upsilon},\tilde{A})\}/n.$$
To investigate the precision of the bias, we introduce the following notations. Let $\mathbb{H}=\mathcal{K}\Gamma_1^{\T}\Sigma_1^{-1}$,
$A=\mathcal{K}\mathcal{W}^2\mathcal{K}^{\T}$, $B^{\T}=\mathcal{K}\mathcal{W}^2\mathbb{H}$, $\mathcal{C}=\mathbb{H}^{\T}\mathcal{W}^2\mathbb{H}$,
$\mathcal{W}=\mathcal{W}(\theta_{10})={\rm diag}(\omega_{11},\ldots,\omega_{qq})$ be an $q\times q$ diagonal matrix, where
$\mathcal{K}$ is defined in Theorem \ref{th2}, and $\omega_{jj}=p'_\gamma(|\theta_{10j}|)$ for $j=1,\ldots,q$.
Denote $\mathcal{A}_1=b_1+a_2+b_2$, $\mathcal{A}_2=c_1+2c_2+d_2$, where the $j$th component of $b_1$ is $b_{1j}={\rm tr}(B^{\T}E\{\partial^2\psi_i(\theta_{10})/\partial\theta_{1j}\partial\theta_1^{\T}\})/2$ for $j=1,\ldots,q$, $a_2=\sum_{j=1}^kE\{B^{\T}e_j\partial^2 \psi_{ij}/\partial\theta_1\partial\theta_1^{\T}\}/2$ in which $\psi_{ij}$ is the $j$th element of $\psi_i$ and $e_j$ is a $k\times 1$ vector whose $j$th element is 1 and 0 elsewhere, $b_2=E({\Gamma_{1i}^j}^{\T}\mathcal{C}\psi_i)$ in which $\Gamma_{1i}^j=\partial^2\psi_i(\theta_{10})/\partial\theta_{1j}\partial\theta_1^{\T}$, the $j$th element of $c_1$ is $c_{1j}={\rm tr}(AE\{\partial^2\psi_{ij}(\theta_0)/\partial\theta_1\partial\theta_1^{\T}\})/2$ for $j=1,\ldots,k$, $c_2=E(\Gamma_{1i}^jB^{\T}\psi_i)$,
and $d_2=-E(\psi_i\psi_i^{\T}\mathcal{C}\psi_i)/2$.

\begin{assumption}\label{ass16} For $1\leq j\leq q$, the second and third derivatives of the penalty function satisfy $\max_{j}p''(|\theta_{10j}|)=O_p(1/n)$ and $p'''(|\theta_{10j}|)=0$, respectively.
\end{assumption}
A lot of penalty functions, for example, Lasso, SCAD (Fan and Li, 2001) and MCP (Zhang, 2010), satisfy Assumption \ref{ass16}.

\begin{corollary}\label{cor1} Under Theorem \ref{th5} and Assumption \ref{ass16}, we have
$${\rm Bias}(\hat{\theta}_1)=\{\mathcal{K}\mathcal{A}_1+\mathbb{H}\mathcal{A}_2\}/n+{\rm Bias}(\hat{\theta}_{1ET}),$$
where ${\rm Bias}(\hat{\theta}_{1ET})$ is derived from Theorem \ref{th4}.2 of Newey and Smith (2004).
\end{corollary}

When there is no penalty function in deriving ET estimator of $\theta_1$ or $\mathcal{W}=0$, we have ${\rm Bias}(\hat{\theta}_1)={\rm Bias}(\hat{\theta}_{1ET})$ because of $\mathcal{A}_1=\mathcal{A}_2=0$ for the considered case. Thus, ${\rm Bias}(\hat{\theta}_1)={\rm Bias}(\hat{\theta}_{1ET})$ for enough large nonzero parameters
for the SCAD penalty function because of $p'_\gamma(|\theta_{10j}|)=\gamma\neq 0$.

\subsection{Implementation}

\noindent Similar to Leng and Tang (2012), a nonlinear optimization procedure can be employed to maximize $\ell_p(\theta)$ given in Equation (\ref{PETL13}).
It is quite difficult to implement the nonlinear optimization procedure because of the nonconcave penalty function $p_\gamma(|\theta_j|)$ involved. To address the issue, we consider
the following local quadratic approximation to the penalty function (Fan and Li, 2001) at a fixed value $\theta_j^{(m)}$ of $\theta_j$:
$p_\gamma(|\theta_j|)\approx p_\gamma(|\theta_j^{(m)}|)+\frac{1}{2}\{p_\gamma '(|\theta_j^{(m)}|)/|\theta_j^{(m)}|\}\{\theta_j^2-(\theta_j^{(m)})^2\}$,
where $\theta_j^{(m)}$ is the estimated value of $\theta_j$ at the $m$th step and $\theta_j$ is the $j$th component of $\theta$. Thus, the nonlinear optimization algorithm given in Owen (2001) can be adopted to maximize Equation (\ref{PETL13}) based on the above local quadratic approximation of $p_\gamma(|\theta_j|)$.
Repeating the nonlinear optimization procedure until convergence yields the PET estimate $\hat\theta$ of $\theta$.

To obtain the PET estimate of $\theta$, it is also necessary to find a data-driven approach to select the penalty parameter $\gamma$.
To select an appropriate penalty parameter $\gamma$, we consider the following adjusted aBIC criterion:
${\rm aBIC}(\gamma)=-2\ell(\hat\theta_\gamma)+C_n\frac{\log(n)}{n}{\rm df}_\gamma$,
where $\ell(\theta)=\ell(\nu,\theta)$ is given in Equation (\ref{PETL22}), $\hat\theta_\gamma$ is the PET estimator of $\theta$ depending on the tuning parameter $\gamma$,
df$_\gamma$ is the number of nonzero components in $\theta$ representing the ``degrees of freedom" of the estimated unconditional moment models,
$C_n$ is a scaling factor diverging to infinity at a slow rate as $p\rightarrow \infty$. When $p$ is fixed, we set $C_n=1$,
otherwise we take $C_n={\rm max}\{\log\log p,1\}$ (Tang and Leng, 2010). The rigorous proof of the consistency of the
aBIC for the PET likelihood function is worth of further investigating.
 Also, the following GCV criterion (Fan and Li, 2001) can be used to select $\gamma$ in a linear regression model:
${\rm GCV}(\gamma)=n^{-1}||Y-X\theta_{\gamma}||^2/\{1-e(\gamma)/n\}^2$,
where $\theta_{\gamma}=\arg\max_{\theta}\{{\rm GCV}(\gamma)\}$, $e(\gamma)={\rm tr}\{P_X(\theta_{\gamma})\}$
and $P_X(\theta_{\gamma})=X\{X^{\T}X+n\mathcal{B}(\theta_{\gamma})\}^{-1}X^{\T}$ and $\mathcal{B}(\theta_{\gamma})={\rm diag}\{p_{\gamma}^{'}(|\theta_{1\gamma}|)/|\theta_{1\gamma}|,\ldots,p_{\gamma}^{'}(|\theta_{p\gamma}|)/|\theta_{p\gamma}|\}$.

\section{Simulation studies}
\renewcommand{\theequation}{3.\arabic{equation}}
\setcounter{equation}{0}

\noindent In this section, three simulation studies are conducted to investigate the finite sample performance of our proposed methodologies.

\noindent {\bf Experiment 1 (The Population Mean Vector)}.
In this experiment, we first generate independent and identically distributed random vector $Z_i\in \mathcal{R}^p$ whose components independently
follow the $\chi_1^2$ distribution (called as a correctly specified model `CM') or $\chi_{1.2}^2$ distribution (called as a misspecified model `MS', which is used to
investigate the robustness of our proposed PET procedure), and then set $X_i=\theta+\mathbb{R}^{1/2}(Z_i-{\bf 1}_p)$,
where the true value $\theta_0$ of parameter vector $\theta\in \mathcal{R}^p$ is set to be $\theta_0=(1,0.6,0.3,0,\ldots,0)^{\T}$,
and the true values of components in $\mathbb{R}=(\rho_{jl})$ are set to be $\rho_{jj}=1$ and $\rho_{jl}=0.3$ or $0.7$ for $j\not= l$, respectively,
which are used to investigate the performance of our presented PET estimator under different correlated structure, and ${\bf 1}_p$ is a $p\times 1$ vector
whose elements are one. Under model CM, we have $E(X_i)=\theta$ and ${\rm Var}(X_i)=\mathbb{R}$,
which implies that components of $X_i$ are not independent of each other. To illustrate our proposed methods, we consider the following estimating equations:
$g(X_i;\theta)=X_i-\theta$, which satisfy the unconditional moment restrictions: $E\{g(X_i;\theta_0)\}=0$ under model CM,
but $E\{g(X_i;\theta_0)\}\not= 0$ under model MS. Clearly, in the experiment, we have $r=p$.

We consider the following four combinations of dimensionality $p$ and sample size $n$: $(n,p)=(50,7)$, $(100,10)$, $(200,14)$ and $(500,19)$,
where $p$ is taken to be the integer of $8(3n)^{1/5.1}-14$, which is used to make comparison with Tang and Leng (2012).
For each of four combinations, $2000$ repetitions are conducted to investigate the accuracy of our proposed estimators in terms of root
mean square errors (RMS) and the performance of our proposed variable selection procedure.
For each replication, $\hat\theta$ (representing the `PET' estimator of $\theta$) is evaluated by maximizing $\ell_p(\theta)$ given in Equation (\ref{PETL13}) via the optimization procedure introduced in Section 2.7 with the initial value of $\theta$ taken to be $\theta^{(0)}=n^{-1}\sum_{i=1}^nZ_i$ for model CM
and $\theta^{(0)}=(1,0.6,0.3,0.01,\ldots,0.01)^{\T}$ for model MS.
Similar to Fan and Li (2001), we set a component of $\hat\theta$ to be zero whenever
it is less than some threshold value, such as 0.001, which is close to zero.

For comparison, we compute the sample mean estimator $\bar{X}=n^{-1}\sum_{i=1}^nX_i$ (denoted as `Mean' method),
the hard-threshold estimator $\hat{\theta}^{HT}_j=\bar{X}_jI(\bar{X}_j<\gamma_1)$ (denoted as `HT' method),
the soft-threshold estimator $\hat{\theta}^{ST}_j={\rm sign}(\bar{Z}_j)\{|\bar{Z}_j|-\gamma_2\}_{+}$ (denoted as `ST' method), and a quadratic-loss-based estimator
$\hat{\theta}_{QL}=\arg\min\limits_\theta\{(\bar{X}-\theta)^{\T}\mathcal{W}_n^{-1}(\bar{X}-\theta)+\gamma_3\sum_{i=1}^{p}|\theta_j|\}$ (denoted as `QL' method),
where $\bar{X}_j$ is the $j$th component of $\bar{X}$, $\gamma_1,\gamma_2$ and $\gamma_3$ are the tuning parameters, which can be obtained
by using a five-fold cross-validation method to minimize the squared predictive error for the mean vector, $\{t\}_+=t$ for $t>0$ and $0$ otherwise, and $\mathbb{W}_n=n^{-1}\sum_{i=1}^{n}(X_i-\bar{X})(X_i-\bar{X})^{\T}$. The aBIC criterion introduced in Section 2.7 is adopted to
select the tuning parameter $\gamma$ in the penalized function (\ref{PETL13}). We evaluate RMS values of nonzero components in $\theta_0$ for the above presented five estimators, and their corresponding average numbers of zero coefficients
that are correctly and incorrectly identified.

\vspace{1mm}
\hfill{$\underline{\overline{\sl Table~ 1~ about~ here}}$}
\vspace{2mm}

Results are presented in Table 1. Examination of Table 1 shows that (i) all the above mentioned four approaches (i.e., Mean, HT, ST and QL methods) to select zero coefficients yield a relatively small average number of false zero coefficients regardless of values of $p$, $n$ and $\rho_{jl}$;
(ii) the PET variable selection method and the hard-threshold variable selection method behave satisfactory
in the sense that their corresponding average numbers of correctly estimated zero components are quite close to the true number $p-3$ of zero components,
 whilst their corresponding average numbers of incorrectly estimated zero coefficients approach 0;
(iii) the RMS value of our proposed PET estimator is smaller than those of other estimators for our considered highly correlated data; (iv) increasing sample size or correlation among components can improve efficiency in terms of the RMS values and the average numbers of false or true zero coefficients.
These results demonstrate that our proposed PET method behaves better than others in terms of variable selection and parameter estimation, especially for highly correlated data, which indicates that our empirical results are consistent with those given in Theorem 1.

Also, to compare the performance of our-used adjusted BIC (aBIC) criterion with the traditional BIC (i.e.,
${\rm BIC}(\gamma)=-2\ell(\hat\theta_\gamma)+\frac{\log(n)}{n}{\rm df}_\gamma$) and AIC
(i.e., ${\rm AIC}(\gamma)=-2\ell(\hat\theta_\gamma)+\frac{2}{n}{\rm df}_\gamma$) criteria in growing dimensionality,
we computed the average model size (i.e., the average value of the number of non-zero coefficients, `AMS')
and the percentage of the correctly identified true model (`PCIM') for the above generated 2000 datasets.
Intuitively, a good model selection procedure should be a procedure whose AMS value is quite close to the true model size $q$ and PCIM value is close to 1.
Results are given in Table 2. Examination of Table 2 shows that (i) the AIC method fails to identify the
true model because of its over-fitting effect in large samples when the true model is of finite dimension; (ii)
the PCIM value of the aBIC method approaches $100\%$ and its AMS value is close to $q=3$ when $p$ is moderate or large;
(iii) the PCIM and AMS values of the BIC method increase as $p$ increases, and are close to $100\%$ and $q=3$ when $p$ is moderate or large, respectively, whilst
the aBIC method slightly outperforms the BIC method. In a word, the aBIC consistently outperforms the BIC and AIC criteria.

\vspace{1mm}
\hfill{$\underline{\overline{\sl Table~ 2~ about~ here}}$}
\vspace{2mm}

To investigate the performance of our proposed PET-likelihood-ratio-based confidence interval of parameter of interest,
we only evaluate the $95\%$ confidence interval of parameter $\theta_2$ for each of the above generated $2000$ data sets.
Table 3 presents the empirical frequencies of $\theta_2\notin R_{\alpha}$ for various true values of $\theta_2$.
Examination of Table 3 shows that (i) the frequency of $\theta_2\notin R_{\alpha}$ at the true value of $\theta_2=0.6$
is quite close to the pre-specified significant level $\alpha=0.05$ as $n$ (or $\rho_{jl}$) is large, for example, $n=500$ (or $\rho_{jl}=0.7$),
which is consistent with the conclusion given in Theorem \ref{th4}; (ii) power increases as $n$ or correlation
coefficient among components increases or
$\theta_2$ deviates more from the true value $0.6$ of $\theta_2$.
These observations show that our presented PET-likelihood-ratio-based test procedure performs well.

\vspace{1mm}
\hfill{$\underline{\overline{\sl Table~ 3~ about~ here}}$}
\vspace{2mm}

To investigate the robustness of the proposed PET procedure to misspecified unconditional moment models (i.e., Model MS), we first generate 2000 data sets from Model MS with the same parameter settings as given in the above simulation study, and then calculate our proposed PET estimates and the penalized empirical likelihood estimates (Leng and Tang, 2012) of $\theta$ and the corresponding average numbers of zero coefficients
that are correctly and incorrectly identified for our proposed variable selection procedure and Leng and Tang's (2012) procedure (denoted as `PEL' method)
for each of 2000 data sets based on estimating equations: $g(X_i,\theta)=X_i-\theta$, which do not satisfy unconditional moment restrictions: $E\{g(X_i,\theta)\}=0$
for $i=1,\ldots,n$ under Model CM. To wit, the fitted unconditional moment restrictions are misspecified.
The values of Bias, RMS and SD for $\theta$ with $(n,p)=(50,7)$ are presented in Table 4,
where `Bias' is the difference between the true value and
the mean of the estimates based on 2000 replications, and `SD' is the standard deviation of 2000 estimates.
From Table 4, we observe that (i) the PET estimates of parameters are robust to misspecified unconditional
moment models in terms of Bias, whilst the penalized empirical likelihood estimates of parameters are
sensitive to misspecified unconditional moment models in the sense that their corresponding Biases deviate from zero;
(ii) the values of `RMS' and `SD' are almost identical under our considered cases,
which indicates that the estimated standard deviation is rather reliable regardless of the PET method or the PEL method; (iii) the PET method behaves better
than the PEL method in terms of RMS values when unconditional moment models are misspecified; (iv) the accuracy of the PET estimator can be improved as the correlation among components of $X_i$ increases; (v) the PET variable selection procedure behaves better than the PEL variable selection method in the sense that the average number of correctly identifying
nonzero components for the PET method is quite close to the true number (i.e., $3$) of nonzero components even when unconditional moment models are misspecified.

\vspace{1mm}
\hfill{$\underline{\overline{\sl Table~ 4~ about~ here}}$}
\vspace{2mm}

\noindent {\bf Experiment 2 (Linear regression model)}. In the experiment, we consider the following linear regression model:
$Y_i=Z_i^{\T}\theta+\epsilon_i$ for $i=1,\ldots,n$, where $Z_i=(z_{i1},\ldots,z_{ip})^{\T}$ is assumed to follow a
multivariate normal distribution with zero mean and covariance matrix $\mathbb{R}=(\rho_{jl})$ with $\rho_{jl}=0.5^{|j-l|}$,
and $\epsilon_i$ follows the standard normal distribution $\mathcal{N}(0,1)$. The true value of $\theta\in\mathcal {R}^p$
is taken to be $\theta_0=(3,1.5,0,0,2,0,\ldots,0)$ including three nonzero components and $p-3$ zero components.
To illustrate our proposed approach to over-identified moment condition models, we introduce an instrumental variable $U_i=(u_{i1},\ldots,u_{ip})^{\T}$,
which are independently generated from $u_{ij}\stackrel{i.i.d.}{\sim}z_{ij}+\mathcal {N}(0,1)$ for $j=1,\ldots,p$, and consider the following unconditional moment restrictions:
$$
g(X_i;\theta)=(z_{i1}(Y_i-Z_i^{\T}\theta),\ldots,z_{ip}(Y_i-Z_i^{\T}\theta),u_{i1}(Y_i-Z_i^{\T}\theta),\ldots,u_{ip}(Y_i-Z_i^{\T}\theta))^{\T},
$$
which satisfy $E\{g(X_i;\theta)\}=0$, where $X_i=(Z_i^{\T},Y_i)^{\T}$ for $i=1,\ldots,n$. In this case, $r=2p$.

Similar to Experiment 1, $2000$ data sets $\{X_i: i=1,\ldots,n\}$ are independently generated from the above specified linear model to
evaluate RMS values of nonzero parameters in $\theta$ and the corresponding average numbers of nonzero coefficients correctly and incorrectly identified.
The tuning parameter $\gamma$ in the PET likelihood (\ref{PETL13}) is selected via the GCV criterion introduced in Section 2.7. For comparison, we evaluate the
RMS values of the least squares estimators of parameters in $\theta$ under the assumption that the true sparsity of the model is known.
Results corresponding to Table 1 are given in Table 5.

\vspace{1mm}
\hfill{$\underline{\overline{\sl Table~ 5~ about~ here}}$}
\vspace{2mm}

Examination of Table 5 shows that (i) the PET method performs well for variable selection in the sense that its corresponding average
number of zero coefficients correctly identified is quite close to $p-3$;
(ii) the RMS values of the PET estimators are slightly larger than those of the generalized least squares estimators
when $n$ is small, whilst their corresponding RMS values are almost identical when $n$ is large, for example, $n=500$;
(iii) the RMS values of two estimators decrease as $n$ increases.

\vspace{5mm}

\noindent {\bf Experiment 3 (Nonparametric structural equation model)}. In the experiment, we consider the following structural equation models:
\begin{equation}\label{ExEE1}
Y_i=\mathbb{Z}\omega_i+\epsilon_i,~~~\omega_i=\mathbb{U}\omega_i+\zeta_i,~~i=1,\ldots,n,
\end{equation}
where $Y_i$ is a $p_y\times 1$ vector of manifest variables, $\omega_i$ is a $q_\omega\times 1$ vector of latent variables, $\mathbb{Z}$ is a $p_y\times q_\omega$ factor loading matrix, $\mathbb{U}$ is a $q_\omega\times q_\omega$ coefficient matrix used to identify the correlation structure among latent variables, and it is assumed that measurement error $\epsilon_i$ is distributed as the multivariate normal distribution with zero mean and covariance $\Phi_\epsilon$, i.e., $\epsilon_i\sim \mathcal{N}(0,\Phi_\epsilon)$ in which $\Phi_\epsilon={\rm diag}(\phi_1,\ldots,\phi_{p_y})$, measurement error $\zeta_i$ follows the multivariate normal distribution $\mathcal{N}(0,\Psi_{\zeta})$ with $\Psi_{\zeta}$=diag$\{\tau_1,\cdots,\tau_{q_\omega}\}$, and $p_y=2q_\omega$.
The data set $\{Y_i: i=1,\ldots,n\}$ is generated from model (\ref{ExEE1}) with the following specifications of $\mathbb{Z}$, $\mathbb{U}$, $\Phi_\epsilon$ and $\Psi_\zeta$:
$$\mathbb{Z}=\left(
\begin{array}{ccccccc}
1& b_{21} & 0 & 0 &\cdots & 0& 0 \\
0&  0    & 1&b_{42}& \cdots&0& 0\\
\vdots&\vdots&\vdots&\vdots& \ddots &\vdots &\vdots \\
0&0 &0 & 0& \cdots &1&b_{p_yq_\omega}\\
\end{array}	
\right),~~
\mathbb{U}=\left(\begin{array}{cccccc}
0& \varphi_{12} &\cdots &\varphi_{1,q_\omega-1}& \varphi_{1,q_\omega}  \\
\varphi_{21}&  0   &\cdots & \varphi_{2,q_\omega-1}& \varphi_{2,q_\omega}\\
\vdots&\vdots & \ddots &\vdots &\vdots  \\
\varphi_{q_\omega,1}&\varphi_{q_\omega,2}  & \cdots  &\varphi_{q_\omega,q_\omega-1}  &0\\
\end{array}
\right)	
$$
in which $1$ and $0$ in $\mathbb{Z}$ and $\mathbb{U}$ are known parameters for model identification. The true values of $b_{2l,l}$, $\phi_j$ and $\tau_l$ are set to be $0.8$, $0.8$ and $0.8$ for $l=1,\ldots,q_\omega$ and $j=1,\ldots,p_y$, respectively; and the true value of $\varphi_{j_1,j_2}$ is taken to $0.8$ for $|j_1-j_2|=1$ and 0 otherwise.
Thus, there are $q_y^2+3q_\omega$ unknown parameters in $\theta$=\{$\mathbb{Z}$, $\mathbb{U}$, $\Phi_{\epsilon}$, $\Psi_{\zeta}$\}. To illustrate our proposed method, we consider the following unconditional moment restrictions: $g(X_i;\theta)={\rm vech}\{Y_iY_i^{\T}-\mathbb{O}(\theta)\}$, which satisfy $E\{g(Y_i;\theta_0)\}=0$, where  $\theta_0$ is the true value of $\theta$, $X_i=Y_i$, $\mathbb{O}(\theta)=\mathbb{Z}(I-\mathbb{U})^{-1}\Psi_{\zeta}(I-\mathbb{U})^{-1}\mathbb{Z}^{\T}+\Phi_{\epsilon}$ and ${\rm vech}(A)$ represents the half-vectorization of matrix $A$. Thus, the number of unconditional moment restrictions is $r=p_y(p_y+1)/2=q_\omega(2q_\omega+1)$.
Clearly, the above considered unconditional moment model is an over-identification case when $q_\omega>2$.
To solve the nonlinear optimization problem related to $\ell_p(\theta)$ given in Equation (\ref{PETL13}), an essential pre-requisite is that zero vector
is the interior point of the convex hull of $\{g(X_i;\theta): i=1,\ldots,n\}$.
Following Zhu et al. (2009), an adjusted PET likelihood $\ell_{ap}(\theta)$ can be used to evaluate the PET estimates of parameters in $\theta$. To this end, we define
$g_{n+1}(\theta)=g(X_{n+1};\theta)=-\frac{a}{n}\sum_{i=1}^ng(X_{i};\theta)$,
where $a=\max\{1, \log(n)/2\}$, and $X_{n+1}$ is introduced purely for notational simplicity. Thus,
the corresponding adjusted PET likelihood ratio function is given by
\begin{equation}\label{EXPLAP}
\ell_{ap}(\theta)=\log\left[\frac{1}{n+1}\sum\limits_{i=1}^{n+1}\exp\{\nu^{\T}(\theta)\mathrm{vech}(Y_iY_i^{\T}-\mathbb{O}(\theta))\}\right]-\sum\limits_{j=1}^{q_\omega}\sum\limits_{l=1,l\neq j}^{q_\omega}p_{\gamma}(|\varphi_{jl}|),
\end{equation}
which indicates that there are $2q_\omega^2+q_\omega=q_\omega(2q_\omega+1)$ moment restrictions but only $\varrho=q_\omega^2-q_\omega=q_\omega(q_\omega-1)$ parameters penalized.

Based on the above presented settings, we consider the following three combinations of the number of the penalized parameters $\varrho$ (i.e., $q_\omega$) and sample size $n$:
($n,q_\omega$) = (185, 3), (392, 4) and (919, 5), where $n$ is taken to be the integer of
$\frac{1}{3}((\varrho+32)/11)^{5.1}$ for $q_\omega$=3, 4 and 5. For comparison, we calculate results corresponding to the penalized empirical likelihood method.
To evaluate estimates of parameters in $\theta$, we choose the sieve space via the following procedure: (i) in the inner loop, given the current estimates $\hat{\mathbb{U}}$, $\hat{\Phi}_\epsilon$ and $\hat{\Psi}_\zeta$ of $\mathbb{U}$, $\Phi_\epsilon$ and $\Psi_\zeta$, we evaluate estimates $\hat{\mathbb{Z}}$ and $\hat\lambda(\hat{\mathbb{Z}})$ of $\mathbb{Z}$ and $\lambda$ based on the selected sieve space:  $\Theta_{s(n)}^{\mathbb{U}}\times\Theta_{s(n)}^{\Phi_{\epsilon}}\times\Theta_{s(n)}^{\Psi_{\zeta}}$ by applying the Newton-Raphison optimization algorithm to the adjusted PET likelihood $\hat{\ell}_{ap}(\mathbb{Z},\lambda(\mathbb{Z})|\hat{\mathbb{U}}, \hat{\Phi}_{\epsilon}, \hat{\Psi}_{\zeta})$, which is defined in Equation (\ref{EXPLAP}) with $\mathbb{U}$, $\Phi_\epsilon$ and $\Psi_\zeta$ replaced by $\hat{\mathbb{U}}$, $\hat{\Phi}_\epsilon$ and $\hat{\Psi}_\zeta$,  where  $\{\Theta_{s(n)}^{\mathbb{F}}\}_{s(n)=1}^{\infty}$ is a sieve space and a sequence of subsets of $\Theta$ for $\mathbb{F}=\mathbb{U}$, $\Phi_\epsilon$ and $\Psi_\zeta$; we can similarly evaluate estimates $\hat{\mathbb{U}}$ and $\hat\lambda(\hat{\mathbb{U}})$ of $\mathbb{U}$ and $\lambda$, estimates $\hat{\Phi}_\epsilon$ and $\hat\lambda(\hat{\Psi}_\epsilon)$ of $\Phi_\epsilon$ and $\lambda$, and estimates $\hat{\Psi}_\zeta$ and $\hat\lambda(\hat{\Psi}_\zeta)$ of $\Psi_\zeta$ and $\lambda$; (ii) in outer loop, we evaluate $(\hat{\mathbb{U}}, \hat{\lambda}(\hat{\mathbb{U}}))\to(\hat{\Phi}_{\epsilon}, \hat{\lambda}(\hat{\Phi}_{\epsilon}))\to (\hat{\Psi}_{\zeta}, \hat{\lambda}(\hat{\Psi}_{\zeta}))\to(\hat{\mathbb{Z}}, \hat{\lambda}(\hat{\mathbb{Z}}))\to\ldots$. Repeating the above procedure until the algorithm convergence yields the PET estimate $\hat\theta$ and penalized empirical likelihood estimate $\hat\theta_{ET}$ of $\theta$. The SD and RMS values of nonzero parameter estimators, the average numbers of correctly identified zero coefficients and incorrectly identified zero coefficients for 2000 replications are presented in Table 6. Inspection of Table 6 shows that (i) the SD and RMS values of two estimators decrease as $n$ increases; (ii) the average number of correctly identified zero coefficients approaches $q_\omega(q_\omega-1)$ for each of our considered two methods,
and the average number of incorrectly identified zero coefficients decreases as $n$ increases.

\vspace{1mm}
\hfill{$\underline{\overline{\sl Table~ 6~ about~ here}}$}
\vspace{2mm}

\section{An example}
\renewcommand{\theequation}{4.\arabic{equation}}
\setcounter{equation}{0}

\noindent In this section, an example taken from the Boston Housing Study is used to illustrate our proposed PET method in R package mlbench.
The data set, which has even been analyzed by Harrison and Rubinfeld (1978),  consists of $506$ observations on $14$ variables.
The main purpose of this study is to identify the effect of clean air on house prices.
Here, we take the logarithm of the median value (LMV) of owner occupied homes to be response variable ($y$),
and other 13 variables to be covariates. These covariates include
per capita crime rate by town (CRIM, $x_1$), proportion of residential land zoned for lots over 25,000 sq.ft (ZN, $x_2$),
proportion of non-retail business acres per town (INDUS, $x_3$), Charles river dummy variable which is $1$
if it is tract bounds river and 0 otherwise (CHAS, $x_4$),
nitric oxides concentration (parts per 10 million, NOX, $x_5$),
average number of rooms per dwelling (RM, $x_6$), proportion of owner-occupied units built prior to 1940 (AGE, $x_7$),
weighted distances to five Boston employment centers (DIS, $x_8$), index of accessibility to radial highways (RAD, $x_9$),
full-value property-tax rate per 10,000 (TAX, $x_{10}$), pupil-teacher ratio by town (OTRATIO, $x_{11}$),
$1000(bk-0.63)^2$ in which bk is the proportion of blacks by town (B, $x_{12}$), and proportion of population that has a lower status (LSTAT, $x_{13}$).

Following Harrison and Rubinfeld (1978), we consider the following linear model for the above introduced Boston Housing data set:
$y_i=x_i^{\T}\theta+\varepsilon_i$, where $x_i=(1,x_{i,1},\ldots,x_{i,91})^{\T}$ in which $x_{i,1},\ldots,x_{i,13}$
are the  above mentioned 13 covariates and $x_{i,14},\ldots,x_{i,91}$ are the interaction effects of any two covariates among 13 covariates,
$\theta=(\theta_0,\theta_1,\ldots,\theta_{91})^{\T}$, and $\varepsilon_i$ is the random error whose distribution is assumed to be unknown
but $E(\varepsilon_i)=0$. Generally, the least square method can be employed to estimate $\theta$. To illustrate our proposed method,
we consider the following unconditional moment restrictions:
$g(X_i;\theta)=x_i(y_i-x_i^{\T}\theta)$, where $X_i=(y_i,x_i^{\T})^{\T}$ for $i=1,\ldots,n$. 
Under the above given model assumption, we have $E\{g(X_i;\theta)\}=0$ for $i=1,\ldots,n$ with $n=506$. In this case, the number of moment restrictions is $r=92$.

The above presented PET method is used to evaluate estimate of $\theta=(\theta_0,\theta_1,\ldots,\theta_{91})^{\T}$ and
select significant covariates with the initial value of $\theta$ taken to be its least square estimate.
For comparison, we calculate the penalized empirical likelihood estimate of $\theta$.
The GCV method introduced in Section 2.7 is adopted to select the tuning parameter $\gamma$ in the PET likelihood function (\ref{PETL13}).
Similar to Fan and Li (2001), we set a component of $\hat\theta$ to be zero whenever its estimate is less than the threshold value 0.001.
Estimates of nonzero regression coefficients in $\theta$ identified by our proposed PET method are presented in Table 7.
Inspection of Table 7 shows that (i) covariates $x_1$, $x_3$, $x_4$, $x_6$, $x_7$, $x_9$, $x_{11}$, $x_{12}$ and $x_{13}$ are
the most significant covaraites in which variables CHAS, RM, AGE, RAD, TAX, OTRATIO and LSTAT have a positive effect on LMV of owner occupied homes,
whilst variables CRIM, INDUS have a negative effect on LMV of owner occupied homes;
(ii) interaction effects $x_1x_4$, $x_1x_6$, $x_1x_{10}$, $x_2x_4$, $x_3x_5$, $x_3x_{6}$, $x_7x_9$
 and $x_{10}x_{11}$ have a positive effect on LMV,
 whilst interaction effects $x_1x_3$, $x_1x_5$, $x_1x_9$, $x_3x_{11}$, $x_4x_5$, $x_4x_6$, $x_5x_7$, $x_6x_7$, $x_6x_{10}$, $x_6x_{11}$, $x_6x_{13}$,
 $x_7x_{11}$, $x_7x_{12}$, $x_7x_{13}$, $x_9x_{11}$, and $x_{10}x_{13}$ have a negative effect on LMV according to our presented PET method;
 (iii) although variable $x_5$ is not detected to be the significant covariate,
 interaction effects $x_1x_5$, $x_3x_5$ and $x_4x_5$ related to $x_5$ are identified to be the significant covariates;
 (iv) the estimated standard errors (SEs) of the proposed PET estimators are smaller than those of the penalized empirical likelihood estimators;
 (v) the PET estimators have shorter confidence intervals than the penalized empirical likelihood estimators.

\hfill{$\underline{\overline{\sl Table~ 7~ about~ here}}$}
\vspace{-3mm}

\section{Discussion}

\noindent This paper presents a PET likelihood procedure for variable selection and parameter estimation in unconditional moment models with a diverging number of parameters
in the presence of correlation among variables and model misspecification.
We show that the PET likelihood possesses some properties analogous to the penalized likelihood such as the oracle properties,
and the PET approach is robust to model misspecification like the exponentially tilted method.
Under some regularity conditions, we show that the constrained PET likelihood ratio statistic for testing contrast hypothesis is asymptotically distributed as the central chi-squared distribution.
Simulation studies are conducted to investigate the finite sample performance of the proposed methodologies, and an example from the Boston Housing Study is used to illustrate our proposed methodologies. Empirical results evidence that the penalized empirical likelihood method leads to an inappropriate conclusion when unconditional moment models are misspecified, but our proposed PET method leads to a desirable conclusion even if unconditional moment models are misspecified.

The proposed PET method in this paper is developed for the completely observed data. In many applications such as economics, data are often subject to missingness due to nonresponse or dropout of participants. In this case, it is interesting to consider the PET likelihood inference on growing dimensional unconditional moment models with missing data. It is well-known that influence analysis is an important step in data analysis. When there are influential observations in a data set, an important issue is how to identify these influential observations for our considered growing dimensional unconditional moment models with missing data. We are working on the two topics.

The developed PET theories in this paper are derived on the basis of the assumption: $p/r<1$. There are theories on generalized EL estimators that allow $p/r>1$ (e.g., see Shi, 2014). Hence, it is interesting to extend the proposed PET theories to the case: $p>r$. To wit, it is interesting to develop some theories and methods to simultaneously select parameter and moment restrictions as done in Cancer, Han and Lee (2016). Also, as a referee pointed out that it is difficult to imagine that the ``pseudo-true'' value $\theta^*$ of parameter vector $\theta$ is also fixed as sample size $n$ varies because the number of moments $r$ increases with sample size $n$. The issue can be incorporated to moment selection problems as done in Cancer, Han and Lee (2016).

\section*{Acknowledgements}

\noindent The authors are grateful to the Editor, an Associate Editor and three referees for their valuable suggestions and comments
that greatly improved the manuscript. The research was supported by grants from the National Natural Science Foundation of China (Grant No.: 11165016).

\vspace{1cm}

\section*{Appendix: Proofs of Theorems}
\renewcommand{\theequation}{A.\arabic{equation}}
\setcounter{equation}{0}

\begin{lemma}\label{lem1}
Under Assumption \ref{ass1}(iii), if $\nu\in\mathcal{V}_n=\{\nu\in\mathcal{R}^{r}:||\nu||\leq \pi_n\}$ and
$\lambda\in\Lambda_n=\{\lambda\in\mathcal{R}^{k}:||\lambda||\leq \rho_{n}\}$, where $\pi_n=o_p(r^{-1/2}n^{-1/\delta})$ and $\rho_{n}=o_p(k^{-1/2}n^{-1/\delta})$, we have
$\max_{1\leq i\leq n}\sup_{\theta\in\Theta_{s(n)}}|\nu^{\T}g(X_i;\theta)|=o_p(1)$ and $\max_{1\leq i\leq n}\sup_{\theta_1\in\Theta_1}|\lambda^{\T}\psi(Z_i;\theta_1)|=o_p(1)$. Also, {\rm w.p.a.1}, $\mathcal{V}_n\subseteq\widehat{\mathcal{V}}_n(\theta)$ for all $\theta\in\Theta_{s(n)}$, and $\Lambda_n\subseteq\widehat{\Lambda}_n(\theta_1)$ for all $\theta_1\in\Theta_1$.
\end{lemma}

\noindent \emph{\bf Proof}.
Assumption \ref{ass1}(iii) implies that $\max_{1\leq i\leq n}\sup_{\theta\in\Theta_{s(n)}}||g(X_i;\theta)||=O_p(n^{1/\delta}r^{1/2})$,
also w.p.a.1 $\max_{1\leq i\leq n}\sup_{\theta_1\in\Theta_1}||\psi(Z_i;\theta_1)||=O_p(n^{1/\delta}k^{1/2})$. Then, we have
\[
\begin{array}{llll}
\max_{1\leq i\leq n}\sup_{\theta\in\Theta_{s(n)}}|\nu^{\T}g(X_i;\theta)|
\leq \pi_n\max_{1\leq i\leq n}\sup_{\theta\in\Theta_{s(n)}}||g(X_i;\theta)||=o_p(1),\\
\max_{1\leq i\leq n}\sup_{\theta_1\in\Theta_1}|\lambda^{\T}\psi(Z_i;\theta_1)|
\leq \rho_n\max_{1\leq i\leq n}\sup_{\theta_1\in\Theta_1}||\psi(Z_i;\theta_1)||=o_p(1).
\end{array}
\]
So, w.p.a.1 $\nu^{\T}g(X_i;\theta)\in \mathcal{E}$ for all $\theta\in\Theta_{s(n)} $ and $||\nu||\leq \pi_n$,
and $\lambda^{\T}\psi(Z_i;\theta_1)\in \mathcal{E}$ for all $\theta_1\in\Theta_1 $ and $||\lambda||\leq \rho_{n}$.

\begin{lemma}\label{lem6} Suppose that Assumptions \ref{ass1} and \ref{ass3} hold. Then, for unconditional moment restrictions $g(X;\theta)$ uniformly in $\theta\in$ $\mathcal{D}_n$,
we have

{\rm (i)} $||\frac{1}{n}\sum\limits_{i=1}^ng(X_i;\Pi_n\theta_0)-E(g(X_i;\Pi_n\theta_0))||=O_p(\sqrt{r/n})$;

{\rm (ii)} $ ||\frac{1}{n}\sum\limits_{i=1}^n\frac{\partial g(X_i;\Pi_n\theta_0)}{\partial\theta}-\Gamma(\Pi_n\theta_0)||=O_p(\sqrt{rp/n})$;

{\rm (iii)} $ ||\frac{1}{n}\sum\limits_{i=1}^ng(X_i;\Pi_n\theta_0)g(X_i;\Pi_n\theta_0)^{\T}-\Sigma(\Pi_n\theta_0)||=O_p(r/\sqrt{n})$;

{\rm (iv)} $ ||\frac{1}{n}\sum\limits_{i=1}^ng(X_i;\Pi_n\theta_0)-E(g(X_i;\Pi_n\theta_0))-\frac{1}{n}\sum\limits_{i=1}^ng(X_i;\theta_0)+E(g(X_i;\theta_0))||=o_p(\sqrt{r/n})$;

{\rm (v)} $||\frac{1}{n}\sum\limits_{i=1}^ng(X_i;\Pi_n\theta_0)g(X_i;\Pi_n\theta_0)^{\T}-\Sigma(\Pi_n\theta_0)-\frac{1}{n}\sum\limits_{i=1}^ng(X_i;\theta_0)g(X_i;\theta_0)^{\T}+
\Sigma(\theta_0)||=o_p(r/\sqrt{n})$.
\end{lemma}

\noindent {\bf Proof}. Since the collection  of components of unconditional moment restrictions $g(X_i;\theta)$ is P-Donsker class (Kosorok, 2008), thus (i) and (ii) hold.
Again, the collection of any products of components of unconditional moment restrictions $g(X_i;\theta)$ is also P-Glivenko-Cantelli (P-G-C) class (Kosorok, 2008), thus (iii) holds.
In fact, (i)-(iii) are the standard uniform consistency results, which are obtained by the law of large numbers. While (iv) and (v) are Bahadur type modulus of continuity results.

\begin{lemma}\label{lem7}
Suppose that Assumptions \ref{ass1}(iii) and \ref{ass3}(i) hold and $r^2=o_p(n)$. Then, for any $\theta\in D_n$, $\nu(\theta)=\arg\min_{\nu \in\widehat{\mathcal{V}}_n(\theta)}\ell(\nu,\theta)$ exists, and $\nu(\theta)=-\{\Sigma(\theta)\}^{-1}\bar{g}(\theta)+o_p(\sqrt{r/n})$, where $\bar{g}(\theta)=n^{-1}\sum_{i=1}^ng(X_i;\theta)$.
\end{lemma}

\noindent {\bf Proof}. Taking $\pi_{n}=o_p(r^{-1/2}n^{-1/\delta})$ yields $\sqrt{r/n}=o_p(\pi_{n})$ because of $r^{2}n^{2/\delta-1}=o_p(1)$. Let $\bar{\nu}=\arg\inf_{\nu\in\mathcal{V}_n}\ell(\nu, \theta)$, where $\mathcal{V}_n$ is defined in Lemma \ref{lem1}. Then, we have $\max_{1\leq i\leq n}\sup_{\theta\in\Theta_{s(n)}}|\nu^{\T}g(X_i;\theta)|=o_p(1)$. Set $v_i=\nu^{\T}g(X_i;\theta)$, $\partial\ell(v_i)/\partial v_i=\frac{1}{n}\sum_{i=1}^n\rho_1(v_i)$, $\partial^2\ell(v_i)/\partial v_i^2=\frac{1}{n}\sum_{i=1}^n\rho_2(v_i)$. For any $\dot{\nu}$ on the line joining $\bar{\nu}$ and 0, it follows from Lemma \ref{lem1} and $\rho_2(0)=1-1/n$ that w.p.a.1 $\rho_2(\dot{\nu}^{\T}g(X_i;\theta))\ge C$. Thus, by Assumption \ref{ass3}(i) and the Taylor expansion at $\nu=0$, we obtain
\[
\begin{array}{llll}
0=\ell(0, \theta)\ge\ell(\bar{\nu}, \theta)&=&\rho_1(0)\bar{\nu}^{\T}\bar{g}(\theta)+\frac{1}{2}\bar{\nu}^{\T}\left\{\frac{1}{n}\sum\limits_{i=1}^n\rho_2(\dot{\nu}^{\T}g(X_i;\theta))g(X_i;\theta)g^{\T}(Z_i;\theta)\right\}
\bar{\nu}\\
&\ge&-||\bar{\nu}||||\bar{g}(\theta)||+C||\bar{\nu}||^2,
\end{array}
\]
which leads to $||\bar{\nu}||\leq ||\bar{g}(\theta)||$ w.p.a.1. By Lemma \ref{lem6}(i) and (iv), we have $||\bar{g}(\theta)||=O_p(\sqrt{r/n})$. Combining the above equations yields  $||\bar{\nu}||=O_p(\sqrt{r/n})=o_p(\pi_{n})$. Therefore, w.p.a.1
$\bar{\nu}\in{\rm int}(\mathcal{V}_n)$, which indicates $\partial\ell(\bar{\nu},\theta)/\partial\nu=0$. By Lemma \ref{lem1} and the convexity of $\ell(\nu,\theta)$ and $\widehat{\mathcal{V}}_n(\theta)$, it follows that $\bar{\nu}=\nu(\theta)$ and $\arg\inf_{\nu\in\widehat{\mathcal{V}}_n(\theta)}\ell(\nu, \theta)$ exists.

Taking the first-order partial derivative of $\ell(\nu,\theta)$ with respect to $\nu$ yields
$$\partial\ell(\nu,\theta)/\partial\nu=\frac{1}{n}\sum\limits_{i=1}^n\frac{\exp\{\nu^{\T}g(X_i;\theta)\}}{1/n\sum_{j=1}^n\exp\{\nu^{\T}g(X_j;\theta)\}}g(X_i;\theta)=0.$$
Since $\max_{1\leq i\leq n}\sup_{\theta\in\Theta_{s(n)}}|\nu^{\T}g(X_i;\theta)|=o_p(1)$, the above equation is equivalent to $$\sum_{i=1}^n\exp\{\nu^{\T}g(X_i;\theta)\}g(X_i;\theta)=\sum\limits_{i=1}^n\{1+\nu^{\T}g(X_i;\theta)(1+o_p(1)\}g(X_i;\theta)=0.$$
Therefore, it follows from Lemma \ref{lem6}(iii) and $r^2=o(n)$ that $\nu(\theta)=-\Sigma^{-1}(\theta)\bar{g}(\theta)+o_p(\sqrt{r/n})$.

\begin{lemma}\label{lem8} If Assumptions \ref{ass1}(iii) and \ref{ass3}(i) hold and $r^2=o(n)$, $rp=o(n)$, for any sequence of sets $\{\theta:||\theta-\theta_0||=O(\sqrt{r/n})\}$, we have
$\ell(\theta)=-(\theta-\theta_0)^{\T}\Gamma\Sigma^{-1}\Gamma^{\T}(\theta-\theta_0)+2(\theta-\theta_0)^{\T}\Gamma\Sigma^{-1}\bar{g}(\theta_0)
-\bar{g}^{\T}(\theta_0)\Sigma^{-1}\bar{g}(\theta_0)+o_p(r/n)$, and
$\hat{\theta}_{ET}-\theta_0=(\Gamma\Sigma^{-1}\Gamma^{\T})^{-1}\Gamma\Sigma^{-1}\bar{g}(\theta_0)$,
where $\hat{\theta}_{ET}$ is the ET estimator of $\theta$.
\end{lemma}
\noindent
\emph{\bf Proof}. It follows from Lemma \ref{lem7} that the ET likelihood $\ell(\theta)$ can be written as
$\ell(\theta)=-\bar{g}^{\T}(\theta)\Sigma^{-1}(\theta)\bar{g}(\theta)+o_p(r/n)$. For any $\theta\in \mathcal{D}_n$,
considering the expansion of $\bar{g}(\theta)$ at $\theta_0$ and using Lemmas \ref{lem6}(i), \ref{lem6}(ii) and \ref{lem6}(v),
we obtain $\bar{g}(\theta)=\bar{g}(\theta_0)+\Gamma^{\T}(\theta-\theta_0)+o_p(\sqrt{r/n})$. Then,
it follows from Lemmas \ref{lem6}(iii) and \ref{lem6}(v) that Lemma \ref{lem8} holds.

\vspace{3mm}

\noindent
\emph{\bf Proof of Theorem \ref{th6}}. From Lemma \ref{lem8}, we obtain
$\ell(\theta)=-(\theta-\theta_0)^{\T}\Gamma\Sigma^{-1}\Gamma^{\T}(\theta-\theta_0)+2(\theta-\theta_0)^{\T}\Gamma\Sigma^{-1}\bar{g}(\theta_0)
-g^{\T}(\theta_0)\Sigma^{-1}\bar{g}(\theta_0)+o_p(r/n)$,
and $(\Gamma\Sigma^{-1}\Gamma^{\T})(\hat{\theta}_{ET}-\theta_0)=\Gamma\Sigma^{-1}\bar{g}(\theta_0)$. Combining the above equations yields
\[
\begin{array}{llll}
\ell(\theta)&=&-(\theta-\theta_0)^{\T}\Gamma\Sigma^{-1}\Gamma^{\T}(\theta-\theta_0)+2(\theta-\theta_0)^{\T}\Gamma\Sigma^{-1}\Gamma^{\T}(\theta-\theta_0)\\[1mm]
&&-\bar{g}^{\T}(\theta_0)\Sigma^{-1}\bar{g}(\theta_0)+o_p(r/n)\\[1mm]
&=&-(\theta-\theta_0)^{\T}\Gamma\Sigma^{-1}\Gamma^{\T}(\theta-2\hat{\theta}_{ET}+\theta_0)-\bar{g}^{\T}(\theta_0)\Sigma^{-1}\bar{g}(\theta_0)+o_p(r/n)\\[1mm]
&=&-(\theta-\hat{\theta}_{ET})^{\T}\Gamma\Sigma^{-1}\Gamma^{\T}(\theta-\hat{\theta}_{ET})+o_p(r/n).
\end{array}
\]
It follows from the local quadratic approximation to the penalty function given in Section 2.7 that
$p_\gamma(\theta)\propto(\theta-\theta_{0})^{\T}J_{0}(\theta-\theta_{0})+o_p(r/n)$ for any $\theta\in \mathcal{D}_n$. Then, we have
\[
\begin{array}{llll}
\ell_p(\theta)&=&\ell(\theta)-p_\gamma(\theta)\\
&=&-(\theta-\hat{\theta}_{ET})^{\T}\Gamma\Sigma^{-1}\Gamma^{\T}(\theta-\hat{\theta}_{ET})-(\theta-\theta_{0})^{\T}J_{0}(\theta-\theta_{0})+R_n\\
&\propto& -(\theta-\hat{\theta})^{\T}\mathfrak{J}(\theta-\hat{\theta})+C_n+R_n,
\end{array}
\]
where $\mathfrak{J}=J_{0}+\Gamma\Sigma^{-1}\Gamma^{\T}$, $\hat{\theta}=\mathfrak{J}^{-1}(J_{0}\theta_{0}+\Gamma\Sigma^{-1}\Gamma^{\T}\hat{\theta}_{ET})$, $R_n=o_p(r/n)$,
and $C_n=-\theta_{0}^{\T}J_{0}\theta_{0}-\hat{\theta}_{ET}\Gamma\Sigma^{-1}\Gamma^{\T}\hat{\theta}_{ET}+\hat{\theta}^{\T}\mathfrak{J}\hat{\theta}$ is some constant that dose not depend on $\theta$.

It follows from Equation (\ref{selec2}) and Lemma \ref{lem8} that
$\hat{\theta}=(\Gamma\Sigma^{-1}\Gamma^{\T}+J_{0})^{-1}(\Gamma\Sigma^{-1}\Gamma^{\T}\hat{\theta}_{ET})$ and
$\hat{\theta}_{ET}=(\Gamma\Sigma^{-1}\Gamma^{\T})^{-1}\Gamma\Sigma^{-1}\bar{g}(\theta_0)+\theta_0$. Then, we obtain
$\hat{\theta}-\theta_0=-(\Gamma\Sigma^{-1}\Gamma^{\T}+J_{0})^{-1}J_{0}\theta_0+
(\Gamma\Sigma^{-1}\Gamma^{\T}+J_{0})^{-1}\Gamma\Sigma^{-1}\bar{g}(\theta_0)$, which yields
$E(\hat{\theta}-\theta_0)=-(\Gamma\Sigma^{-1}\Gamma^{\T}+J_{0})^{-1}J_{0}\theta_0=m$ and
var$(\hat{\theta}-\theta_0)=(\Gamma\Sigma^{-1}\Gamma^{\T}+J_{0})^{-1}\Gamma\Sigma^{-1}\Gamma^{\T}(\Gamma\Sigma^{-1}\Gamma^{\T}+J_{0})^{-1}=V$.
Let $m_j$ be the $j$th component of $m$, $V_{jj}$ be the $j$th diagonal element of $V$, and $\mathcal{L}$ be standard normal random variable. Denote $m^*=\max_j\{|m_j|\}$ and $V^*=\mathbb{E}_{\rm max}\{V\}$.
Then, we have
\begin{equation}{\label{tail}}
           \begin{aligned}
\sum_{j=1}^p{\rm Pr}(|\hat{\theta}_{j}-\theta_{0j}|\ge \gamma)&=\sum_{j=1}^p{\rm Pr}\left((|\hat{\theta}_{j}-\theta_{0j}|-|m_j|)/\sqrt{V^*}\ge (\gamma-|m_j|)/\sqrt{V^*}\right)\\
&\leq\sum_{j=1}^p{\rm Pr}\left((|\hat{\theta}_{j}-\theta_{0j}-m_j|)/\sqrt{V_{jj}}\ge (\gamma-m^*)/\sqrt{V^*}\right)\\
&\leq p{\rm Pr}\left(|\mathcal{L}|\ge (\gamma-m^*)/\sqrt{V^*}\right)\\
&\leq \frac{2 p\sqrt{V^*}}{\gamma-m^*}\exp\left\{-\frac{(\gamma-m^*)^2}{2V^*}\right\}.
 \end{aligned}
\end{equation}
Note that it is necessary to assume $\gamma>m^*$ because of ${\rm Pr}(|\hat{\theta}_{j}-\theta_{0j}-m_j|/\sqrt{V_{jj}}\ge (\gamma-m_j^*)/\sqrt{V^*})=1$ when $\gamma\leq m^*$, which indicates that the bound of $\sum_{j=1}^p{\rm Pr}(|\hat{\theta}_{j}-\theta_{0j}|\ge \gamma)$ is $p$ and we can not obtain the selection consistency. Although $m^*$ should not be larger than $\gamma$, we want to gain the smaller tail probability by increasing $m^*$.
The same situation appears for $V^*$. Hence, it is necessary to get the bounds of the maximum bias $m^*$ and the maximum variance $V^*$.

Note that $V=\{(\Gamma\Sigma^{-1}\Gamma^{\T}+J_{0})^{-1}-J_{0}(\Gamma\Sigma^{-1}\Gamma^{\T}+J_{0})^{-2}\}$. Thus, we have
$V^*=\mathbb{E}_{\max}\{(\Gamma\Sigma^{-1}\Gamma^{\T}+J_{0})^{-1}-J_{0}(\Gamma\Sigma^{-1}\Gamma^{\T}+J_{0})^{-2}\}
\leq\mathbb{E}_{\max}\{(\Gamma\Sigma^{-1}\Gamma^{\T}+J_{0})^{-1}\}\leq \mathcal{E}_n^{-1}$,
where $\mathcal{E}_n$ is the smallest eigenvalue of matrix $\Gamma\Sigma^{-1}\Gamma^{\T}+\phi_nI$ in which $\phi_n$ is some diagonal element of  $J_{0}$. If $\mathcal{E}_n\ge \phi_n$, thus we have $V^*\leq \phi_n^{-1}$, which indicates that we should select the larger $\phi_n=\min_{j\in\mathbb{J}^c}p_{\gamma}^{'}(|(\Pi_n\theta_{0})_j|)/\gamma$, otherwise, $V^*$ will be magnified to $1$ that would not have selection consistency. Hence, 
we choose
$V^*\leq \gamma/\min_{j\in\mathbb{J}^c}p_{\gamma}^{'}(|(\Pi_n\theta_{0})_j|)$.

To obtain the bound of the bias $m^*$, we note that
$m^*=||(\Gamma\Sigma^{-1}\Gamma^{\T}+J_{0})^{-1}J_{0}\theta_0||_{\infty}=
||(\Gamma\Sigma^{-1}\Gamma^{\T}+\phi_nI)^{-1}\phi_n\theta_0||_{\infty}=||Q^*\mathbb{D}_n^*Q^{*\T}\theta_0||_{\infty}$,
where $Q^*$ is orthogonal matrix derived from the eigenvalue decomposition of $\Gamma\Sigma^{-1}\Gamma^{\T}$,
$\mathbb{D}_n^*$=Diag$(\frac{\phi_n}{d_1^*+\phi_n},\ldots,\frac{\phi_n}{d_p^*+\phi_n})$,
and $d_1^*\geq d_2^*\geq \cdots\geq d_p^*$ are $p$ eigenvalues of matrix $\Gamma\Sigma^{-1}\Gamma^{\T}$. Further, we obtain
\[
\begin{array}{llll}
||Q^*\mathbb{D}_n^*Q^{*\T}\theta_0||_{\infty}&=&||Q^*\mathbb{D}_n^*Q^{*\T}Q_{1}^*\eta_1+Q^*\mathbb{D}_nQ^{*\T}Q_{2}^*\eta_2||_{\infty}\\
&\leq&||Q^*\mathbb{D}_n^*Q^{*\T}Q_{1}^*\eta_1||_{\infty}+||Q^*\mathbb{D}_n^*Q^{*\T}Q_{2}^*\eta_2||_{\infty}\\
&=&||Q^*\mathbb{D}_n^*(\eta_1^{\T},0^{\T})^{\T}||_{\infty}+||Q^*\mathbb{D}_n^*(0^{\T},\eta_2^{\T})^{\T}||_{\infty}\\
&=&||Q^*\mathbb{D}_n^*(\eta_1^{\T},0^{\T})^{\T}||_{\infty}+||Q_{2}^*\eta_2||_{\infty}
\end{array}
\]
By Assumption \ref{ass4}(iv), we have $||Q_{2}^*\eta_2||_{\infty}=O_p(\max_{j\in\mathbb{J}}p_{\gamma}^{'}(|(\Pi_n\theta_{0})_j|)\sqrt{q}/\gamma)$. For the first term, we have
\[
\begin{array}{llll}
||Q^*\mathbb{D}_n^*(\eta_1^{\T},0^{\T})^{\T}||_{\infty}&\leq& ||Q^*\mathbb{D}_n^*(\eta_1^{\T},0^{\T})^{\T}||=||\mathbb{D}_n^*(\eta_1^{\T},0^{\T})^{\T}||\\
&\leq& \frac{\phi_n}{d_{t}^*+\phi_n}||\eta_1||=\frac{\phi_n}{d_{t}^*+\phi_n}||Q_{1}^*\eta_1||\leq\frac{\phi_n}{d_{t}^*+\phi_n}||\theta_0||.
\end{array}
\]
Because all components of $\theta_0$ are finite and there are $q$ nonzero components in $\theta_0$,
we should choose $\phi_n=\max_{j\in\mathbb{J}}p_{\gamma}^{'}(|(\Pi_n\theta_{0})_j|)/\gamma$, otherwise,
$\phi_n/(d_{q}^*+\phi_n)$ will be magnified to $1$ that would not have selection consistency. Thus,
 we should choose $m^*\leq \max_{j\in\mathbb{J}}p_{\gamma}^{'}(|(\Pi_n\theta_{0})_j|)\sqrt{q}/\gamma$.
It follows from Assumption \ref{ass4}(ii) and \ref{ass4}(iii) that $\gamma>2m^*$ should be selected for sufficiently large $n$,
and Equation (\ref{tail}) can be rewritten as
$$
           \begin{aligned}
{\rm Pr}(\hat{\mathbb{J}}\neq \mathbb{J})&\leq2\frac{p\sqrt{V^*}}{\gamma/2}\exp\left\{-\frac{(\gamma/2)^2}{2V^*}\right\}\\
&\leq2\frac{p}{\sqrt{\gamma\min_{j\in\mathbb{J}^c}p_{\gamma}^{'}(|(\Pi_n\theta_{0})_j|)}}
\exp\left\{-\frac{\gamma\min_{j\in\mathbb{J}^c}p_{\gamma}^{'}(|(\Pi_n\theta_{0})_j|)}{8}\right\}\rightarrow 0,
 \end{aligned}
$$
which implies that Theorem \ref{th7} holds.

\begin{lemma}\label{lem2}
Suppose that Assumptions \ref{ass1} and \ref{ass3}(i) hold. Then, we have
$\lambda(\theta_{10})=\arg\inf_{\lambda\in\widehat{\Lambda}_n(\theta_{10})}\bar{\ell}(\lambda, \theta_{10})$ {\rm w.p.a.1.}, $\inf_{\lambda\in\widehat{\Lambda}_n(\theta_{10})}\bar{\ell}(\lambda,\theta_{10})= O_p(k/n)$, $\lambda(\theta_{10})=O_p(\sqrt{k/n})$ and $\hat{\nu}(\theta_{0})=O_p(\sqrt{r/n})$.
\end{lemma}

\noindent \emph{\bf Proof}.
Taking $\rho_{n}=o_p(k^{-1/2}n^{-1/\delta})$ yields $\sqrt{k/n}=o_p(\rho_{n})$ because of $k^{2}n^{2/\delta-1}=o_p(1)$, which is derived from Assumption \ref{ass1}(iii): $r^{2}n^{2/\delta-1}=o_p(1)$. Let $\bar{\lambda}=\arg\inf_{\lambda\in\Lambda_n}\bar{\ell}(\lambda, \theta_{10})$, where $\Lambda_n$ is defined in Lemma \ref{lem1}. Thus, we have $\max_{1\leq i\leq n}\sup_{\theta_1\in\Theta_1}|\lambda^{\T}\psi(Z_i;\theta_1)|=o_p(1)$. Denote $\tilde{v}_i=\lambda^{\T}\psi(Z_i;\theta_1)$,
$\partial\bar{\ell}(\tilde{v}_i)/\partial \tilde{v}_i=\frac{1}{n}\sum_{i=1}^n\tilde\rho_1(\tilde{v}_i)$ and
$\partial^2\bar{\ell}(\tilde{v}_i)/\partial \tilde{v}_i^2=\frac{1}{n}\sum_{i=1}^n\tilde\rho_2(\tilde{v}_i)$. For any $\lambda_d$ lying in the joining line between $\bar{\lambda}$ and the original point $0$, it follows from Lemma \ref{lem1} and $\tilde\rho_2(0)=1-1/n$ that w.p.a.1 $\tilde\rho_2(\lambda_d^{\T}\psi(Z_i;\theta_{10}))\ge C$.
Then, by Assumption \ref{ass3}(i), we obtain $a_0\leq \sup_{\theta_1\in \Theta_1}\mathbb{E}\{\frac{1}{n}\sum_{i=1}^n\psi(Z_i;\theta_1)\psi^{\T}(Z_i;\theta_1)\}\leq b_0<\infty$ w.p.a.1. The Taylor expansion at $\lambda=0$ with Lagrange remainder leads to
\begin{equation}\label{EQUA1}
\begin{array}{llll}
0=\bar{\ell}(0, \theta_{10})&\ge&\bar{\ell}(\bar{\lambda}, \theta_{10})=\tilde\rho_1(0)\bar{\lambda}^{\T}\bar{\psi}(\theta_{10})\\
&&+\frac{1}{2}\bar{\lambda}^{\T}
\left\{\frac{1}{n}\sum\limits_{i=1}^n\tilde\rho_2(\lambda_d^{\T}\psi(Z_i;\theta_{10}))\psi(Z_i;\theta_{10})\psi^{\T}(Z_i;\theta_{10})\right\}\bar{\lambda}\\
&\ge&-||\bar{\lambda}||||\bar{\psi}(\theta_{10})||+C||\bar{\lambda}||^2,
\end{array}
\end{equation}
which yields $||\bar{\lambda}||\leq ||\bar{\psi}(\theta_{10})||$ w.p.a.1. Therefore, $||\bar{\psi}(\theta_{10})||=O_p(\sqrt{k/n})$,
which yields  $||\bar{\lambda}||=O_p(\sqrt{k/n})=o_p(\rho_n)$. Thus, w.p.a.1 $\bar{\lambda}\in{\rm int}(\Lambda_n)$,
which yields $\partial\bar{\ell}(\bar{\lambda},\theta_{10})/\partial\lambda=0$.
By Lemma \ref{lem1} and the convexity of $\bar{\ell}(\lambda,\theta_{10})$ and $\widehat{\Lambda}_n(\theta_{10})$, it follows that $\bar{\lambda}=\lambda(\theta_{10})$, and $\arg\inf_{\lambda\in\widehat{\Lambda}_n(\theta_{10})}\bar{\ell}(\lambda, \theta_{10})$ exists.
By the last inequality given in Equation (\ref{EQUA1}), we obtain $\inf_{\lambda\in\widehat{\Lambda}_n(\theta_{10})}\bar{\ell}(\lambda, \theta_{10})= O_p(k/n)$.

By the same argument as done above, it follows from Lemma \ref{lem1} and Assumption \ref{ass3}(i) that $\hat{\nu}(\theta_{0})=O_p(\sqrt{r/n})$.

\vspace{2mm}

\begin{lemma}\label{lem3}
If Assumptions \ref{ass1}, \ref{ass3}(i) and \ref{ass7} hold, we have $||\bar{\psi}(\hat{\theta}_1)||=O_p(\sqrt{k/n})$ and $||\lambda(\hat{\theta}_1)||=O_p(\sqrt{k/n})$.
\end{lemma}
\noindent \emph{\bf Proof}.
For $\rho_{n}$ defined in Lemma \ref{lem1}, let $\tilde{\lambda}=-\rho_{n}\bar{\psi}(\hat{\theta}_1)/||\bar{\psi}(\hat{\theta}_1)||$, which yields $\tilde{\lambda}\in\Lambda_n$. Then, we have $\max_{1\leq i\leq n}|\tilde{\lambda}^{\T}\psi(Z_i;\hat{\theta}_1)|=o_p(1)$
and $\tilde{\lambda}\in\widehat{\Lambda}_n(\hat{\theta}_1)$ w.p.a.1. For any $\lambda_d$ lying in the joining line between $\tilde{\lambda}$ and 0, it follows from Lemma \ref{lem1} that w.p.a.l. $\max_{1\leq i\leq n}\tilde\rho_2(\lambda_d^{\T}\psi(Z_i;\hat{\theta}))\leq C^*$. Taking Taylor's expansion of $\bar{\ell}(\tilde{\lambda},\hat{\theta}_1)$ yields
\[
\begin{array}{llll}
\bar{\ell}(\tilde{\lambda}, \hat{\theta}_{1})&=&\tilde{\lambda}^{\T}\bar{\psi}(\hat{\theta}_{1})+\frac{1}{2}\tilde{\lambda}^{\T}
\left\{\frac{1}{n}\sum\limits_{i=1}^n\tilde\rho_2(\lambda_d^{\T}\psi(Z_i;\hat{\theta}_{1}))\psi(Z_i;\hat{\theta}_{1})\psi^{\T}(Z_i;\hat{\theta}_{1})\right\}\tilde{\lambda}\\
&\leq&-\rho_n||\bar{\psi}(\hat{\theta}_{1})||+C^*\rho_{n}^2.
\end{array}
\]
On the other hand, we have
$\bar{\ell}_p(\bar{\lambda}, \hat{\theta}_{1})\ge\inf_{\lambda\in\widehat{\Lambda}_n(\hat{\theta}_1)}\bar{\ell}_p(\lambda, \hat{\theta}_{1})\ge \inf_{\lambda\in\widehat{\Lambda}_n(\theta_{10})}\bar{\ell}_p(\lambda, \theta_{10})$.
By $\gamma=O_p(k/(nq))$ given in Assumption \ref{ass7} and Lemma \ref{lem2}, we obtain
$$ \inf\limits_{\lambda\in\widehat{\Lambda}_n(\theta_{10})}\bar{\ell}_p(\lambda, \theta_{10})=\inf\limits_{\lambda\in\widehat{\Lambda}_n(\theta_{10})}\bar{\ell}(\lambda, \theta_{10})-\sum\limits_{j=1}^qp_\gamma(|\theta_{10j}|)=O_p(k/n).$$
It follows from $\bar{\ell}_p(\lambda, \theta_{1})\leq \bar{\ell}(\lambda, \theta_{1})$ for any $\theta_1\in\Theta_1$ and $\lambda\in\widehat{\Lambda}_n(\theta_{1})$ that
\[
\begin{array}{llll}
-\rho_n||\bar{\psi}(\hat{\theta}_{1})||+C^*\rho_n^2\ge \bar{\ell}(\bar{\lambda}, \hat{\theta}_{1})\ge \inf\limits_{\lambda\in\widehat{\Lambda}_n(\theta_{10})}\bar{\ell}_p(\lambda, \theta_{10})= O_p(k/n),
\end{array}
\]
which indicates $||\bar{\psi}(\hat{\theta}_{1})||= O_p(\rho_n)$. Now we consider any $\varepsilon_n\rightarrow 0$.
Let $\ddot{\lambda}=-\varepsilon_n\bar{\psi}(\hat{\theta}_{1})$. It follows from $\ddot{\lambda}=o_p(\rho_n)$ that $\ddot{\lambda}\in\Lambda_n$ w.p.a.1. Using the same argument given above yields $O_p(k/n)\leq\ddot{\lambda}^{\T}\bar{\psi}(\hat{\theta}_{1})+C||\ddot{\lambda}||^2
=-\varepsilon_{n}||\bar{\psi}(\hat{\theta}_{1})||^2+C\varepsilon_{n}^2||\bar{\psi}(\hat{\theta}_{1})||^2$.
For enough large $n$, $1-\varepsilon_nC$ is bounded away from zero. Thus, it follows that $\varepsilon_{n}||\bar{\psi}(\hat{\theta}_{1})||^2=O_p(k/n)$, which leads to $||\bar{\psi}(\hat{\theta}_{1})||=O_p(\sqrt{k/n})$.

Using the same augment given in Lemma \ref{lem2}, it follows from $||\bar{\psi}(\hat{\theta}_1)||=O_p(\sqrt{k/n})$ that $||\lambda(\hat{\theta}_1)||=O_p(\sqrt{k/n})$.

\begin{lemma}\label{lem4}
If Assumptions \ref{ass1}, \ref{ass3}(i)and \ref{ass7} hold, we have $||\nu(\hat{\theta})||=O_p(\sqrt{r/n})$, where $\hat{\theta}$ is defined in Theorem \ref{th1}.
\end{lemma}
\noindent \emph{\bf Proof}.
By $\max_{1\leq i\leq n}\sup_{\theta\in\Theta_{s(n)}}|\nu^{\T}g(X_i;\theta)|=o_p(1)$ given in Lemma \ref{lem1} and the similar proof of Lemma \ref{lem2}, it follows from Assumptions \ref{ass1} and \ref{ass3}(i) that $\inf_{\nu\in\widehat{\mathcal{V}}_n(\theta_{0})}\ell(\nu, \theta_{0})= O_p(r/n)$.
Following the same argument as given in Lemma \ref{lem3} and by Assumption \ref{ass7}: $\gamma=O_p(k/(nq))=o_p(r/(nq))$,
we can obtain  $||\bar{g}(\hat{\theta})||=O_p(\sqrt{r/n})$. Again, following the same arguments as given in Lemma \ref{lem2} and Lemma \ref{lem3},
we have  $||\nu(\hat{\theta})||=O_p(\sqrt{r/n})$.

\vspace{5mm}

\noindent  \emph{\bf Proof of Theorem \ref{th1}}.  Following Fan and Lv (2011), we divide the proof procedure of Theorem \ref{th1} into two steps.
First, we prove the consistency of the proposed PET estimator in the $q$-dimensional subspace. To this end,
we consider restricting $\ell_p(\theta)$ into the $q$-dimensional subspace $\{\theta\in\mathcal{R}^p:\theta_{\mathbb{J}^c}=0\}$ of $\mathcal{R}^p$. The corresponding constrained PET likelihood function is $\bar{\ell}_p(\theta_1)$ given in Equation (\ref{cPET}). For the constrained subspace, it follows from $p_\gamma(0)=0$ and $L(\theta)=L(\theta_1)$ that $\ell_p(\theta)=\bar{\ell}_p(\theta_1)$, where $L(\theta)$ is defined in Equation (2.1). Hence, $(\hat\theta_1,0)^{\T}$ is the maximizer of $\ell_p(\theta)$ on the constrained subspace, where $\hat{\theta}_{\mathbb{J}}=\hat{\theta}_1=\arg\max_{\theta_1\in\Theta_1}\bar{\ell}_p(\theta_1)$.
From Lemma \ref{lem3}, we have $||\bar{\psi}(\hat{\theta}_1)||=O_p(\sqrt{k/n})$. Following the argument of Chang, Chen and Chen (2013), if $||\hat{\theta}_1-\theta_{10}||$ does not converge to zero in probability, there exists a subsequence $\{n^{*},k^{*},q^{*}\}$ such that $||\hat{\theta}_{1n^{*}}-\theta_{10}||\ge\varepsilon$ a.s. for some positive constant $\varepsilon$. By Assumption \ref{ass6}, we have $||E\{\psi(Z_i;\hat{\theta}_{1n^{*}})\}||=o_p\{\zeta_1(k^{*},q^{*})\}+O_p(\sqrt{k^{*}/n^{*}})$, which is in conflict with
$||E\{\psi(Z_i;\hat{\theta}_{1n^{*}})\}||\ge\zeta_1(k^{*},q^{*})\zeta_2(\varepsilon)$ because of $\lim\inf_{k,q\rightarrow\infty}\zeta_1(k,q)>0$. Therefore, we obtain $||\hat{\theta}_1-\theta_{10}||\rightarrow 0$ as $n\rightarrow \infty$. Assumption \ref{ass5}(i) implies $||\bar{\psi}(\hat{\theta}_1)-\bar{\psi}(\theta_{10})||\ge C||\hat{\theta}_1-\theta_{10}||$ w.p.a.1. Then, we have $||\hat{\theta}_1-\theta_{10}||=O_p(\sqrt{k/n})$.

On the other hand, the sparsity property of the proposed PET estimator can be concluded from Theorem \ref{th7}, Hence, we have proved Theorem \ref{th1}.

\vspace{2mm}

\noindent  \emph{\bf Proof of Theorem \ref{th2}}. By Theorem \ref{th1}, we only need to prove asymptotic
normality of $\hat{\theta}_1$. Let $S(\lambda,\theta_1)=\log n^{-1}\sum_{i=1}^{n}\exp\{\lambda^{\T}\psi(Z_i;\theta_1)\}-\sum_{j=1}^{q}p_{\gamma}(|\theta_{1j}|)$, where $\theta_{1j}$ is the $j$th component of $\theta_1$.
Then, the constrained PET Likelihood $\bar{\ell}_{p}(\theta_1)$ given in Equation (\ref{cPET}) can be written as $\bar{\ell}_{p}(\theta_1)=S(\lambda,\theta_1)$. Let
$S_{1}(\lambda,\theta_1)=\partial S(\lambda,\theta_1)/\partial \lambda=\sum_{i=1}^{n}\pi_i\psi(Z_{i};\theta_1)$,
$S_{2}(\lambda,\theta_1)=\partial S(\lambda,\theta_1)/\partial \theta_1=\sum_{i=1}^{n}\pi_i\{\partial_{\theta_1}
\psi(Z_{i};\theta_1)\}^{\T}\lambda-W(\theta_1)$,
where $\partial_{\theta_1}\psi(Z_i;\theta_1)=\partial\psi(Z_i;\theta_1)/\partial\theta_1^{\T}$, $\pi_i=\exp\{\lambda^{\T}\psi(Z_{i};\theta_1)\}/$ $\sum_{j=1}^{n}\exp\{\lambda^{\T}\psi(Z_j;\theta_1)\}$,
and the $j$th component of vector $W(\theta_1)$ is $p_{\gamma}^{'}(|\theta_{1j}|){\rm sign}(\theta_{1j})$ for $j=1,\ldots,q$.
Thus, it follows from the definitions of $\hat{\lambda}$ and $\hat{\theta}_1$ that $\hat{\lambda}$ and $\hat{\theta}_1$ satisfy $S_{k}(\hat{\lambda},\hat{\theta}_1)=0$ for $k=1,2$.

Let $\bar{\Sigma}_1(\theta_{10})=\frac{1}{n}\sum_{i=1}^{n}\psi(Z_{i};\theta_{10})\psi^{\T}(Z_{i};\theta_{10})$ and $\bar{\Gamma}_1(\theta_{10})
=\frac{1}{n}\sum_{i=1}^{n}\partial_{\theta_1}^{\T}\psi(Z_{i};\theta_{10})$. Then, we have
\[
\begin{array}{llll}
S_{11}(0,\theta_{10})=\partial S(0,\theta_{10})/\partial \lambda\lambda^{\T}=\bar{\Sigma}_1(\theta_{10})-n^{-2}\left\{\sum\limits_{i=1}^{n}\psi(Z_i;\theta_{10})\right\}\left\{\sum\limits_{i=1}^{n}\psi(Z_i;\theta_{10})\right\}^{\T},\\
S_{12}(0,\theta_{10})=\partial S(0,\theta_{10})/\partial\lambda\theta_1^{\T}=\bar{\Gamma}_1^{\T}(\theta_{10}),~~~
S_{21}(0,\theta_{10})=\partial S(0,\theta_{10})/\partial\theta_1\lambda^{\T}=\bar{\Gamma}_1(\theta_{10}),\\
S_{22}(0,\theta_{10})=\partial S(0,\theta_{10})/\partial\theta_1\theta_1^{\T}=0.
\end{array}
\]
Let $\Sigma_1=E\{\bar{\Sigma}_1(\theta_{10})\}$ and $\Gamma_1=E\{\bar{\Gamma}_1(\theta_{10})\}$. Taking Taylor's expansion of $S_{k}(\hat{\lambda},\hat{\theta}_1)=0$ ($k=1,2$) at $(0,\theta_{10})$ yields
\begin{equation}\label{PETLA1}
\left(
\begin{array}{*{8}c}
-S_1(0,\theta_{10})\\
0
\end{array}
\right)
=
\left(
\begin{array}{*{8}c}
\Sigma_1&\Gamma_1^{\T}\\
\Gamma_1&0
\end{array}
\right)
\left(
\begin{array}{*{8}c}
\hat{\lambda}-0\\
 \hat{\theta}_1-\theta_{10}
\end{array}
\right)
+R_n,
\end{equation}
where $R_{n}=\sum_{j=1}^{5}R_{jn}$, $R_{1n}=(R_{1n}^{\T(1)}$, $R_{1n}^{\T(2)})^{\T}$ in which $R_{1n}^{(1)}\in\mathcal{R}^{k}$ and $R_{1n}^{(2)}\in\mathcal{R}^{q} $ and the $j$th component of $R_{1n}^{(l)}$ is $R_{1n,j}^{(l)}=\frac{1}{2}(\hat{\Delta}-\Delta_{0})^{\T}\partial_\Delta^{2}S_{l,j}(\Delta^*)(\hat{\Delta}-\Delta_{0})$ for $l=1,2$ and $j=1,\ldots,k+q$,  $\Delta=(\lambda^{\T}, \theta_1^{\T})^{\T}$, $\partial_\Delta^2S_l=\partial^2S_l/\partial\Delta\partial\Delta^{\T}$
and $\Delta^*=({\lambda^*}^{\T}, {\theta_1^*}^{\T})^{\T} $ satisfying $ ||\lambda^*||\leq||\hat{\lambda} ||$ and
$||\theta_1^*-\theta_{10}||\leq||\hat{\theta}_1-\theta_{10}||$. Other terms $R_{2n},\ldots,R_{5n}$ are shown as follows.

Define
$$Q_n=
\left(
\begin{array}{*{8}c}
S_{11}(0,\theta_{10})&S_{12}(0,\theta_{10})\\
S_{21}(0,\theta_{10})&S_{22}(0,\theta_{10})
\end{array}
\right),~~~
Q=
\left(
\begin{array}{*{8}c}
\Sigma_1&\Gamma_1^{\T}\\
\Gamma_1&0
\end{array}
\right).$$
Following the argument of Fan and Peng (2004),
it is easily shown that
$$P(||Q_n-Q||>\epsilon)\leq \frac{1}{\epsilon^2}\sum_{i,j=1}^{k+q}E\left\{ \frac{\partial^2S(0,\theta_{10})}{\partial\Delta_i\partial\Delta_j}-E\frac{\partial^2S(0,\theta_{10})}{\partial\Delta_i\partial\Delta_j}\right\}^2= O_p\left(\frac{(k+q)^2}{n}\right).$$
Let $\tilde{G}_n=(0_{d\times k},G_n\mathcal{K}^{-1/2})$, where $G_n$ is defined in Theorem \ref{th2}, and $\mathbb{S}^*=(-S_1^{\T}(0,\theta_{10}),0^{\T})^{\T}$.  From $Q_n^{-1}-Q^{-1}=-Q_n^{-1}(Q_n-Q)Q^{-1}$, we have
\[
\begin{array}{llll}
||\sqrt{n}\tilde{G}_n(Q_n^{-1}-Q^{-1})\mathbb{S}^*||^2&\leq&
n\mathbb{E}_{\rm max}(G_nG_n^T)\mathbb{E}_{\rm max}(\mathcal{K}^{-1})\mathbb{E}_{\rm max}^{-2}(Q_n^{-1})||(Q_n-Q)Q^{-1}\mathbb{S}^*||^2\\
&\leq&n\mathbb{E}_{\rm max}(G_nG_n^T)\mathbb{E}_{\rm max}(\mathcal{K}^{-1})\mathbb{E}_{\rm max}^{-2}(Q_n^{-1})||(Q_n-Q)||^2||Q^{-1}\mathbb{S}^*||^2\\
&=&O_p(k(k+q)^4/n^2)=o_p(1).
\end{array}
\]
From the implicit theorem and envelope theorem,  we have $\hat{\lambda}=\lambda(\hat{\theta}_1)=O_p(\sqrt{k/n})$. It follows from the Cauchy-Schwarz inequality and $k^2(k+q)^3=o_p(n)$ that
$$||R_{1n}^{(1)}||^2\le n^{-2}||\hat{\Delta}-\Delta_0||^4n^2\sum_{i,j,l=1}^{k+q}\partial^2_{\Delta}S_{1,l}(\Delta^*)/\partial\Delta_i\partial\Delta_j=O_p(k^2(k+q)^3/n^2)=o_p(1/n).$$
By the definitions of $R_{1n}^{(1)}$ and $R_{1n}^{(2)}$, we obtain $||R_{1n}^{(2)}||^2=||R_{1n}^{(1)}||^2$, which leads to $||R_{1n}||=o_p(\sqrt{1/n})$.

When $k\sqrt{k+q}=o_p(\sqrt{n})$, which is a relaxed condition of $k^2(k+q)^3=o_p(n)$, it follows from
$\sup_{\Delta\in\mathcal{P}}\mathbb{E}_{\max}\{\partial_\Delta^{2}S_{l,j}(\Delta)\}=O_p(1)$ given in the remark  of  Theorem \ref{th2}(ii) that $||R_{1n}^{(l)}||_\infty\leq\sup_{\Delta\in\mathcal{P}}\max_{j\in \{1,\ldots,k+q\}}\mathbb{E}_{\max}\{\partial_\Delta^{2}S_{l,j}(\Delta)\}
||\hat{\Delta}-\Delta_0||^2=O_p(k\sqrt{k+q}/n)=o_p(1/\sqrt{n})$ for $l=1,2$.

By Assumption \ref{ass8}, we obtain $||R_{2n}||=||(0,W^{\T}(\theta_{10}))^{\T}||=o_p(\sqrt{1/n})$ and $||R_{3n}||=||(0,W^{'}(\theta_1^*)(\hat{\theta}_1-\theta_{10}))^{\T}||=o_p(\sqrt{1/n})$.
Following Assumption 4, Lemma \ref{lem3}, Theorem \ref{th1} and $k^2(k+q)^3=o_p(n)$, we have
$||R_{4n}||=||(\{(\bar{\Sigma}_1(\theta_{10})-\Sigma_1)\hat{\lambda}\}^{\T}+
\{(\bar{\Gamma}_1(\theta_{10})-\Gamma_1)(\hat{\theta}_1-\theta_{10})\}^{\T},0^{\T})^{\T}||=o_p(\sqrt{1/n})$ and
$||R_{5n}||=||(0^{\T},\hat{\lambda}^{\T}(\bar{\Gamma}_1(\theta_{10})-\Gamma_1))^{\T}||=o_p(\sqrt{1/n})$. Combining the above equations yields $||R_{n}||=o_p(\sqrt{1/n})$.

It follows from Equation (\ref{PETLA1}) that
\begin{equation} \label{PETLA2}
\left(
\begin{array}{*{8}c}
\hat{\lambda}-0\\
\hat{\theta}_1-\theta_{10}
\end{array}
\right)
=Q^{-1}
\left\{\left(
\begin{array}{*{8}c}
-S_1(0,\theta_{10})\\
0\\
\end{array}
\right)+R_n
\right\},
\end{equation}
which leads to
\begin{equation}\label{PETLA4}
\hat\theta_1-\theta_{10}=-\mathcal{K}\Gamma_{1}\Sigma_{1}^{-1}\{S_1(0,\theta_{10})+R_{2n}\},
\end{equation}
where $||R_{2n}||=o_p(\sqrt{1/n})$, $\mathcal{K}=(\Gamma_{1}\Sigma_{1}^{-1}\Gamma_{1}^{\T})^{-1}$ and $S_1(0,\theta_{10})=\bar{\psi}=\frac{1}{n}\sum_{i=1}^{n}\psi(Z_{i};\theta_{10})$.

Let $B=-\mathcal{K}\Gamma_{1}\Sigma_{1}^{-1}$,
$X_{ni}=n^{-1/2}G_n\mathcal{K}^{-1/2}B\psi(Z_i;\theta_{10})=n^{-1/2}Y_{ni}$. Then, for any $\epsilon>0$, we have
\[
\begin{array}{llll}
\sum_{i=1}^{n}E||X_{ni}||^2I(||X_{ni}||>\epsilon)
&=&nE||X_{n1}||^2I(||X_{n1}||>\epsilon)\\
&\leq&n\left\{E||X_{n1}||^4\right\}^{1/2}\left\{P(||X_{n1}||>\epsilon\right\}^{1/2}.
\end{array}
\]
It follows from $G_nG_n^{\T}\rightarrow V$ that $P(||X_{n1}||>\epsilon)\le E||Y_{n1}||^2/(n\epsilon^2)=O_p(n^{-1})$ and
\[
\begin{array}{llll}
E||X_{n1}||^4&=&n^{-2}E\left\{\psi^{\T}(Z_i;\theta_1)B\mathcal{K}^{-1/2}G_n^{\T}G_n\mathcal{K}^{-1/2}B\psi(Z_i;\theta_1)\right\}^2\\
&\leq&n^{-2}\mathbb{E}^2_{\rm max}(G_nG_n^{\T})\mathbb{E}^2_{\rm max}(\mathcal{K}^{-1})E(\psi^{\T}\psi)^2\\
&=& O_p(\frac{q^2}{n^2}).
\end{array}
\]
Combining the above equations yields
$\sum_{i=1}^{n}E||X_{ni}||^2I(||X_{ni}||>\epsilon)=O_p(\frac{1}{\sqrt{n}})=o_p(1)$.

Since $\sum_{i=1}^{n}{\rm cov}(X_{ni})=n{\rm cov}(X_{n1})=
G_nG_n^{\T}\rightarrow V$ as $n\rightarrow \infty$ and $BS_1(0,\theta_{10})=(\hat{\theta}_1-\theta_{10})$, it follows from the central limit theorem that
$\sqrt{n}G_n\mathcal{K}^{-1/2}(\hat{\theta}_1-\theta_{10})\stackrel{{\cal L}}{\rightarrow}\mathcal {N}(0,V)$.
Thus, we have proved Theorem \ref{th2}.

\vspace{5mm}

\noindent  \emph{\bf Proof of Theorem \ref{th3}}. The method given in Theorem \ref{th1}0 of Schennach (2007) can be used to prove Theorem \ref{th3}(i) based on
$\ell_p(\theta)=\log\frac{1}{n}\sum_{i=1}^{n}\exp\{\nu^{\T}(\theta)g(X_i;\theta)\}-\sum_{j=1}^{p}p_{\gamma}(|\theta_j|)$ and the assumption on the continuity of the penalty function $p_\gamma(\cdot)$.

Now we prove the sparsity of the proposed PET estimator. The first-order partial derivative of the PET likelihood $\ell_p(\theta)$ with respect to $\theta_j$ for $j\notin\mathbb{J}$ is given by
\[
\begin{array}{llll}
\partial\ell_{p}(\theta)/\partial\theta_j
&=&\frac{\frac{1}{n}\sum\limits_{i=1}^{n}\exp\{\nu^{\T}g(X_{i};\theta)\}\partial_{\theta_j}g^{\T}(X_i;\theta)}
{\frac{1}{n}\sum\limits_{i=1}^{n}\exp\{\nu^{\T}g(X_i;\theta)\}}\nu-p_{\gamma}^{'}(|\theta_{j}|){\rm sign}(\theta_{j})\\
&=&\frac{\{h_{\nu j}+O_p(\sqrt{1/n})\}\nu}{E\exp\{\nu^{\T}g(X_i;\theta)\}+O_p(\sqrt{1/n})}-p_{\gamma}^{'}(|\theta_{j}|){\rm sign}(\theta_{j})\\[1mm]
&\stackrel{\Delta}{=}&\mathcal{J}_{1}+\mathcal{J}_{2},
\end{array}
\]
where  $\partial_{\theta_j}g^{\T}(X_i;\theta)=\partial g^{\T}(X_i;\theta)/\partial\theta_j$, and
$h_{\nu j}=E[\exp\{\nu^{\T}g(X_{i};\theta)\}\partial_{\theta_j}g^{\T}(X_i;\theta)]$.
By Assumption \ref{ass14}, combining the above equations yields $\mathcal{J}_{1}\leq O_p(1)||\nu||$.
Thus, for $j\notin\mathbb{J}$, we have
\[
\partial\ell_{p}(\theta)/\partial\theta_j
=O_p(1)||\nu||-p_{\gamma}^{'}(|\theta_{j}|){\rm sign}(\theta_{j})\stackrel{\Delta}{=}\gamma\{-\frac{p_{\gamma}^{'}(|\theta_{j}|)}{\gamma}{\rm sign}(\theta_{j})\}+\frac{||\nu||}{\gamma}O_p(1)\},
\]
which leads to $\partial\ell_p(\theta)/\partial\theta_j=\gamma\{-p_{\gamma}^{'}(|\theta_j|){\rm sign}(\theta_j)/\gamma+o_p(1)\}$ implying that the sign of $\partial\ell_{p}(\theta)/\partial\theta_j$ is dominated by the sign of $\theta_j$. Then, for any $j\notin\mathbb{J}$ and as $n\rightarrow\infty$,
we have
$\partial\ell_{p}(\theta)/\partial\theta_j<0$ when $\theta_j>0$, and $\partial\ell_{p}(\theta)/\partial\theta_j>0$ when $\theta_{j}<0$ with probability tending to one,
which means $\hat{\theta}_{2}=0$ with probability tending to one. Thus, Theorem \ref{th3}(ii) holds.

To prove Theorem \ref{th3}(iii), we require the following Lemma. For convenience, we define
$I_p =(H^{\T}_1, H^{\T}_2)^{\T}$, where $H_1\in\mathcal{R}
^{q\times p}$
and $H_2\in\mathcal{R}^{(p-q)\times p}$.

\begin{lemma}\label{lem5} Let $\hat{\phi}=(\hat{\rho},\hat{\tau},\hat{\nu},\hat{\theta})$ be the PET estimator of $\phi=(\rho,\tau,\nu,\theta)$,
where $\rho=n^{-1}\sum_{i=1}^n\rho_i$ with $\rho_i=\exp(\nu^{\T}g(X_i;\theta))$. Then, $\hat\phi$ is the solution to $n^{-1}\sum_{i=1}^n\Psi(X_i;\phi)=0$, where $\Psi^{\T}(X_i;\phi)=(\rho_i-\rho,\rho_iH_2\theta,\rho_ig(X_i;\theta),\rho_i\Gamma_i^{*\T}\nu-\rho_iW(\theta)+\rho_iH_2^{\T}\tau)$
with $\Gamma_i^*=\partial_{\theta}g(X_i;\theta)$.

\end{lemma}

\noindent  \emph{\bf Proof}. Define
$S(\nu,\theta,\tau)=\log n^{-1}\sum_{i=1}^{n}\exp\{\nu^{\T}g(X_i;\theta)\}-\sum_{j=1}^{p}p_{\gamma}(|\theta_{j}|)+\tau^{\T}H_{2}\theta$.
Following Tang and Leng (2012) and using the sparsity of $\hat{\theta}$, we obtain that
$\hat{\theta}$ satisfies $\sum_{i=1}^{n}\hat{\pi}_i\hat{\Gamma}_i^{*\T}\hat{\nu}-W(\hat{\theta})+H_2^{\T}\hat{\tau}=0$, which leads to $\sum_{i=1}^n(\hat\rho_i\hat\Gamma_i^{*\T}\hat\nu-\hat\rho_iW(\hat\theta)+\hat\rho_iH_2^{\T}\hat\tau)=0$,
where $\pi_i=\rho_i/\sum_{j=1}^{n}\rho_j$; and $\hat\nu$ and $\hat\tau$ satisfy
$\sum_{i=1}^{n}\hat\pi_ig(X_i;\hat\theta)=0$ and $H_2\hat\theta=0$, which lead to $\sum_{i=1}^n\hat\rho_ig(X_i;\hat\theta)=0$ and $\sum_{i=1}^n\hat\rho_iH_2\hat\theta=0$, respectively. Also, it follows from the definition of $\rho$ that $\sum_{i=1}^n(\rho_i-\rho)=0$. The above equations show that $\hat\phi$ is the solution to
$n^{-1}\sum_{i=1}^n\Psi(X_i;\phi)=0$.

\vspace{5mm}

\noindent  \emph{\bf Proof of Theorem \ref{th3}(iii)}. Lemma \ref{lem5} presents a just-identified GMM estimator, thus we can apply Theorem \ref{th3}.4 of Newey and McFadden (1994) to the just-identified case if we can show that (i) $E\{\sup_{\phi\in\mathcal{N}_\Phi}||\partial\Psi^{\T}(X_i;\phi)/\partial\phi||\}<\infty$ holds for some neighborhood $\mathcal{N}_\Phi$ of $\phi$,
and (ii) $E\{\Psi(X_i;\phi)\Psi^{\T}(X_i;\phi)\}$ exists.

Consider components of matrix $\partial\Psi^{\T}(X_i;\phi)/\partial\phi$, which is given by $\omega\exp\{k_{1}\nu^{\T}g(X_i;\theta)\}$
$g^{k_{g}}\Gamma^{k_{\Gamma}}\Omega^{k_{\Omega}}$ for $0\leq k_{g}+
k_{\Gamma}+k_{\Omega}\leq 2$ and $k_{1}=0,1$, where $g, \Gamma$ and $\Omega$ denote elements of $g(X_i;\theta),\Gamma(X_i;\theta)$ and $\Omega_{jl}(X_i;\theta)$, respectively, and $\omega$ denotes the product of elements of $\phi$ that is necessarily bounded for $\phi\in\mathcal{N}_\Phi$ and also includes the first and second partial derivatives of penalty function such as $p_\gamma^{'}(\theta)$ and $p_\gamma^{''}(\theta)$. By Assumption \ref{ass14}, we can obtain (i). It follows from $\exp\{k_{1}\nu^{\T}g(X_i;\theta)\}|g|^{k_{g}}|\Gamma|^{k_{\Gamma}}|\Omega|^{k_{\Omega}}\leq\exp\{k_{1}\nu^{\T}g(X_i;\theta)\}
|f(X_i)|^{k_{g}+k_{\Gamma}+k_{\Omega}}$ that
$$
\begin{array}{lll}
E[\sup_{\phi\in\mathcal{N}_\Phi}\exp\{k_{1}\nu^{\T}g(X_i;\theta)\}|g|^{k_{g}}|\Gamma|^{k_{\Gamma}}|\Omega|^{k_{\Omega}}]\leq E[\sup_{\phi\in\mathcal{N}_\Phi}\exp\{k_{1}\nu^{\T}g(X_i;\theta)\}|f(X_i)|^{k_2}]\\
=E[\sup_{\theta\in\cal N}\sup_{\nu\in\Lambda(\theta)}\exp\{k_{1}\nu^{\T}g(X_i;\theta)\}(f(X_i))^{k_2}]<\infty.
\end{array}
$$ Matrix $\Psi(X_i;\phi)\Psi^{\T}(X_i;\phi)$ has elements of the form $\omega\exp\{k_{1}\nu^{\T}g(X_i;\theta)\}g^{k_{g}}\Gamma^{k_{\Gamma}}$ with $k_{1}=0,1,2$, and $0\leq
k_{g}+k_{\Gamma}\leq 2$. Similar argument implies (ii).

\vspace{5mm}

\noindent  \emph{\bf Proof of Theorem \ref{th4}}.  It follows from the argument of Theorem \ref{th2} that $\hat{\lambda}$ and $\hat{\theta}_1$ can be obtained by maximizing $\bar{\ell}_p(\theta_1)$  under $H_0\cup H_1$, which indicates that
\begin{equation}\label{PETL7}
\hat\theta_1-\theta_{10}=-\mathcal{K}\Gamma_{1}^{\T}\Sigma_{1}^{-1}\bar{\psi}+R_{2n},~~~
\hat{\lambda}=\{\Sigma_1^{-1}\Gamma_1\mathcal{K}\Gamma_1^{\T}\Sigma_1^{-1}-\Sigma_1^{-1}\}\bar{\psi}+R_{1n}.
\end{equation}
Let $\tilde{w}_i=\hat{\lambda}^{\T}\psi(z_i;\hat{\theta}_1)$. It follows from Lemma \ref{lem1} that $\max_{1\leq i\leq n}|\tilde{w}_i|=o_p(1)$. Taking Taylor expansion of $\bar{\ell}(\hat{\theta}_1)$ at $\tilde{w}_i$ leads to
$$\bar{\ell}(\hat{\theta}_1)=\log\left\{1+\frac{1}{n}\sum\limits_{i=1}^{n}\tilde{w}_i(1+o_p(1))\right\}=\frac{1}{n}\sum\limits_{i=1}^{n}\tilde{w}_i(1+o_p(1)).$$
Substituting expressions of $\hat{\lambda}$ into $\bar{\ell}_p(\hat{\theta}_1)$ yields
$2n\bar{\ell}_p(\hat{\theta}_1,\hat{\lambda})=-n\bar{\psi}^{\T}\{\Sigma_1^{-1}\Gamma_1\mathcal{K}\Gamma_1^{\T}\Sigma_1^{-1}-\Sigma_1^{-1}\}\bar{\psi}+o_p(1)$, where $o_p(1)$ includes penalty function. It follows from $B_{n}\theta_{1}=0$ and $B_nB_n^{\T}=I_d$
and Theorem 3 of Tang and Leng (2012) that the constrained PET estimators $\tilde{\lambda}$ and $\tilde{\theta}_1$ of $\lambda$ and $\theta_1$ under $H_0$ can be obtained by maximizing
$\tilde{\ell}_p(\lambda,\theta_1,\tau) =\log n^{-1}\sum_{i=1}^{n}\exp\{\lambda^{\T}\psi(Z_i$; $\theta_1)\}-\sum_{j=1}^{q}p_{\gamma}(|\theta_{1j}|)+\tau^{\T}B_{n}\theta_1$.
It is easily shown that $\tilde{\lambda}=(\Sigma_1^{-1}\Gamma_1\mathbb{P}\Gamma_1^{\T}\Sigma_1^{-1}-\Sigma_1^{-1})\bar{\psi}+\tilde{R}_{1n}$, where $\mathbb{P}=\mathcal{K}B_{n}^{\T}(B_{n}\mathcal{K}B_{n}^{\T})^{\T}B_{n}\mathcal{K}-\mathcal{K}$. Thus, the constrained maximum PET likelihood is given by
$2n\bar{\ell}_p(\tilde{\theta}_1, \tilde{\lambda})=-n\bar{\psi}^{\T}(\Sigma_1^{-1}\Gamma_1\mathbb{P}\Gamma_1^{\T}\Sigma_1^{-1}-\Sigma_1^{-1})\bar{\psi}+o_p(1)$.
Then, the constrained PET likelihood ratio statistic for testing $H_0:B_n\theta_{10}=0$ is given by
$$\hat{\ell}(B_n)=2n\bar{\ell}_p(\hat{\theta}_1, \hat{\lambda})-2n\bar{\ell}_p(\tilde{\theta}_1, \tilde{\lambda})=n\bar{\psi}^{\T}\Sigma_1^{-1/2}(\mathbb{T}_1-\mathbb{T}_2)\Sigma_1^{-1/2}\bar{\psi}+o_p(1),$$
where $\mathbb{T}_1=\Sigma_1^{-1/2}\Gamma_1\mathbb{P}\Gamma_1^{\T}\Sigma_1^{-1/2}$ and $\mathbb{T}_2=\Sigma_1^{-1/2}\Gamma_1\mathcal{K}\Gamma_1^{\T}\Sigma_1^{-1/2}$.
It is easily shown that $\mathbb{T}_1$ and $\mathbb{T}_2$ are symmetric idempotent matrices, and the rank of matrix $\mathbb{T}_1-\mathbb{T}_2$ is $d$,
which indicates that there is a matrix $\mathcal{T}$ such that $\mathbb{T}_1-\mathbb{T}_2=\mathcal{T}^{\T}\mathcal{T}$ and $\mathcal{T}\mathcal{T}^{\T}=I_d$
(Fan and Peng, 2004). Also, it follows from the center limit theorem that
$\sqrt{n}\mathcal{T}\Sigma^{-1/2}\bar{\psi}\stackrel{{\cal L}}{\rightarrow}\mathcal {N}(0,I_d)$, which leads to
$n\bar{\psi}^{\T}\Sigma^{-1/2}(\mathbb{T}_1-\mathbb{T}_2)\Sigma^{-1/2}\bar{\psi} \stackrel{{\cal L}}{\rightarrow}\mathcal {\chi}_d^2$.

\vspace{5mm}

\noindent  \emph{\bf Proof of Theorem \ref{th5}}. For simplicity, we denote $M_j^*=E\{\partial^2m(Z_i;\eta_0)/\partial\eta_j\partial\eta^{\T}\}$ in which $\eta_j$ is the $j$th element of $\eta$, $M_{jt}^*=E\{\partial^3m(Z_i;\eta_0)/\partial\eta_j\partial\eta_t\partial\eta^{\T}\}$,
$\tilde{A}=\frac{1}{\sqrt{n}}\sum_{i=1}^n\partial m(Z_i;\eta_0)/\partial\eta^{\T}$ $-\sqrt{n}\mathbb{M}$,
$\tilde{B}_j=\frac{1}{\sqrt{n}}\sum_{i=1}^n\partial^2m(Z_i;\eta_0)/\partial\eta_j\partial\eta^{\T}-\sqrt{n}M_j^*$,
$\tilde{\upsilon}=-\frac{1}{\sqrt{n}}\sum_{i=1}^n\mathbb{M}^{-1}m(Z_i;\eta_0)$,
$\tilde{a}=\mathbb{M}^{-1}\sum_{j=1}^\mathcal{S}\tilde{\upsilon}_jM_j^*\tilde{\upsilon}$ in which $\tilde{\upsilon}_j$ is the $j$th component of $\tilde{\upsilon}$, $\tilde{b}=\mathbb{M}^{-1}\tilde{A}\tilde{\upsilon}$,
$\widehat{M}(\eta)=\frac{1}{n}\sum_{i=1}^n\partial m(Z_i;\eta)/\partial\eta^{\T}$.

Taking the Taylor expansion of $n^{-1}\sum_{i=1}^nm(Z_i;\hat\eta)$ at $\eta_0$ yields
\begin{equation}\label{high1}
\begin{array}{llll}
0=&\widehat{m}(\eta_0)+\widehat{M}(\eta_0)(\hat{\eta}-\eta_0)+\frac{1}{2}\sum\limits_{j=1}^{\mathcal{S}}(\hat{\eta}_j-\eta_{0j})\{\partial \widehat{M}(\eta_0)/\partial\eta_j\}(\hat{\eta}-\eta_0)\\
&+\frac{1}{6}\sum\limits_{j,t=1}^{\mathcal{S}}(\hat{\eta}_j-\eta_{0j})(\hat{\eta}_t-\eta_{0t})\{\partial^2 \widehat{M}(\bar{\eta})/\partial\eta_j\eta_t\}(\hat{\eta}-\eta_0),
\end{array}
\end{equation}
where $\widehat{m}(\eta_0)=n^{-1}\sum_{i=1}^nm(Z_i;\eta_0)$ and $\bar{\eta}$ lies in the jointing line between $\hat\eta$ and $\eta_0$.
Let $\widehat{M}=\widehat{M}(\eta_0)$.  Then, it is easily shown from Equation (\ref{high1}) that
\begin{equation}\label{high2}
\begin{array}{llll}
\hat{\eta}-\eta_0&=&\tilde{\upsilon}/\sqrt{n}-\mathbb{M}^{-1}\Big\{\tilde{A}(\hat{\eta}-\eta_0)/\sqrt{n}+\frac{1}{2}\sum\limits_{j=1}^{\mathcal{S}}(\hat{\eta}_j-\eta_{0j})
M_j^*(\hat{\eta}-\eta_0)\\
&&+\frac{1}{2}\sum\limits_{j=1}^{\mathcal{S}}(\hat{\eta}_j-\eta_{0j})\frac{\tilde{B}_j}{\sqrt{n}}(\hat{\eta}-\eta_0)\\
&&+\frac{1}{6}\sum\limits_{j,t=1}^{\mathcal{S}}(\hat{\eta}_j-\eta_{0j})(\hat{\eta}_t-\eta_{0t})M_{jt}^*(\hat{\eta}-\eta_0)\\
&&+\frac{1}{6}\sum\limits_{j,t=1}^{\mathcal{S}}(\hat{\eta}_j-\eta_{0j})(\hat{\eta}_t-\eta_{0t})(\partial^2 \widehat{M}(\bar{\eta})/\partial\eta_j\eta_t-M_{jt}^*)(\hat{\eta}-\eta_0)\Big\}.\\
\end{array}
\end{equation}
Denote $\hat{\eta}-\eta_0=\tilde{\upsilon}/\sqrt{n}-\mathbb{M}^{-1}\Big\{\mathcal{S}_1^*+\mathcal{S}_2^*+\mathcal{S}_3^*+\mathcal{S}_4^*+\mathcal{S}_5^*\Big\}$.
By the definitions of $\mathcal{S}_1^*,\ldots,\mathcal{S}_6^*$, we have $||\mathcal{S}_1^*||\leq ||\tilde{A}/\sqrt{n}||\cdot ||\hat{\eta}-\eta_0||=O_p(\sqrt{q^2/n})O_p(\sqrt{q/n})=O_p(\frac{q^{3/2}}{n})$,
$||2\mathcal{S}_2^*||\leq ||\hat{\eta}-\eta_0||\surd\{\sum_{l_1,l_2=1}^{\mathcal{S}}(\sum_{j=1}^{\mathcal{S}}(\hat{\eta}_j-\eta_{0j})M_{jl_1l_2}^*)^2\}
\leq ||\hat{\eta}-\eta_0||\surd\{\sum_{l_1,l_2=1}^{\mathcal{S}}||\hat{\eta}-\eta_0||^2(\sum_{j=1}^{\mathcal{S}}M_{jl_1l_2}^{*2})\}=O_p(\frac{q^{5/2}}{n})$,
$||2\mathcal{S}_3^*||\leq ||\hat{\eta}-\eta_0||\surd\{\sum_{l_1,l_2=1}^{\mathcal{S}}||\hat{\eta}-\eta_0||^2(\sum_{j=1}^{\mathcal{S}}(\widehat{M}_{jl_1l_2}^*-M_{jl_1l_2}^*)^2)\}
=O_p(\frac{q^{5/2}}{n\sqrt{n}})$,
$||6\mathcal{S}_4^*||\leq ||\hat{\eta}-\eta_0||\surd\{\sum_{l_1,l_2=1}^{\mathcal{S}}(\sum_{j,t=1}^{\mathcal{S}}(\hat{\eta}_j-\eta_{0j})(\hat{\eta}_t-\eta_{0t})M_{jtl_1l_2}^*)^2\}
\leq ||\hat{\eta}-\eta_0||\surd\{\sum_{l_1,l_2=1}^{\mathcal{S}}||\hat{\eta}-\eta_0||^4(\sum_{j=1}^{\mathcal{S}}(\sum_{t=1}^{\mathcal{S}}M^{*2}_{jtl_1l_2})^2)\}
=O_p(\frac{q^{4}}{n\sqrt{n}})$, and
$||6\mathcal{S}_5||\leq||\hat{\eta}-\eta_0||\surd\{\sum_{l_1,l_2=1}^{\mathcal{S}}||\hat{\eta}-\eta_0||^4(\sum_{j=1}^{\mathcal{S}}$ $(\sum_{t=1}^{\mathcal{S}}
(\widehat{M}_{jtl_1l_2}^*-M_{jtl_1l_2}^*)^2)^2)\}=O_p(\frac{q^{5}}{n^2})$, where $M_{jl_1l_2}^*$ is the $(l_1,l_2)$th element of matrix $M_{j}^*$.

Based on the first and second-order conditions for $q=o_p(n^{1/5})$ and $q/k\rightarrow \kappa$, we can obtain the consistency and oracle properties of the proposed PET estimator. Thus, the order of $\mathcal{S}_2^*$ is the largest among $\mathcal{S}_1^*,\ldots,\mathcal{S}_5^*$. Combining the above results yields $\hat{\eta}-\eta_0=\tilde{\upsilon}/\sqrt{n}+O_p(\frac{q^{5/2}}{n})$.
Using $\tilde{\upsilon}/\sqrt{n}$ to replace $\hat{\eta}-\eta_0$ in $\mathcal{S}_1^*$ and $\mathcal{S}_2^*$ yields
$\hat{\eta}-\eta_0=\tilde{\upsilon}/\sqrt{n}-\mathbb{M}^{-1}\sum_{j=1}^{\mathcal{S}}\tilde{\upsilon}_jM_j^*\tilde{\upsilon}/2n-\mathbb{M}^{-1}
\tilde{A}\tilde{\upsilon}/n+O_p(\frac{q^4}{n\sqrt{n}})$.
Replacing $\hat{\eta}-\eta_0$ in $\mathcal{S}_3^*$ and $\mathcal{S}_4^*$ by $\tilde{\upsilon}/\sqrt{n}-\mathbb{M}^{-1}\sum_{j=1}^{\mathcal{S}}\tilde{\upsilon}_jM_j^*\tilde{\upsilon}/2n-\mathbb{M}^{-1}\tilde{A}\tilde{\upsilon}/n$ and in $\mathcal{S}_1^*$ and $\mathcal{S}_2^*$ by $\tilde{\upsilon}/\sqrt{n}$ leads to Equation (\ref{high01}).

\vspace{5mm}

\noindent {\bf Proof of Corollary \ref{cor1}}. Let $\hat{\theta}_1$ be the PET estimator of nonzero parameter vector $\theta_1$. Denote $\eta=(\theta_1^{\T},\lambda^{\T})^{\T}$, $\eta_0=(\theta_{10}^{\T},0^{\T})^{\T}$, $\Gamma_{1i}(\theta_1)=\partial \psi_i(\theta_1)/\theta_1$, and
\begin{equation}{\label{coro1}}
           \begin{aligned}
m(Z_i;\eta)=\rho_1(\lambda^{\T}\psi_i(\theta_1))
\begin{pmatrix}
\Gamma_{1i}(\theta_1)\lambda \\
\psi_i(\theta_1)
\end{pmatrix}
-\begin{pmatrix}
 W(\theta_1)\\
0
\end{pmatrix},
\end{aligned}
\end{equation}
where $\rho_1(\lambda^{\T}\psi_i(\theta_1))=n\pi_i=n\exp\{\lambda^{\T}\psi(Z_{i};\theta_1)\}/\sum_{j=1}^{n}\exp\{\lambda^{\T}\psi(Z_j;\theta_1)\}$, and the components of vector $W(\theta_1)$ is $p_{\gamma}^{'}(|\theta_{1j}|){\rm sign}(\theta_{1j})$ for $j=1,\ldots,q$. Let $\rho_1=1$, $\rho_2=1-1/n$, $\rho_3=1-3/n+2/n^2$, and \begin{equation}{\label{coro2}}
           \begin{aligned}
\frac{\partial m(Z_i;\eta_0)}{\partial\eta}=
\begin{pmatrix}
-\dot{W}(\theta_{10})&\Gamma_{1i}^{\T} \\
\Gamma_{1i}&\psi_i\psi_i^{\T}
\end{pmatrix}, \mathbb{M}=\begin{pmatrix}
 -\dot{W}(\theta_{10})&\Gamma_1^{\T} \\
\Gamma_1&\Sigma_1
\end{pmatrix}, \mathbb{M}^{-1}=\begin{pmatrix}
\mathcal{K}&\mathbb{H}^{\T} \\
\mathbb{H}&P
\end{pmatrix}+O_p(1/n),
\end{aligned}
\end{equation}
where $\dot{W}(\theta_{10})_{jj}=\partial^2 p_\gamma(\theta_{10j})/\partial\theta_{1j}^2.$ It follows from Assumption \ref{ass16} that $\max_{1\leq j\leq q}\dot{W}(\theta_{10})_{jj}$ $=O_p(1/n)$. Denote  $\mathcal{K}=(\Gamma_1^{\T}\Sigma_1^{-1}\Gamma_1)^{-1}$, $\mathbb{H}=\mathcal{K}\Gamma_1^{\T}\Sigma_1^{-1}$, and $P=\Sigma_1^{-1}-\Sigma_1^{-1}\Gamma_1 \mathcal{K}\Gamma_1^{\T}\Sigma_1^{-1}$. Let $\Gamma^j_{1i}=\partial^2\psi_i(\theta_{10})/\partial\theta_{1j}\partial\theta_1^{\T}$, $\psi_i^j=\partial \psi_i(\theta_{10})/\partial\theta_{1j}$, $w'^{j}=p'''(|\theta_{10j}|)=0$, $\tau=j-q$ for $j>q$, $e_\tau$ be a $k\times 1$ vector whose $\tau$th component is 1 and 0 elsewhere. Here, $\psi_{i\tau}$ is the $\tau$th component of $\psi_i$. Then, we have
\begin{equation}{\label{coro3}}
           \begin{aligned}
\frac{\partial^2 m(Z_i;\eta_0)}{\partial\eta_j\partial\eta}=\left\{
\begin{array}{llll}
\begin{pmatrix}
0&{\Gamma_{1i}^j}^{\T} \\
\Gamma_{1i}^j&\psi_i^j\psi_i^{\T}+\psi_i{\psi_i^j}^{\T}
\end{pmatrix}   & {\rm if}~ j\leq q,\\ [5mm]
\begin{pmatrix}
\partial^2\{e_\tau^{\T}\psi_i(\theta_{10})\}/\partial\theta_1\partial\theta_1^{\T}&{\Gamma_{1i}^j}^{\T}e_\tau\psi_i^{\T}+\psi_{i\tau}{\Gamma_{1i}^j}^{\T} \\
\psi_ie_\tau^{\T}\Gamma_{1i}^j+\psi_{i\tau}\Gamma_{1i}^j&-\rho_3\psi_{i\tau}\psi_i\psi_i^{\T}
\end{pmatrix} & {\rm if}~ j> q.
\end{array}
\right.
\end{aligned}
\end{equation}
Let $\upsilon_i=-\mathbb{M}^{-1}m(Z_i;\eta_0)$, $\mathbb{V}_i=\partial m(Z_i;\eta_0)/\partial\eta-E\{\partial m(Z_i;\eta_0)/\partial\eta\}$ and $\mathcal{W}=\mathcal{W}(\theta_{10})$. It follows from Equation (\ref{coro2}) that
\[
E(\upsilon_i\upsilon_i^{\T})=
\begin{pmatrix}
\mathcal{K}\mathcal{W}^2\mathcal{K}^{\T}+\mathcal{K}  &\mathcal{K}\mathcal{W}^2\mathbb{H}\\
\mathbb{H}^{\T}\mathcal{W}^2\mathcal{K}^{\T} &\mathbb{H}^{\T}\mathcal{W}^2\mathbb{H}+P
\end{pmatrix}
=\begin{pmatrix}
 A&B^{\T} \\
B&\mathcal{C}
\end{pmatrix}
+\begin{pmatrix}
 \mathcal{K}&0 \\
0&P
\end{pmatrix},
\]
\[
E(\mathbb{V}_i\upsilon_i)=
\begin{pmatrix}
-E(\Gamma_{1i}^jP\psi_i)\\
-E(\Gamma_{1i}^j\mathbb{H}\psi_i+\psi_i\psi_i^{\T}P\psi_i)
\end{pmatrix}=
\begin{pmatrix}
\varphi\\
f
\end{pmatrix},
\]
where $A=\mathcal{K}\mathcal{W}^2\mathcal{K}^{\T}$, $B=\mathbb{H}^{\T}\mathcal{W}^2\mathcal{K}^{\T}$, $\mathcal{C}=\mathbb{H}^{\T}\mathcal{W}^2\mathbb{H}$, $\varphi=-E(\Gamma_{1i}P\psi_i)$ and $f=-E(\Gamma_{1i}^j\mathbb{H}\psi_i+\psi_i\psi_i^{\T}P\psi_i)$.

Combining the above equations yields
\[
\begin{array}{llll}
\sum\limits_{j=1}^{q+k}M_j^*E(\upsilon_i\upsilon_i^{\T})e_j/2&=&\frac{1}{2}\sum\limits_{j=1}^{q}M_j^*(A^{\T},B)^{\T}e_j+
\frac{1}{2}\sum\limits_{j=1}^{k}M_{j+q}^*(B^{\T},\mathcal{C}^{\T})^{\T}e_j\\
&&+\frac{1}{2}\sum\limits_{j=1}^{q}M_j^*(\mathcal{K},0)^{\T}e_j+
\frac{1}{2}\sum\limits_{j=1}^{k}M_{j+q}^*(0,P)^{\T}e_j\\
&=&
\begin{pmatrix}
a_1+b_1\\
c_1+d_1
\end{pmatrix}+
\begin{pmatrix}
a_2+b_2\\
c_2+d_2
\end{pmatrix}+
\begin{pmatrix}
E(\Gamma_{1i}^jP\psi_i)\\
\tilde{d}+\rho_3E(\psi_i\psi_i^{\T}P\psi_i)/2
\end{pmatrix},
\end{array}
\]
where  $a_{1}=0$, the $j$th component of $b_1$ is $b_{1j}={\rm tr}(B^{\T}E\{\partial^2\psi_i(\theta_{10})/\partial\theta_{1j}\partial\theta_1^{\T}\})/2$ for $j=1,\ldots,q$, $b_2=E\{{\Gamma_{1i}^j}^{\T}\mathcal{C}\psi_i\}$,
$a_2=\sum_{j=1}^kE\{\partial^2 \psi_{ij}/\partial\theta_1\partial\theta_1^{\T}B^{\T}e_j\}/2$, the $j$th component of $c_1$ is $c_{1j}={\rm tr}(AE\{\partial^2\psi_{ij}(\theta_{10})/\partial\theta_1\partial\theta_1^{\T}\})/2$
for $j=1,\ldots,k$, $d_1=c_2=E\{\Gamma_{1i}^jB^{\T}\psi_i\}$, $d_2=-E(\psi_i\psi_i^{\T}\mathcal{C}\psi_i)/2$,
the $j$th component of $\tilde{d}$ is $\tilde{d}_{j}={\rm tr}(\mathcal{K}E\{\partial^2\psi_{ij}(\theta_{10})/\partial\theta_1\partial\theta_1^{\T}\})/2$
for $j=1,\ldots,k$.
Then, Bias$(\hat{\theta}_1)$ is the first $p$ elements of
$\mathbb{L}_B=E\{Q_1(\tilde{\upsilon})+Q_2(\tilde{\upsilon},\tilde{A})\}/n$, which is given by
$$
           \begin{aligned}
\mathbb{L}_B&=-(n\mathbb{M})^{-1}\left(E(\mathbb{V}_i\upsilon_i)+\sum\limits_{j=1}^{q+k}M_j^*E(\upsilon_i\upsilon_i^{\T})e_j/2\right)\\
&=-(n\mathbb{M})^{-1}\left\{
\begin{pmatrix}
\mathcal{A}_1\\
\mathcal{A}_2
\end{pmatrix}- \begin{pmatrix}
0\\
\tilde{d}-E(\Gamma_{1i}^j\mathbb{H}\psi_i)+(\rho_3/2-1)E(\psi_i\psi_i^{\T}P\psi_i)
\end{pmatrix}\right\},
\end{aligned}
$$
where $\mathcal{A}_1=a_1+b_1+a_2+b_2$ and $\mathcal{A}_2=c_1+d_1+c_2+d_2$. Therefore, Bias$(\hat{\theta}_1)=\{\mathcal{K}\mathcal{A}_1+\mathbb{H}\mathcal{A}_2\}/n$+Bias$(\hat{\theta}_{1ET}).$

\newpage
\section*{References}

\begin{description}

\item Ai, C. and Chen, X. (2003). Efficient estimation of models with conditional moment restrictions containing unknown functions. {\sl Econometrica} {\bf 71} 1795-1843.

\item Bondell, H. D. and Reich, B. J. (2012). Consistent high-dimensional Bayesian variable selection via penalized credible regions. {\sl Journal of the American Statistical Association} {\bf 107} 1610-1624.

\item Bradic, J., Fan, J. and Wang, W. (2011). Penalized composite quasi-likelihood for ultrahigh-dimensional variable selection. {\sl Journal of the Royal Statistical Society} Series B {\bf 73} 325-349.

\item Caner, M. (2010). Exponential tilting with weak instruments: estimation and testing. {\sl Oxford Bulletin of Economics and Statistics} {\bf 72} 307-325.

\item Caner, M., Han, X. and Lee, Y. (2016). Adaptive elastic net GMM estimation with many invalid moment conditions: simultaneous model and moment selection. {\sl Journal of Business and Economic Statistics} (in press).

\item Caner, M. and Zhang, H. H. (2014). Adaptive elastic net GMM estimator. {\sl Journal of Business and Economics Statistics} {\bf 32} 30-47.

\item Chang, J., Chen, S. X. and Chen, X. (2015). High dimensional generalized empirical likelihood for moment restrictions with dependent data. {\sl Journal of Econometrics} {\bf 185} 283-304.

\item Chen, X. and Pouzo, D. (2012). Estimation of nonparametric conditional moment models with possibly nonsmooth generalized residuals. {\sl Econometrica} {\bf 80} 277-321.

\item Darolles, S., Fan, Y., Florens, J. and Renault, E. (2011). Nonparametric instrumental regression. {\sl Econometrica} {\bf 79} 1541-1566.


\item Fan, J. and Li, R. (2001). Variable selection via nonconcave penalized likelihood and its oracle properties. {\sl Journal of the American Statistical Association} {\bf 96} 1348-1360.


\item Fan, J. and Lv, J. (2011). Nonconcave penalized likelihood with NP-dimensionality. {\sl  IEEE Transactions on Information Theory} {\bf 57} 5467-5484.

\item Fan, J. and Peng, H. (2004). Nonconcave penalized likelihood with a diverging number of parameters. {\sl Annals of Statistics} {\bf 32} 928-961.

\item Hall, P. and Horowitz, J. (2005). Nonparametric methods for inference in the presence of instrumental variables. {\sl Annals of Statistics} {\bf 33} 2904-2929.

\item Harrion, D. and Rubinfeld, D. L. (1978). Hedonic prices and the demand for clean air. {\sl Journal of Environmental Economics and Management} {\bf 5} 81-102.

\item Imbens, G. W., Spady, R. H. and Johnson, P. (1998). Information theoretic approaches to inference in moment condition models. {\sl Econometrica} {\bf 66} 333-357.

\item Kitamura, Y. (2000). Comparing misspecified dynamic econometric models using nonparametric likelihood. Department of Economics, University of Wisconsin.

\item  Kosorok, M. R. (2008). {\sl Introduction to Empirical Processes and Semiparametric Inference.} Springer, New York.

\item Lam, C. and Fan, J. (2008). Profile-kernel likelihood inference with diverging number of parameters. {\sl Annals of Statistics} {\bf 36} 2232-2260.

\item Leng, C. and Tang, C. Y. (2012). Penalized empirical likelihood and growing dimensional general estimating equations. {\sl Biometrika} {\bf 99} 703-716.

\item Li, G., Peng, H. and Zhu, L. (2011). Nonconcave penalized M-estimation with a diverging number of parameters. {\sl Statistica Sinica} {\bf 21} 391-419.


\item Lu, W., Goldberg, Y. and Fine, J. P. (2012). On the robustness of the adaptive lasso to model misspecification. {\sl Biometrika} {\bf 99} 717-731.

\item Lv, J. and Fan, J. (2009). A unified approach to model selection and sparse recovery using regularized least squares. {\sl Annals of Statistics} {\bf 37} 3498-3528.

\item Newey, W. and McFadden, D. (1994). Large sample estimation and hypothesis testing.
In {\sl Handbook of Econometrics 4} (R. F. Engle and D. L. McFadden, eds.) 2111-2245. North-Holland, Amsterdam.


\item Owen, A. B. (2001). {\sl Empirical Likelihood}. Chapman and Hall, New York.

\item Qin, J. and Lawless, J. (2008). Estimating equations, empirical likelihood and constraints on parameters. {\sl Canadian Journal of Statistics} {\bf 23} 145-159.

\item Schennach, S. M. (2005). Bayesian exponentially tilted empirical likelihood. {\sl Biometrika} {\bf 92} 31-46.

\item Schennach, S. M. (2007). Point estimation with exponentially tilted empirical likelihood.  {\sl Annals of Statistics} {\bf 35} 634-672.

\item Shi, Z. (2015). Econometric estimation with high-dimensional moment equalities. Available at {\sl SSRN}: http://ssrn.com/abstract=2491102.

\item Staiger, D. and Stock, J. (1997). Instrumental variables regression with weak instruments. {\sl Econometrica} {\bf 65} 557-586.

\item Tang, C. Y. and Leng, C. (2010). Penalized high dimensional empirical likelihood. {\sl Biometrika} {\bf 97} 905-920.


\item Wang, H., Li, B. and Leng, C. (2009). Shrinkage tuning parameter selection with a diverging number of parameters. {\sl Journal of the Royal Statistical Society} Series B {\bf 71} 671-683.

 \item Zhang, C. H. (2010). Nearly unbiased variable selection under minimax concave penalty. {\sl The Annals of statistics} {\bf 38} 894-942.

\item Zhu, H., Zhou, H., Chen, J., Li, Y., Lieberman, J. and Styner, M. (2009). Adjusted exponentially tilted likelihood with applications to brain morphology.
{\sl Biometrics} {\bf 65} 919-927.

\item Zou, H. and Zhang, H. (2009). On the adaptive elastic-net with a diverging number of parameters. {\sl Annals of Statistics} {\bf 37} 1733-1751.

\item Zou, H. (2006). The adaptive lasso and its oracle properties. {\sl Journal of the American Statistical Association} {\bf 101} 1418-1429.

\end{description}


\newpage
{\renewcommand{\arraystretch}{0.6} \tabcolsep0.04 in
\begin{center}
{Table 1. Performance of the PET likelihood in Experiment 1 for Model C}\\[2mm]
 \begin{tabular}{ccccccccccccc}\hline \hline
&&\multicolumn{5}{c}{$d_{jl}=0.3$} &&\multicolumn{5}{c}{$d_{jl}=0.7$}\\
\cline{3-7} \cline{9-13}
&& \multicolumn{3}{c}{RMS} &  && \multicolumn{3}{c}{RMS} & \\
\cline{3-5} \cline{9-11}
$(n,p)$ &Method&$\hat\theta_1$&$\hat\theta_2$& $\hat\theta_3$&T&F &&$\hat\theta_1$&$\hat\theta_2$&$\hat\theta_3$&T&F\\ \hline
$(50,7)$&Mean&0.082 &	0.093& 0.091 &&&&0.096&0.098&0.101&& 	\\
        &PET &0.009 &	0.009& 0.009 &3.89&0.08&&0.007&0.007&0.007&3.98&0.00\\
        &HT  &0.083 &	0.093& 0.081 &3.05&0.74&&0.077&0.089&0.097&2.87&0.42\\
        &ST  &0.198 &	0.176& 0.171 &3.29&0.33&&0.156&0.139&0.148&3.55&0.23\\
        &QL  &0.079 &	0.085& 0.090 &2.90&0.28&&0.079&0.088&0.080&0.98&0.16\\
$(100,10)$&Mean&0.023&	0.037& 0.018 &&&&0.025&0.039&0.022&&\\	
          &PET &0.002&	0.001& 0.002 &6.89&0.05&&0.001&0.000&0.001&6.99&0.00\\
          &HT  &0.031&	0.028& 0.021 &6.21&0.56&&0.013&0.026&0.021&6.98&0.24\\
          &ST  &0.089&	0.078& 0.080 &2.81&0.22&&0.043&0.031&0.044&4.67&0.09\\
          &QL  &0.036&	0.019& 0.016 &4.45&0.17&&0.011&0.009&0.016&5.11&0.08\\
$(200,14)$&Mean&0.010&	0.008& 0.009 &&&&0.010&0.010&0.009&&	\\
          &PET &0.000&	0.000& 0.000 &11.02&0.01&&0.000&0.000&0.001&11.06&0.00\\
          &HT  &0.008&	0.007& 0.008 &10.23&0.09&&0.002&0.004&0.003&10.92&0.02\\
          &ST  &0.014&	0.018& 0.010 &7.65&0.01&&0.011&0.009&0.014&8.17&0.00\\
          &QL  &0.006&	0.009& 0.008 &8.88&0.02&&0.003&0.005&0.006&9.10&0.00\\
$(500,19)$&Mean&0.004&	0.002& 0.003 &&&&0.003&0.005&0.001&&	\\
          &PET &0.000&	0.000& 0.000 &16.00&0.00&&0.000&0.000&0.000&16.00&0.00\\
          &HT  &0.002&	0.003& 0.002 &15.13&0.07&&0.001&0.001&0.001&15.32&0.00\\
          &ST  &0.006&	0.008& 0.006 &11.35&0.00&&0.004&0.002&0.007&12.13&0.00\\
          &QL  &0.001&	0.003& 0.000 &9.81&0.00&&0.001&0.000&0.001&11.80&0.00\\
\hline
\end{tabular} \\ \vspace{3mm}
\footnotesize{Note: `T' represents the average number of correctly estimated zero coefficients, `F' denotes the average number of incorrectly estimated zero coefficients.}
\end{center}
}

{\renewcommand{\arraystretch}{0.6} \tabcolsep0.04 in
\begin{center}
{Table 2. Performance of the SCAD-PET under different Criteria of tuning parameter selection.}\\[2mm]
\begin{tabular}{cccccccccccccccc}\hline \hline
 &&\multicolumn{3}{c}{aBIC} && \multicolumn{3}{c}{BIC} &&\multicolumn{3}{c}{AIC}\\
 \cline{3-5} \cline{7-8} \cline{10-11}
&(n,p) &&MS & CM &&   MC & CM && MS& CM\\
\cline{2-13}
 &(50,7) &&2.4&$75\%$&&2.2&$70\%$&&5.2&$82\%$\\
 &(100,10) &&2.7&$86\%$&&2.5&$81\%$&&6.8&$78\%$\\
 &(200,14) &&2.9&$94\%$&&2.6&$92\%$&&10.3&$70\%$\\
 &(500,19) &&3&$100\%$&&2.9&$95\%$&&12.0&$51\%$\\
 \hline
\end{tabular}\\ \vspace{3mm}
\end{center}
}

\newpage
{\renewcommand{\arraystretch}{0.7} \tabcolsep0.08 in
\begin{center}
{Table 3. Frequency ($\%$) that the true value of $\theta_2$ does not fall in the $95\%$ PET-likelihood-ratio-based confidence interval in Experiment 1 for Model C}\\[2mm]
\begin{tabular}{cccccccccccccccccccccc}\hline \hline
    &     &\multicolumn{5}{c}{$d_{jl}=0.3$} && \multicolumn{5}{c}{$d_{jl}=0.7$}\\ \cline{3-7} \cline{9-13}
    &     &\multicolumn{5}{c}{True value of $\theta_2$} && \multicolumn{5}{c}{True value of $\theta_2$}\\ \cline{3-7} \cline{9-13}
$n$ & $p$ & 0.4 & 0.5 & 0.6& 0.7 & 0.8  && 0.4 & 0.5 & 0.6 & 0.7  & 0.8\\ \hline
$50$& $7$ & 18.2& 10.2& 6.4& 8.8 & 23.3 && 37.6 & 13.2 & 6.2 & 11.8 & 39.5\\
$100$&$10$& 36.8& 17.6& 6.0& 15.6 & 39.0&& 58.2  & 32.3 & 5.8 & 38.1& 63.2\\
$200$&$14$&79.2& 37.1& 6.3& 41.1 & 69.9 && 92.2  & 53.2 & 5.7 & 62.2& 89.2\\
$500$&$19$&98.8&71.2&5.4 & 79.4 & 97.3  && 100.0 & 87.8 & 5.2 & 81.6 & 100.0\\ \hline
\end{tabular} \\ \vspace{3mm}
\end{center}
}

\vspace{5mm}
{\renewcommand{\arraystretch}{0.7} \tabcolsep0.04 in
\begin{center}
{Table 4. Performance of the PET and PEL estimates in Experiment 1 for Model M}\\[2mm]
\begin{tabular}{cccccccccccccccc}\hline \hline
&\multicolumn{7}{c}{PET} && \multicolumn{7}{c}{PEL}\\ \cline{2-8} \cline{10-16}
&\multicolumn{3}{c}{$d_{jl}=0.3$} &&\multicolumn{3}{c}{$d_{jl}=0.7$} && \multicolumn{3}{c}{$d_{jl}=0.3$} &&\multicolumn{3}{c}{$d_{jl}=0.7$}\\ \cline{2-4} \cline{6-8} \cline{10-12}\cline{14-16}
& Bias & RMS & SD && Bias & RMS & SD && Bias & RMS & SD && Bias & RMS & SD\\ \cline{2-8} \cline{10-16}
$\hat\theta_1$ &0.022&0.079&0.082&&0.018&0.076&0.074&&0.154&0.240&0.142&&0.145&0.213&0.123\\
$\hat\theta_2$ &0.037&0.088&0.092&&0.023&0.079&0.080&&0.168&0.229&0.139&&0.137&0.207&0.127\\
$\hat\theta_3$ &0.023&0.083&0.082&&0.019&0.073&0.076&&0.152&0.232&0.133&&0.146&0.209&0.119\\
$\hat\theta_4$ &0.002&0.001&0.001&&0.001&0.001&0.001&&0.137&0.197&0.137&&0.126&0.199&0.125\\
$\hat\theta_5$ &0.001&0.002&0.001&&0.001&0.001&0.001&&0.145&0.237&0.136&&0.127&0.225&0.125\\
$\hat\theta_6$ &0.001&0.001&0.001&&0.000&0.001&0.001&&0.143&0.196&0.136&&0.126&0.206&0.126\\
$\hat\theta_7$ &0.001&0.001&0.002&&0.000&0.001&0.001&&0.138&0.234&0.135&&0.119&0.184&0.124\\
T              &\multicolumn{3}{c}{3.42}&&\multicolumn{3}{c}{3.69}&&\multicolumn{3}{c}{1.53}&&\multicolumn{3}{c}{1.72}\\
F              &\multicolumn{3}{c}{0.11}&&\multicolumn{3}{c}{0.08}&&\multicolumn{3}{c}{1.11}&&\multicolumn{3}{c}{0.93}\\
\hline
\end{tabular}
\end{center}
}

\vspace{5mm}
{\renewcommand{\arraystretch}{0.6} \tabcolsep0.1 in
\begin{center}
{Table 5. Performance of the PET and least squares estimates in Experiment 2}\\[2mm]
 \begin{tabular}{cccccccccccccccccccccc}\hline
 &&\multicolumn{5}{c}{PET}     && \multicolumn{3}{c}{Least squares method}\\ \cline{3-7} \cline{9-11}
 && \multicolumn{3}{c}{RMS} && && \multicolumn{3}{c}{RMS}\\ \cline{3-5} \cline{9-11}
$n$&$p$&$\theta_1$&$\theta_2$&$\theta_5$&T&F   &&$\theta_1$&$\theta_2$&$\theta_5$\\ \hline
 50&7  &0.138         &	0.103        & 	0.121       &3.83&0&&0.153           &	0.112         &	0.156          \\
100&10 &0.112         &	0.092        &	0.101       &6.78&0&&0.133           &	0.101         &	0.121          \\
200&14 &0.018         &	0.019 	     &  0.018       &10.88&0&&0.021          &	0.028	      & 0.023          \\
500&19 &0.007         &	0.006 	     &  0.006       &15.85&0&&0.009          &	0.008	      & 0.008          \\ \hline
\end{tabular} \\
\end{center}
}
\vspace{4mm}

\newpage
{\renewcommand{\arraystretch}{0.6} \tabcolsep0.04 in
\begin{center}
{Table 6. Performance of the PET and penalized empirical likelihoods in Experiment 3}\\[2mm]
 \begin{tabular}{ccccccccccccc}\hline \hline
&&\multicolumn{5}{c}{PET} &&\multicolumn{5}{c}{PEL}\\
\cline{3-7} \cline{9-13}
$(n,q^2-q)$ &par.& Bias& SD & RMS& T&F  &&  Bias& SD & RMS& T&F \\ \hline
$(185,6)$&         &  & & &   2.89 &  0.24      &&   &   &   &2.77&0.31\\
        &$\varphi_{12}$&0.086 &	0.310& 0.299 & & &&-0.092&0.321&0.302&  &  \\
       &$\varphi_{23}$&-0.132 &	1.212& 1.228 & & &&-0.146&1.315&1.330&  &  \\
$(392,12)$&         &  & & &   8.19 &  0.11      &&   &   &         &8.22&0.14	\\
        &$\varphi_{12}$&-0.055 &0.16& 0.18& & &&-0.057&0.18&0.21&  &  \\
       &$\varphi_{23}$&0.089 &	0.82& 0.84 & & &&-0.083&0.79&0.77&  &  \\
       &$\varphi_{34}$&-0.102 &	0.78& 0.74 & & &&0.110&0.89&0.86&  &  \\

$(919,20)$&         &  & & &   15.72 &  0.00      &&   &   &   &15.81&0.02	\\
        &$\varphi_{12}$&0.016 &	0.090& 0.093 & & &&-0.011&0.089&0.087&  &  \\
       &$\varphi_{23}$&-0.040 &	0.198& 0.199 & & &&0.032&0.199&0.197&  &  \\
        &$\varphi_{34}$&0.041 &	0.277& 0.280 & & &&0.047&0.287&0.285&  &  \\
         &$\varphi_{45}$&-0.011 &	0.106& 0.109 & & &&0.009&0.096&0.097&  &  \\ \hline
$(185,6)$ &$b_{21}$&0.052 &	0.299& 0.302 & & &&0.061&0.304&0.314&  &  \\
          &$b_{42}$&-0.057 &0.281& 0.278 & & &&-0.059&0.292&0.287&  &  \\
          &$b_{63}$&0.029&	0.292& 0.274 & & &&0.032&0.295&0.279&  &  \\
          &$\phi_1$&0.045 &	0.310& 0.321 & & &&0.052&0.333&0.341&  &  \\
          &$\phi_2$&-0.032 &	0.298& 0.292 & & &&-0.036&0.312&0.309&  &  \\
          &$\phi_3$&0.076 &	0.331& 0.329 & & &&0.081&0.389&0.392&  &  \\
          &$\phi_4$&0.051 &	0.289& 0.285 & & &&0.057&0.288&0.290&  &  \\
          &$\phi_5$&-0.031 &	0.189& 0.187 & & &&-0.031&0.188&0.189&  &  \\
          &$\phi_6$&-0.041 &	0.213& 0.216 & & &&-0.047&0.230&0.226&  &  \\
          &$\tau_1$&0.057 &	0.312& 0.316 & & &&0.059&0.320&0.324&  &  \\
          &$\tau_2$&0.061 &	0.381& 0.375 & & &&0.063&0.385&0.382&  &  \\
          &$\tau_3$&-0.031 &	0.236& 0.233 & & &&-0.032&0.235&0.237&  &  \\
\hline
\end{tabular} \\ \vspace{3mm}
\footnotesize{Note: `T' represents the average number of correctly estimated zero coefficients, `F' denotes the average number of incorrectly estimated zero coefficients.}
\end{center}
}

 \newpage
 {\renewcommand{\arraystretch}{0.6} \tabcolsep0.1 in
\begin{center}
{Table 7. Estimates (Est), standard errors (SE), 95\% confidence intervals (CI) of nonzero parameters in the Boston Housing data.}\\[2mm]
\begin{tabular}{ccccccccccccccccc}\hline
&\multicolumn{3}{c}{PET} && \multicolumn{3}{c}{PEL}  \\ \cline{2-4} \cline{6-8}
Variable  & Est   & SE & CI & &Est & SE & CI \\ \hline
$x_{1}$   &-8.51 & 0.318   &  (-9.13,-7.88)  &&  -7.82   & 1.354   & (-10.47,-5.16)  \\
$x_{3}$   &-1.56 & 0.067   & (-1.69,-1.43)   && -1.50    & 0.128   &(-1.75,-1.24)    \\
$x_{4}$   &0.829  & 0.061   & (0.709,0.949)  &&  0.293   & 0.114   & (0.041,0.490)  \\
$x_{6}$ &	2.040 & 0.022   &  (1.997,2.084) &&  2.13   & 0.045   & (2.038,2.216)   \\
$x_{7} $  &	2.766 &  0.057  &  (2.653,2.879) &&  3.28   & 0.111   & (3.061,3.496)   \\
$x_{9}$   &	1.191 &  0.075  & (1.042,1.340)  && 1.48    &  0.228  & (1.031,1.930)   \\
$x_{11}$   &0.567 &  0.026   & (0.514,0.619) &&  0.673   & 0.079   & (0.540,0.804) \\
$x_{12}$   &0.384 &  0.020  & (0.344,0.424)  && 0.642    & 0.026   &  (0.589,0.693) \\
$x_{13}$   &2.076 &  0.047   &  (1.983,2.170)&& 2.135    & 0.074   & (1.989,2.279) \\
$x_{1}x_{3}$&-5.89&  0.246  &  (-6.37,-5.40) &&  -6.51   & 0.658  & (-7.81,-5.22)   \\
$x_{1}x_{4}$&0.222&  0.028  & (0.166,0.278)  &&  0.165  &  0.017  & (0.131,0.199)  \\
$x_{1}x_{5}$&-1.01& 0.052  & (-1.11,-0.91)   &&  -1.13   &0.071   & (-1.26,-1.98)   \\
$x_{1}x_{6}$ &0.423& 0.025   & (0.372,0.474) &&  0.767  & 0.073   &  (0.623,0.910)     \\
$x_{1}x_{9}$&-15.8& 0.451   &  (-16.67,-14.90)&& -15.6    & 1.238   & (-17.98,-13.12) \\
$x_{1}x_{10}$&29.46& 0.894   & (27.70,31.22)  &&  27.28   & 2.262   & (22.84,31.72)   \\
$x_{3}x_{5}$ &1.006& 0.039   &   (0.929,1.083)&&   1.348  &   0.049 & (1.251,1.445)  \\
$x_{3}x_{6}$ &1.193& 0.049   & (1.095,1.290)  &&  1.024   &  0.103  & (0.821,1.227)   \\
$x_{3}x_{11}$ &-0.34& 0.028   & (-0.40,-0.28) &&   -0.51  &  0.062  &  (-0.63,-0.39) \\
$x_{4}x_{5}$ &-0.72&  0.062  &  (-0.84,-0.59) &&   -0.45  &  0.029  &  (-0.51,-0.39)   \\
$x_{4}x_{6}$  &-0.19& 0.023   & (-0.24,-0.15) && -0.28    & 0.026   &  (-0.33,-0.22) \\
$x_{5}x_{7}$  &-0.95&  0.034  & (-1.02,-1.88) && -1.23    & 0.049   &  (-1.33,-1.13)  \\
$x_{6}x_{7}$  &-0.80&  0.035  &  (-0.87,-0.73)&& -0.91    & 0.052   & (-1.08,-0.88)   \\
$x_{6}x_{10}$  &-1.63& 0.050   & (-1.73,-1.53)&& -1.60    &0.090    &(-1.88,-1.46)   \\
$x_{6}x_{11}$  &-0.95& 0.026   & (-0.99,-0.89)&& -1.04    &0.061    &(-1.10,-0.80)    \\
$x_{6}x_{13}$  &-1.42& 0.039   & (-1.49,-1.34)&& -1.43    & 0.061   &  (-1.62,-1.34)  \\
$x_{7}x_{9}$  &0.320& 0.017   & (0.281,0.352) && 0.778    & 0.066   & (0.572,0.823)    \\
$x_{7}x_{11}$  &-0.70& 0.028   & (-0.75,-0.64)&& -0.81    & 0.069   & (-1.01,-0.76)   \\
$x_{7}x_{12}$  &-0.62& 0.029   & (-0.67,-0.56)&& -0.84    & 0.037   &  (-0.99,-0.76)  \\
$x_{7}x_{13}$  &-0.20& 0.019   & (-0.23,-0.16)&&  -0.41   & 0.024   &  (-0.53,-0.32)  \\
$x_{9}x_{11}$  &-1.28& 0.082   & (-1.44,-1.12)&&  -1.93   & 0.200   &  (-2.27,-1.61)  \\
$x_{10}x_{11}$  &2.27& 0.053   & (2.162,2.373)&&  2.647   &  0.119  & (2.433,2.899)  \\
$x_{10}x_{13}$  &-0.92& 0.032  & (-0.98,-0.85)&& -0.95    & 0.091  & (-1.11,-0.82)     \\ \hline
\end{tabular} \\ \vspace{3mm}
\end{center}
}

\end{document}